\documentclass[3p,11pt]{elsarticle}

\usepackage{amssymb}
\usepackage{amsmath}
\usepackage{amsfonts}
\usepackage{indentfirst}
\usepackage{graphicx}
\usepackage{xcolor}
\usepackage[english,french]{babel}

\usepackage{tikz}
\usepackage{pgfplots}

\numberwithin{equation}{section}

\newtheorem{tm}{Theorem}
\newtheorem{lm}{Lemma}
\newtheorem{df}{Definition}
\newtheorem{prop}{Proposition}
\newtheorem{cor}{Corollary}

\newtheorem{rem}{Remark}
\newcommand {\Om} {\Omega}

\newcommand   {\ep}  {\varepsilon}
\newcommand   {\vps}{\varepsilon}

\newcommand {\noame} {\noalign{\medskip}}
\newcommand {\dis} {\displaystyle}

\newcommand{\gb}{\mathbf{g}}
\newcommand{\nb}{\mathbf{n}}
\newcommand{\ub}{\mathbf{u}}

\newcommand{\ubeps}{\ub_{\varepsilon}}
\newcommand{\vb}{\mathbf{v}}
\newcommand{\vbeps}{\vb_{\varepsilon}}
\newcommand{\wb}{\mathbf{w}}

\newcommand{\wbeps}{\mathbf{w}_\vps}
\newcommand{\Jb}{\mathbf{J}}
\newcommand{\Vb}{\mathbf{V}}

\newcommand{\dv}{\mbox{div}\,}

\newcommand{\G}{\Gamma}
\newcommand{\Geps}{\G_{\varepsilon}}

\newcommand{\Omeps}{\Omega_\ep}

\newcommand{\Veps}{\Vb_{\ep}}

\newcommand{\Thalf}{T_{\textrm{half}}}

\def \ZZ {\mathbb Z}
\def \RR {\mathbb R}

\newcommand{\Ac}{\mathcal{A}}
\newcommand{\Lc}{\mathcal{L}}

\def \beq {\begin{equation}}
\def \eeq {\end{equation}}
\def \ba {\begin{array}}
\def \ea {\end{array}}

\def \noame{\noalign{\medskip}}

\def \ep {\varepsilon}
\def \eps {\varepsilon}
\def \om {\omega}
\def \Om {\Omega}

\begin{document}

\begin{frontmatter}


\selectlanguage{english}
\title{A generalized Reynolds equation for micropolar flows past a ribbed surface with nonzero boundary conditions}

\selectlanguage{english}
\author[authorlabel1]{M. Bonnivard}
\ead{mathieu.bonnivard@u-paris.fr}
\author[authorlabel2]{I. Pa\v{z}anin}
\ead{pazanin@math.hr}
\author[authorlabel3]{F.J. Su\'arez-Grau}
\ead{fjsgrau@us.es}

\address[authorlabel1]{Universit\'e de Paris and Sorbonne Universit\'e, CNRS, Laboratoire J-L. Lions / LJLL, , F-75006 Paris, France}
\address[authorlabel2]{Department of Mathematics, Faculty of Science, University of Zagreb,
Bijeni\v{c}ka 30, 10000 Zagreb, Croatia}
\address[authorlabel3]{Departamento de Ecuaciones Diferenciales y An\'alisis Num\'erico, Facultad de Matem\'aticas, Universidad de Sevilla, \\ C/ Tarfia s/n,  41012 Sevilla, Spain}

\medskip

\begin{abstract}
Inspired by the lubrication framework, in this paper we consider a micropolar fluid flow through a rough thin domain, whose thickness is considered as the small parameter $\ep$ while the roughness at the bottom is defined by a periodical function with period of order $\ep^{\ell}$ and amplitude $\ep^{\delta}$, with $\delta>\ell>1$. Assuming nonzero boundary conditions on the rough bottom and by means of a version of the unfolding method, we identify a critical case $\delta={3\over 2}\ell-{1\over 2}$ and obtain three macroscopic models coupling the effects of the rough bottom and the nonzero boundary conditions. In every case we provide the corresponding micropolar Reynolds equation. We apply these results to carry out a numerical study of a model of squeeze-film bearing lubricated with a micropolar fluid. Our simulations reveal the impact of the roughness coupled with the nonzero boundary conditions on the performance of the bearing, and suggest that the introduction of a rough geometry may contribute to enhancing the mechanical properties of the device.
\vskip 0.5\baselineskip

}
\end{abstract}
\end{frontmatter}

\selectlanguage{english}

\section{Introduction}
Microfluidics is a multidisciplinary field intersecting engineering, physics, chemistry, microtechnology and biotechnology, with practical applications to the design of systems in which such small volumes of fluids will be used. Microfluidic area emerged in the beginning of the 1980s and is used in the development of inkjet printheads, DNA chips, lab-on-a-chip technology, micro-propulsion, and micro-thermal technologies.

Microfluidics deals with the manipulation of lubricants that are  geometrically constrained to a small (typically sub-millimetre) scale, and with the experimental and theoretical study of their mechanical behaviour. 
This behaviour can differ from ``macrofluidic" behaviour since, at the microscale, factors such as surface tension, energy dissipation, and fluidic resistance start to dominate the system. In particular, when a lubricant is in contact with a solid, at a small scale, surfacic effects may become preponderant. As a result, in order to reduce the energy dissipation of microscaled fluid-solid systems, one needs to understand and quantify very precisely the behaviour of the fluid near a solid wall. 

 From an experimental point of view, an efficient method to reduce the friction consists
in using a certain type of rough boundaries, that are called riblets.
These riblets are characterized by fast oscillations in the transversal direction, with a low
amplitude, and by their constancy in the  direction of the flow; they are essentially one-dimensional perturbations of the boundary of the solid.  The aim of the use of riblets is to prevent vortices to appear in
the neighborhood of the solid wall, and thus to reduce the momentum transfer from the
vortices to the solid boundary. By using homogenization techniques, the influence of riblets on the slip behaviour of viscous fluids  has been studied recently. In \cite{Bucur1}, starting with perfect slip condition at a highly ribbed surface, it is showed that when the oscillating parameter goes to zero, no-slip condition appears in the transversal direction while perfect slip still holds in the  direction of the flow. This means that  riblets tend to prevent
the fluid from slipping laterally, whereas the motion in the  direction of the flow is allowed with no constraint. In the same spirit,  it was proved in \cite{CLS-SIMA} that surfaces with low amplitude riblets give rise to a friction parameter in the transversal direction and no roughness effects in the direction of the flow.

The mathematical models for describing the motion of the lubricant in a device with small volume usually result from the simplification of the
geometry of the lubricant film, i.e.\! its thickness. Using the  film thickness as a small parameter $\varepsilon$, a
simple asymptotic approximation can be easily derived providing a well-known Reynolds equation
for the pressure of the fluid. Formal derivation goes back to the 19th century and the celebrated
work of Reynolds \cite{Reynolds}. The justification of this approximation, namely the proof that it can be obtained as the limit of the Stokes system (as thickness tends to zero) is provided in \cite{Bayada1} for a Newtonian flow between two plain surfaces. Different Reynolds equations for Newtonian fluids including roughness effects have been obtained  for example in
 \cite{Bayada2,Benh, BonSG, Bresch,CLS-CRM, CLS-SIMA,Chupin, PazSua-CRMec, PazSua_BMMS}.    \\


Nevertheless, most of the modern lubricants are no longer
Newtonian fluids, since the use of additives
in lubricants has become a common practice in
order to improve their performance.
Therefore, several microcontinuum theories \cite{Eringen1} have
been proposed to account for the effects of additives.
Eringen micropolar fluid theory \cite{Eringen} ignores the deformation of the microelements
and allows for the particle micromotion to take
place.
 From a mathematical point of view, a micropolar Reynolds equation was obtained in \cite{Bayada3}  for a micropolar flow in a thin film with a plain bottom assuming zero boundary conditions for microrotation. Other related results on the lubrication with a micropolar fluid with zero boundary condition can be found in \cite{Ja6,Sinha2}, and some others references including roughness effects  in  \cite{BonSGPaz, Boukrouche,Paz-Sua, SuaBMMS}.
 
 In the previously mentioned references, a zero boundary condition for the microrotation  is assumed, implying that the fluid elements cannot rotate on the fluid-solid interface. If $s$ is the  horizontal velocity of the boundary, these conditions are written as follows:
{\color{black} \begin{eqnarray}
{\bf u}={\bf s} & (\ub\ \hbox{velocity}),\label{veloc_bound}\\
{\bf w}={\bf 0} & (\wb\ \hbox{microrotation}). \label{zero_micro_bc}
\end{eqnarray}
}
However, more general boundary conditions for the microrotation were introduced to take into account the rotation of the microelements on the solid boundary. In the case where the boundary is flat, these conditions read 
 \begin{equation}\label{cond1}
{\alpha\over 2} (\nabla\times {\bf u})\times{{\bf n}} = {\bf w}\times{ {\bf n}}\, ,\quad {\color{black} {\bf w}\cdot {\bf n}=0}\, ,
 \end{equation}
where ${\bf n}$ is a normal unit vector to the boundary. Conditions \eqref{cond1} were effectively proved to be in good accordance with experiments, see \cite{Bessonov2, Ref14_Bessonov, Ref8_Bessonov, Ref2_Bessonov}.  The coefficient $\alpha$  describes the interaction between the given fluid and solid; it characterizes microrotation retardation on the solid surfaces.

In \cite{Bessonov2}, a generalized micropolar Reynolds equation is derived by using conditions (\ref{veloc_bound}), (\ref{cond1}), and the relevance of the new parameter $\alpha$ regarding the performance of  lubricated devices for both load and friction, is established by numerical computations. Nevertheless, it was mathematically proved in  \cite{Bayada4} that it is not possible to consider the boundary condition (\ref{cond1}) and simultaneously retain the no-slip condition (\ref{veloc_bound}) for  the velocity. This would be like considering  simultaneously, at the same boundary, a Neumann and a Dirichlet boundary condition. In order to obtain a well-posed variational formulation of the micropolar system, it is straightforward to confirm (see e.g.~\cite{Bayada4}) that a velocity condition compatible with (\ref{cond1}) needs to be introduced. This condition allows a slippage in the tangential direction and retains a non-penetration condition in the normal direction ${\bf n}$ ($\delta_0$ is a real parameter)
\begin{equation}\label{cond2}
({\bf u}-{\bf s})\times{{\bf n}}=\delta_0 (\nabla\times {\bf w})\times{{\bf n}}\, ,\quad {\bf u}\cdot {{\bf n}}=0\, .
\end{equation}
It is worth stressing that in most lubrication studies, it is assumed that the speed of the lubricant at the surface equals that of the solid surface. However, it has been found that wall slip occurs, not only in non-Newtonian flows \cite{Bair, Hat, Jacob, Meng, Spikes,Wilson}, but also in hydrodynamic lubrication or elasto-hydrodynamic lubrication \cite{Bona, Craig, Hervet, Kaneta}.  It seems that such phenomenon is linked to physical or chemical interactions of the solid surfaces with the lubricant. Several boundary conditions have been considered in those works to model the observed slippage. Most of them include limited yield stress or retain slippage value proportional to the shear stress.
In that context, condition (\ref{cond2}) appears as a new interpretation of the slippage observed in lubrication with micropolar fluids, expressed in terms of  the microrotation field ${\bf w}$.

 In \cite{Bayada4}, by using the  nonzero boundary conditions (\ref{cond1})-(\ref{cond2})  described above, in a 2-dimensional thin domain without roughness (see also \cite{Ja8} for the 3D flow), Bayada \emph{et al.}\! derive rigourously  a generalized version of the Reynolds equation taking such boundary conditions into account. They perform their study in the critical case where one  the non-Newtonian characteristic parameters of the micropolar fluid has specific (small) order of magnitude. The authors provide a comparison with the model in \cite{Bessonov2} that uses the no-slip condition \eqref{veloc_bound} for the velocity field, and observe that the introduction of slippage may enhance the performance of a bearing (that is, increase the load and reduce the friction coefficient) if the coupling number of the micropolar fluid and the nondimensional coefficient describing its slippage on the wall, are above a certain value.  
 \\
 
Observe that in previous studies, the nonzero boundary condition has been considered on a plain bottom. In this paper, we impose this condition on a surface covered by riblets with low amplitude, and use asymptotic analysis to derive a micropolar Reynolds equation coupling the effects of the nonzero boundary conditions and the riblets. Since we are interested in the effect of the roughness, we adopt a simple geometric setting where the top boundary is plane, given by $\varepsilon h$ with $h>0$. At the bottom we consider a surface covered by periodically distributed riblets with low amplitude, associated with a small parameter $\eps$, where $\ep^\delta$ is the amplitude and $\ep^\ell$ is the period, where $\delta>\ell>1$.    This type of rough surface has been treated in \cite{BonSGTierra, CLS-SIMA,SuaZAMM, SuaNA} for fluid flows with Navier slip boundary conditions.

First, we identify a range of values of the coupling parameter $N^2$, namely $ N^2 \leq 1/2 $, under which there is existence and uniqueness of solution (Theorem \ref{thm_existence}). Later, by means of homogenization and reduction of dimension techniques,  we identify the critical regime, i.e. $\delta={3\over 2}\ell-{1\over 2}$, in which the nonzero boundary conditions make appear two friction parameters reflecting the riblets effect on both the effective velocity and micropolar fields (Theorem \ref{thm_effective}). Finally, we also obtain a precise description of the corresponding Reynolds equation which implicitly contains the effective nonzero boundary conditions describing the roughness effects (see (\ref{bc_botom_case2_system_1}) for more details).  Moreover, we give the corresponding Reynolds equations corresponding to the sub-critical and super-critical regimes.  This constitutes a generalization of the results of \cite{Bayada4} to domains with rough bottom (Theorem \ref{thm_reynolds}). \\

The paper is organized as follows. In Section 2 and 3, we formulate the problem and introduce some notation. In Section 4, we state our main results providing the homogenized model and the generalized Reynolds equation, which are proved in Section 5. The details of certain explicit computations or asymptotic developments are postponed to the Appendix. Finally, in Section 6 we conduct numerical simulations based on the generalized micropolar Reynolds equation obtained for a particular lubrication device: a squeze-film bearing.

%
\section{Position of the problem}\label{Sec:PositionOfPb}
\begin{figure}[h!]
\begin{center}
\includegraphics[width=9.5cm]{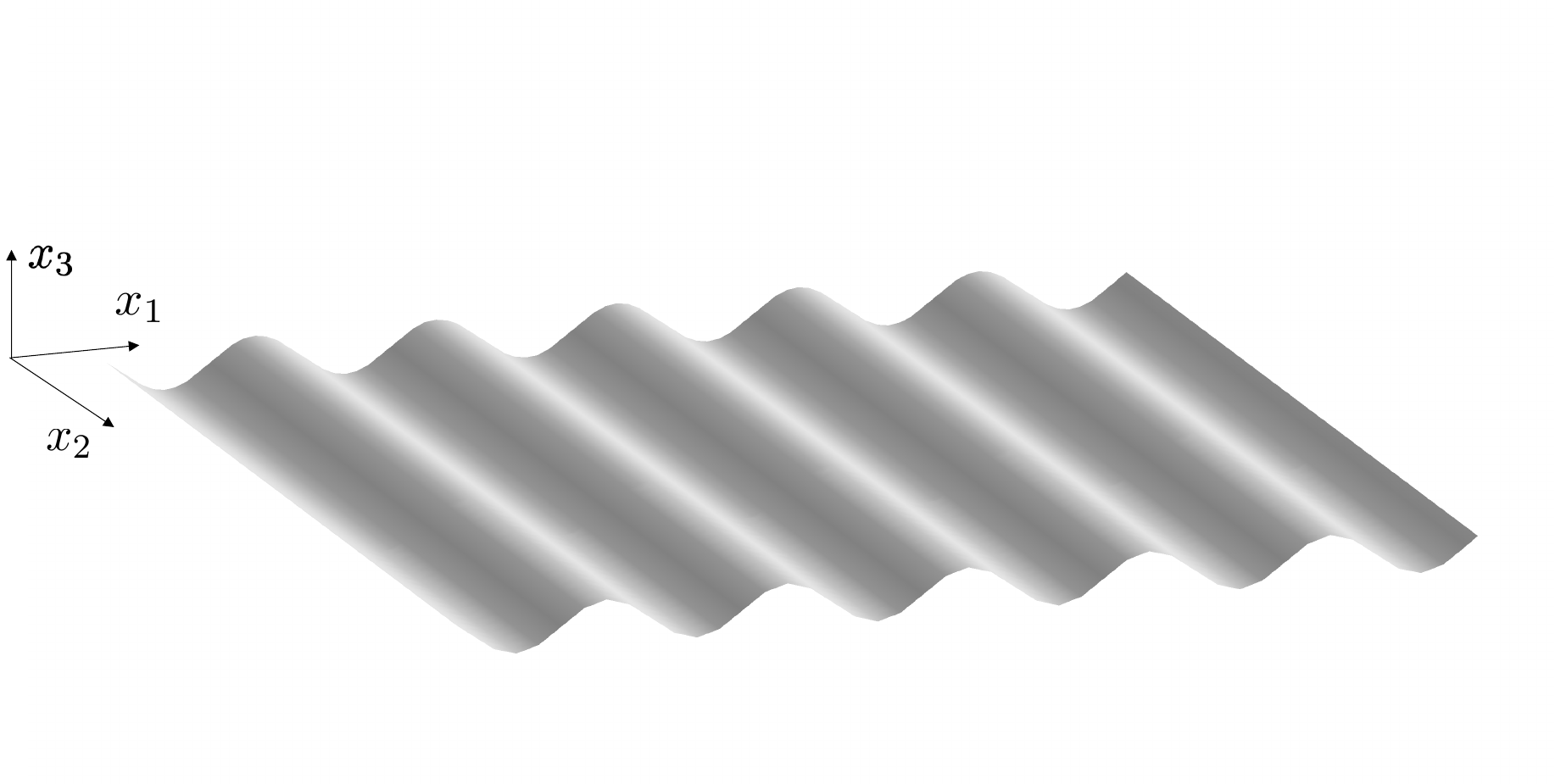}\qquad
\includegraphics[width=9cm]{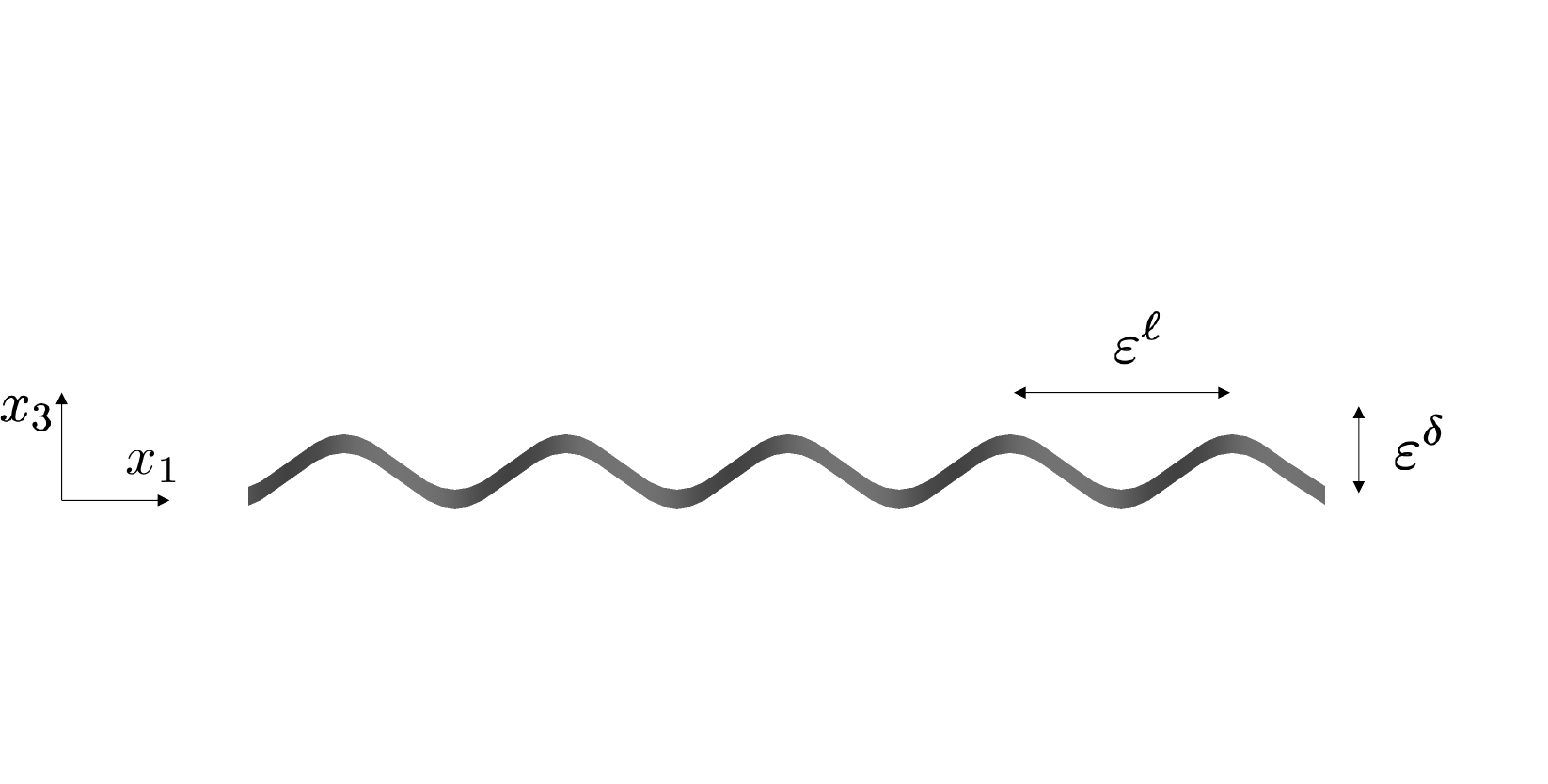}
\end{center}
\caption{Representation of the geometry of riblets (upper) and of their cross section (lower). The riblets are periodic with period $\ep^{\ell}$ in the $x_1$ direction, constant in the $x_2$ direction and oscillate with an amplitude of order $\ep^{\delta}$ in the $x_3$ direction.}
\label{fig:roughness}
\end{figure}

In the following,  $x\in \RR^3$ is decomposed as $x=(x',x_3)$ where
$x'=(x_1,x_2)\in \RR^2$ and $x_3\in \RR$. We take $e_1$, $e_2$ and $e_3$ to be the vectors of the canonical basis in $\mathbb{R}^3$, and $e_1'$, $e_2'$ to be the vectors of the canonical basis in $\mathbb{R}^2$.
The domain under consideration has the following form 
$$ \overline \Omega_\ep=\left\{(\overline x',\overline x_3)\in\mathbb{R}^2\times \mathbb{R}\ :\  \overline x'=(\overline x_1, \overline x_2)\in L\,\omega,\quad - \overline \Psi_{\varepsilon}(\overline x')<\overline x_3<h\, c\right\}.$$
Here $L$ is the characteristic length of the domain, $\omega\subset\mathbb{R}^2$ is an open subset with smooth boundary, $c$ is the characteristic distance between the plates, $h>0$ is an adimensional constant, $\ep$ is the ratio $\ep =  {c\over L}$ and
  $\overline \Psi_\ep$ is defined by
\begin{eqnarray}  \frac{1}{L} \overline \Psi_\varepsilon(\overline x') =  \lambda \ep^{\delta}\Psi\left( {1 \over L \ep^\ell}  \overline{x}'\cdot e'_1  \right)  \label{rugprofile}\end{eqnarray}
see Figure \ref{fig:roughness}, where $\lambda>0$ is an amplitude parameter and $\delta, \ell>0$ satisfy  
\begin{equation}\label{relation_parameters}
1<\ell<\delta.
\end{equation}
In definition \eqref{rugprofile}, $ \Psi\in W^{2,\infty}_\#(\RR)$ is a $\RR$-valued function with period $1$  (we use the index $\#$ to mean periodicity of period $1$), that models the roughness profile on the lower surface,  and that is normalized in the sense that 
\begin{equation}\label{Psi-normalized}
	\int_0^1 |\partial_{z_1}\Psi(z_1)|^2\, dz_1 = 1\, .
	\end{equation}
 Let $\overline\Gamma^0_\ep$, $\overline\Gamma^{1}_\ep$ and $\overline{\Gamma}^{\ell}_\ep$ denote the lower, upper and lateral boundaries on $\overline\Omega_\ep$, namely
$$\begin{array}{l}
\dis 
\overline{\Gamma}^0_\varepsilon=\left\{(\overline x',\overline x_3)\in\mathbb{R}^2\times \mathbb{R}\ :\  \overline x'\in L\, \omega,\quad \overline x_3=-\overline\Psi_\ep( \overline x')\right\},\\
\noame\dis
\overline {\Gamma}^{1}_\ep=\left\{(\overline x',\overline x_3)\in\mathbb{R}^2\times \mathbb{R}\ :\  \overline x'\in L\,\omega,\quad {\color{black} \overline  x_3=\varepsilon h L}\right\},\\
\noame\dis
\overline{\Gamma}^{\ell}_\ep=\partial\overline{\Omega}_\varepsilon-(\overline{\Gamma}^0_\varepsilon\cup \overline{\Gamma}^{1}_\ep).
\end{array}$$
The exterior normal $\overline{{\bf n}}_\ep$ to $\overline{\Gamma}^0_\varepsilon$ is defined by
\begin{equation}\label{Def:normal}
\forall \overline{x}'\in L\omega\quad \overline{{\bf n}}_\ep(\overline{x}',-\overline{\Psi}_\ep( \overline{x}')) = \frac{1}{[1+\partial_{\overline x_1}\overline{\Psi}_\ep( \overline{x}')^2]^{1/2}} \left(-\partial_{\overline x_1} \overline{\Psi}_\ep( \overline{x}'),0,-1\right)\, .
\end{equation}
For any vector field $\xi$ defined on $\overline{\Gamma}^0_\varepsilon$, we note $[\xi]_{tan}$ its tangential part, i.e.\! is the vector field defined on $\overline{\Gamma}^0_\varepsilon$ by
$[\xi]_{tan}=\xi-(\xi\cdot \overline{{\bf n}}_\ep)\overline{{\bf n}}_\ep$.

\subsection{The equations and boundary conditions}\label{subsec:equations_boundary_conditions}
The micropolar fluid  flow is described by the following equations expressing the balance of momentum, mass and angular momentum:
\begin{eqnarray}
 &&-(\nu+\nu_r)\Delta\overline{\mathbf{u}}_\ep + \nabla \overline p_\ep =2 \nu_r(\nabla\times\overline{\mathbf{w}}_{\ep})\;,\label{1,1,1}\\[0.2cm]
 &&\mbox{div}\,\overline{\mathbf{u}}_\ep =0\;,\label{21,1,1}\\[0.2cm]
 &&-{\color{black}c_r}\,\Delta\overline{\mathbf{w}}_\ep+4\nu_r\overline{\mathbf{w}}_{\ep}=2\nu_r(\nabla\times\overline{\mathbf{u}}_{\ep})\;.\label{31,1,1}
\end{eqnarray}
In the above system, velocity $\overline{\mathbf{u}}_\ep$, pressure $\overline p_\ep$ and microrotation $\overline{\mathbf{w}}_\ep$ are unknown. $\nu$ is the Newtonian 
viscosity, while $\nu_r$ and {\color{black} $c_r$} are microrotation viscosities resulting from the asymmetry of the stress tensor. All viscosity coefficients are assumed to be positive constants. 

Let {\color{black} $\overline{V}_\ep$} be the velocity of the upper plate, and {\color{black} $\overline{\mathbf{g}}_\ep$} be velocity of the fluid on the lateral boundaries of the domain. As discussed in the Introduction, the following boundary conditions are imposed
\begin{eqnarray}
&&\overline{\mathbf{u}}_\varepsilon = (0,0,-{\color{black} \overline{V}_\ep}),\quad \overline{\mathbf{w}}_\varepsilon = 0 \quad\hbox{    on    }\quad\overline \Gamma^{1}_\ep\, ,\label{bc1}\\[0.2cm]
&&\overline{\mathbf{u}}_\varepsilon ={\color{black} \overline{\mathbf{g}}_\ep},\quad \overline{\mathbf{w}}_\varepsilon = 0\quad \hbox{    on    }\quad\overline\Gamma^\ell_\ep\, ,\label{bc2_bar}\\[0.2cm]
&&\overline{\mathbf{u}}_\varepsilon\cdot  \overline{{\bf n}}_\ep = 0,\quad \overline{\mathbf{w}}_\varepsilon\cdot  \overline{{\bf n}}_\ep = 0\quad \hbox{    on     }\quad\overline \Gamma_\ep^0\, ,\label{bc3}\\[0.2cm]
&&{\alpha\over 2}  [D\overline{\ub}_\ep\,  \overline{{\bf n}}_\ep]_{tan}  = \overline{\mathbf{w}}_\varepsilon\times  \overline{{\bf n}}_\ep \quad \hbox{    on     }\quad\overline \Gamma_\ep^0\, ,\label{bc4}\\[0.2cm]
&&[D\overline{\wb}_\ep\,  \overline{{\bf n}}_\ep]_{tan}={2\nu_r\over {\color{black}c_r}}\beta\,\overline{\mathbf{u}}_\ep\times  \overline{{\bf n}}_\ep
\quad \hbox{    on     }\quad\overline \Gamma_\ep^0\, .\label{bc5}
\end{eqnarray}

{\color{black} We remark that along the paper, $D\ub$ denotes the gradient of a vectorial function $\ub=(u_i)_{1\leq i\leq 3}$, defined by $(D\ub)_{i,j}=\partial_j u_i$, and should not be confused with the symmetric part of the gradient.} Notice that the usual (Dirichlet) boundary conditions (\ref{bc1})-(\ref{bc2_bar}) for the velocity and microrotation are prescribed on $\Gamma_\ep^1\cup \Gamma_\ep^\ell$. However, on the lower part $\Gamma_\ep^0$ (corresponding to the rough boundary), new type of boundary conditions (\ref{bc4})-(\ref{bc5}) are imposed, together with the non-penetration conditions (\ref{bc3}). {\color{black} Finally, coefficient $\beta\in\RR_+$ in \eqref{bc5} is a friction coefficient that  controls the slippage of the fluid at the wall.} \\

{\color{black} Let us stress that conditions~\eqref{bc4}-\eqref{bc5} are adaptations of conditions~\eqref{cond1}-\eqref{cond2} to the present case of an oscillating boundary $\overline \Gamma_\ep^0$. Since system~\eqref{1,1,1}--\eqref{31,1,1} couples a Stokes equation on $\overline \ub_\ep$ with an elliptic system on $\overline \wb_\ep$, in the present context of slip boundary conditions, the normal conditions~\eqref{bc3} must be completed by tangential conditions on $[D\overline \ub_\ep \overline \nb_\ep]_{tan}$ and $[D\overline \wb_\ep \overline \nb_\ep]_{tan}$, in the aim of obtaining a well-posed problem. 
	
To obtain conditions~\eqref{bc4}-\eqref{bc5}, we have interpreted the rotational terms $(\nabla\times \ub)\times \nb$, $(\nabla\times \wb)\times \nb$ appearing in the initial formulatation of the tangential boundary conditions~\eqref{cond1}-\eqref{cond2}, as being respectively equal to $[D\ub\, \nb]_{tan}$ and $[D\wb\, \nb]_{tan}$. 
This is indeed the case for a flat boundary $\Gamma$ of normal ${\bf n}$, since for regular vector fields $\bf u, \bf v$ satisfying ${\bf u}\cdot {\bf n}= {\bf v}\cdot {\bf n}=0$ on $\Gamma$, there holds
\[
\int_{\Gamma}[(\nabla \times {\bf u})\times {\bf n} ]\cdot {\bf v}\,d\sigma = \int_{\Gamma} [D {\bf u}\, {\bf n}]_{tan}\cdot {\bf v}\,d\sigma\, .
\]
Last equality is obtained by writing $(\nabla\times {\bf u})\times {\bf n} = D{\bf u}\,{\bf n} - (D{\bf u})^T\,{\bf n}$ and using that  ${\bf v}\cdot \nabla({\bf u}\cdot {\bf n})=0$ on $\Gamma$, which gives
\begin{align*}
\int_{\Gamma} [(\nabla\times {\bf u})\times  {\bf n}]\cdot {\bf v} \,d\sigma& = \int_{\Gamma}[D{\bf u}\,{\bf n}]\cdot {\bf v}  - \int_{\Gamma}[(D{\bf u})^T\, {\bf n}]\cdot {\bf v}\,d\sigma\\
& = \int_{\Gamma}[D{\bf u}\,{\bf n}]_{tan}\cdot {\bf v}\,d\sigma - \int_{\Gamma} [({\bf v}\cdot\nabla){\bf u}]\cdot {\bf n}\,d\sigma\\
& = \int_{\Gamma}[D {\bf u}\, {\bf n} ]_{tan}\cdot {\bf v}\,d\sigma + \int_{\Gamma} [({\bf v}\cdot\nabla){\bf n}]\cdot {\bf u}\,d\sigma\\
& = \int_{\Gamma}[D{\bf u}\,{\bf n}]_{tan}\cdot {\bf v}\,d\sigma
\end{align*}
because $({\bf v}\cdot\nabla){\bf n}=0$ since ${\bf n}$ is a constant vector. 

Hence, conditions~\eqref{bc4}-\eqref{bc5} and~\eqref{cond1}-\eqref{cond2} are equivalent in the case of a flat boundary, so the tangential conditions~\eqref{bc4}-\eqref{bc5} can be seen as a generalization of~\eqref{cond1}-\eqref{cond2} to the case of a non flat boundary.
}\\

In \cite{Ref14_Bessonov}, it was proposed to define the parameter  $\alpha$ appearing in (\ref{bc4})  as a microrotation retardation at the boundary and to connect it with the different viscosity coefficients. It has been shown experimentally \cite{Jacob_BO, Leger} that there are chemical interactions between a solid surface and the nearest fluid layer. This had to be taken into account, especially for a non-Newtonian fluid and a very thin-film thickness. This can be done by introducing a viscosity $\nu_b$ near the surface which is different from $\nu$ and $\nu_r$. In \cite{Ref14_Bessonov}, it was proposed to define $\alpha$ by means of this boundary viscosity $\nu_b$ by
\begin{equation}\label{Def:alpha}
	\alpha={\nu+\nu_r-\nu_b\over \nu_r}.
\end{equation}
Following \cite{Ref14_Bessonov}, it is possible to give physical limits to $\nu_b$, inducing limits on $\alpha$:
\begin{equation}
0\leq \nu_b\leq \nu+\nu_r\Rightarrow 0\leq \alpha \leq {\nu+\nu_r\over \nu_r}. \label{Condition:viscosities-alpha}
\end{equation}
 The condition $\alpha=0$ is equivalent to strong adhesion of the fluid particles to the boundary surface so that they do not rotate relative to the boundary, i.e. ${\bf w}=0$. Thus, from now on, we consider $\alpha>0$ so that the stress tensor and the micro-rotation are coupled on the boundary. 
\\

It has been observed (see e.g. \cite{Bayada3,Bayada4,Ja9}) that the magnitude of the viscosity coefficients appearing in (\ref{1,1,1})-(\ref{31,1,1}) may influence the effective flow. Thus, it is reasonable to work with the system written in a non-dimensional form. In view of that, we introduce the {\color{black} characteristic velocity $V_0$ of the fluid}, and define:
\begin{equation}\label{non-dim}
\begin{array}{l}\dis
x'={\overline x' \over L},\quad x_3={\overline x_3\over L},\quad \Psi_\ep={\overline\Psi_\ep\over L}, \\
\noame\dis
{\bf u}_\ep={\overline{\bf u}_\ep\over V_0},\quad p_\ep={L\over V_0(\nu+\nu_r)}\overline p_\ep,\quad {\bf w}_\ep={L\over V_0}\overline{\bf w}_\ep,\quad
{\bf g}_\vps={{\color{black} \overline{\mathbf{g}}_\ep}\over V_0},\quad {V_\vps}={{\color{black} \overline V_\ep}\over V_0},\\
\noame\dis
N^2={\nu_r\over \nu+\nu_r},\quad R_M={{\color{black}c_r}\over \nu+\nu_r}{1\over L^2}.
\end{array}
\end{equation}
Dimensionless (non-Newtonian) parameter $N^2$ characterizes the coupling between the equations for the velocity and the microrotation, and is of order $\mathcal{O}(1)$ with respect to small parameter $\ep$.  Notice that assumption~\eqref{Condition:viscosities-alpha} yields
\begin{equation}\label{Condition:N-alpha}
{1\over\alpha} \geq N^2\, .
\end{equation}
 The second dimensionless parameter denoted by $R_M$ is related to the characteristic length of the microrotation effects and will be compared with $\ep$. {\color{black} We also assume that friction parameter $\beta$ is of order $\mathcal{O}(1)$.}

In view of the above changes of variables, the fluid domain becomes
$$  \Omega_\ep=\left\{( x', x_3)\in\mathbb{R}^2\times \mathbb{R}\ :\   x' \in \omega,\quad    -\Psi_\ep \big(x'\big)< x_3<\ep\, h\right\},$$
 where according to \eqref{rugprofile}, $\Psi_{\varepsilon}$ is given by
\begin{equation}\label{rugprofile-nondim}\Psi_\ep(x')=  \lambda\ep^{\delta}\Psi\Big({1\over \ep^\ell}x'\cdot e'_1\Big) 
\end{equation}
and the lower, upper and lateral boundaries are now described by
\begin{align}
\Gamma^0_\varepsilon &=\left\{( x', x_3)\in\mathbb{R}^2\times \mathbb{R}\ :\   x'\in   \omega,\quad x_3=-  \Psi_\ep\big(x'\big) \right\},\nonumber\\
\Gamma^{1}_\ep &=\left\{( x', x_3)\in\mathbb{R}^2\times \mathbb{R}\ :\   x'\in \omega,\quad x_3=\ep h\right\},\nonumber\\
  \Gamma^{\ell}_\ep & =\partial \Omega_\varepsilon-(\Gamma^0_\ep \cup \Gamma^{1}_\ep)\, .\nonumber
\end{align}
The exterior normal ${\bf n}_\ep$ to $\Geps^0$ is now defined by
\begin{equation*}
\forall x'\in \omega\quad {\bf n}_\ep\Big(x',- \Psi_\ep\big(x'\big)\Big) = \frac{1}{\Big[1+ \lambda^2\ep^{2(\delta-\ell)}\partial_{x_1} \Psi\left(\frac{1}{\ep^{\ell}}x'\cdot e'_1\right)^2\Big]^{1/2}} \Big(- \lambda\ep^{\delta-\ell}\partial_{x_1} \Psi\left(\frac{1}{\ep^{\ell}}x'\cdot e'_1\right),0,-1\Big)\, .
\end{equation*}
The tangential part of a vector field $\xi$ defined on $\Geps^0$ is accordingly given by
$[\xi]_{tan}=\xi-(\xi\cdot {\bf n}_\ep){\bf n}_\ep$.
%
%
The flow equations \eqref{1,1,1}-\eqref{31,1,1} now have the following form
\begin{eqnarray}
 &&-\,\Delta\mathbf{u}_\ep + \nabla p_\ep =2 N^2(\nabla\times\mathbf{w}_{\ep})\quad\hbox{in}\quad\Omega_\ep\;,\label{1,1}\\[0.2cm]
 &&\mbox{div}\,\mathbf{u}_\ep =0\quad \hbox{in}\quad\Omega_\ep,\label{21,1}\\[0.2cm]
 &&-R_M\,\Delta \mathbf{w}_\ep+4N^2\mathbf{w}_{\ep}=2N^2(\nabla\times\mathbf{u}_{\ep})\quad \hbox{in}\quad\Omega_\ep,\label{31,1}
\end{eqnarray}
with boundary conditions
\begin{eqnarray}
&&\mathbf{u}_\varepsilon = -V_\ep e_3, \quad \mathbf{w}_\varepsilon = 0 \quad \hbox{on}\quad\Gamma^1_\ep,\label{bc-1}\\[0.2cm]
&&{\mathbf{u}}_\varepsilon = {\bf g}_\vps,\quad {\mathbf{w}}_\varepsilon = 0\quad \hbox{on}\quad\Gamma^\ell_\ep,\label{bc-2}\\[0.2cm]
&&\mathbf{u}_\varepsilon\cdot {\bf n}_\varepsilon = 0,\quad \mathbf{w}_\varepsilon\cdot {\bf n}_\varepsilon = 0\quad \hbox{on}\quad\Gamma_\ep^0,\label{bc-3}\\[0.2cm]
&&  {\alpha\over 2}[D\mathbf{u}_\ep\, {\bf n}_\ep]_{tan}   = \mathbf{w}_\varepsilon\times {\bf n}_\varepsilon \quad \hbox{on}\quad\Gamma_\ep^0,\label{bc-4}\\[0.2cm]
&&R_M[D\mathbf{w}_\ep\, {\bf n}_\ep]_{tan}=2N^2\beta\,\mathbf{u}_\ep\times {\bf n}_\varepsilon\quad \hbox{on}\quad\Gamma_\ep^0.\label{bc-5}
\end{eqnarray}
The divergence-free condition \eqref{21,1} imposes the following compatibility condition on the boundary data:
\begin{equation}
\int_{\G^{\ell}_\ep}\gb_\vps\cdot {\bf n}_\ep\, d\sigma = V_\ep |\omega|\, , \label{CompatibilityCondition-geps}
\end{equation}
where $\sigma$ stands for the Hausdorff measure of dimension 2, and $|\omega|$ is the area of $\omega$.\\

 In the present paper the aim is to derive the macroscopic law describing the effective flow in $\Omega_\ep$ by using rigorous asymptotic analysis with respect to the small parameter $\ep$. In particular, we shall focus on detecting the roughness-induced effects together with the effects of nonzero boundary conditions.

Let us start by defining the notion of weak solution to system~\eqref{1,1}--\eqref{bc-5}.

\paragraph{Weak formulation of problem \eqref{1,1}--\eqref{bc-5}}

Let us introduce the functional spaces $\Veps$ and $\Veps^0$ defined by
\begin{align*}
\Vb_\vps & = \{ \varphi\in H^1(\Omeps)^3,\ \varphi_{|\G^1_{\vps}\cup\G^{\ell}_\vps}=0,\ \varphi\cdot {\bf n}_\vps = 0 \ \textrm{on}\ \G^0_\vps \},\\
\Vb_\vps^0 & = \{\varphi\in \Veps,\ \mbox{div}\,\varphi = 0\ \textrm{in}\ \Omeps\},
\end{align*}
endowed with the norm $\|D\varphi\|_{L^2(\Omeps)^3}$.  Assume that $(\ubeps, \wbeps, p_\vps)$ is a classical solution to system \eqref{1,1}--\eqref{bc-5}. Multiplying \eqref{1,1} by a test function $\varphi\in \Veps$, integrating by parts and taking into account the boundary conditions and the free divergence condition satisfied by $\varphi$, we obtain
\begin{align*}
& \int_{\Omeps} D\ub_\eps : D\varphi \,dx - \int_{\Geps^0}[D\ub_\eps \, {\bf n}_\ep]_{tan}\cdot \varphi \,d\sigma -\int_{\Omeps} p_\ep\, \dv \varphi \,dx
 - 2N^2 \int_{\Omeps}(\nabla\times \wb_\ep)\cdot \varphi  \,dx = 0\, .
 \end{align*}
Hence, boundary condition~\eqref{bc-4} yields
\begin{align*}
\int_{\Omeps} D\ub_\eps  : D\varphi \,dx    -\int_{\Omeps} p_\ep\, \dv \varphi  \,dx - 2N^2 \int_{\Omeps}(\nabla\times \wb_\ep)\cdot \varphi \,dx  - {2\over\alpha} \int_{\Geps^0}(\wb_\ep\times {\bf n}_\ep)\cdot \varphi  \,d\sigma = 0 \, .
\end{align*}
Using the integration by part formula
\begin{equation}\label{bypartsform}
\int_{\Omeps}(\nabla\times \wbeps)\cdot\varphi  \,dx= \int_{\Omeps}(\nabla\times \varphi)\cdot\wbeps  \,dx - \int_{\Geps^0}(\wbeps \times {\bf n}_\vps)\cdot \varphi \,d\sigma\,,
\end{equation}
the previous equality can be rewritten as
\begin{align}
& \int_{\Omeps} D\ub_\eps  : D\varphi \,dx -\int_{\Omeps} p_\ep\, \dv \varphi  \,dx - 2N^2 \int_{\Omeps}\wb_\ep\cdot (\nabla\times \varphi) \,dx  \nonumber \\
& -2\left({1\over\alpha}-N^2\right) \int_{\Geps^0}(\wb_\ep\times {\bf n}_\ep)\cdot \varphi  \,d\sigma = 0\, .\label{WeakForm-u-Step1}
\end{align}
Multiplying Equation \eqref{31,1} by another test function $\psi\in \Veps$, integrating by parts and using boundary condition \eqref{bc-5}, we obtain
\begin{align}
& R_M\int_{\Omeps} D\wb_\ep : D\psi  \,dx - 2N^2\beta \int_{\Geps^0}(\ub_\ep\times {\bf n}_\ep)\cdot \psi  \,d\sigma + 4N^2\int_{\Omeps} \wb_\eps\cdot \psi  \,dx \nonumber \\
& - 2N^2\int_{\Omeps}(\nabla \times \ub_\eps)\cdot \psi  \,dx = 0 \, .\label{WeakForm-u-Step2}
\end{align}
Summing relations \eqref{WeakForm-u-Step1} and \eqref{WeakForm-u-Step2} yields 
\begin{align}
& \int_{\Omeps} D\ub_\eps  : D\varphi  \,dx  -\int_{\Omeps} p_\ep\, \dv \varphi  \,dx + R_M\int_{\Omeps} D\wb_\ep : D\psi  \,dx  - 2N^2\int_{\Omeps}(\nabla \times \ub_\eps)\cdot \psi  \,dx  \nonumber \\ & 
- 2N^2 \int_{\Omeps}\wb_\ep\cdot (\nabla\times \varphi)  \,dx  + 4N^2\int_{\Omeps} \wb_\eps\cdot \psi \,dx  -2\left({1\over\alpha}-N^2\right) \int_{\Geps^0}(\wb_\ep\times {\bf n}_\ep)\cdot \varphi \,d\sigma \label{WeakForm-u} \\ & 
  - 2N^2\beta \int_{\Geps^0}(\ub_\ep\times {\bf n}_\ep)\cdot \psi  \,d\sigma 
=  0 . \nonumber
\end{align}
This leads to the following definition. 
\begin{df}\label{Def:WeakSolution-u-w-p}
We say that $(\ubeps,\wbeps,p_\ep)\in H^1(\Omeps)^3\times H^1(\Omeps)^3\times L^2_0(\Omeps)$ is a weak solution to system~\eqref{1,1}--\eqref{bc-5} if $(\ubeps,\wbeps)$ satisfy boundary conditions~\eqref{bc-1}--\eqref{bc-3}, ${\rm div}\,\ubeps=0$ in $\Omega$ and relation~\eqref{WeakForm-u} holds for any $(\varphi,\psi)\in \Veps\times\Veps$. 
\end{df}

\section{Notation}\label{notation}
 The unitary cube  of $\RR^2$ will be denoted by  $Z'=(-{1\over
2},{1\over 2})^2$, and  we set  $\widehat Q=Z'\times
(0,+\infty)$. For any $M>0$, we define $\widehat Q_M=Z'\times
(0,M)$.
We introduce the space $L^2_{\#}(Z')$, which is defined by the functions $u$ in $L^2_{loc}(\RR^2)$ and $Z'$-periodic. The space $L^2_{\#}(\widehat Q)$ is defined by the functions $\widehat u$ in $L^2_{loc}(\RR^2\times
(0,+\infty))$ and
$$\int_{\widehat Q}|\widehat u|^2dz<+\infty,\quad \widehat u(z'+k',z_3)=\widehat
u(z),\quad \forall k'\in \ZZ^2,\quad \hbox{a.e. } z\in \RR^2\times
(0,+\infty).$$\par 

 We define $L^2_0(\mathcal{O})$, with   $\mathcal{O}$ a bounded and measurable subset of $
\RR^N$, by the functions of
$L^2(\mathcal{O})$ with zero integral.\par

For every $\theta'=(\theta_1,\theta_2)$, we define
$$[ \theta']^\perp=(- \theta_2, \theta_1),\quad {\rm rot}_{x_3} \theta'=\partial_{x_3}[\theta]^\perp,\quad {\rm Rot}_{x'}\theta'=\partial_{x_1}\theta_2-\partial_{x_2}\theta_1\,.$$

We define  the sets
$$ \Om_\ep^-=(\om\times (-\infty,0))\cap \Om_\ep,
\quad \Om_\ep^+=(\om\times (0,+\infty))\cap \Om_\ep.
$$

\par
Given $k'\in \ZZ^2$ and $\tau>0$, we define
$$C_{\tau}^{k'}= \tau Z'+\tau k',\quad Q_\tau^{k'}=(C_{\tau}^{k'}\times \RR)\cap \Theta_\ep,$$
where $\Theta_\varepsilon=\{x\in \mathbb{R}^2\times \mathbb{R}\ : \ .\Psi_\varepsilon(x')<x_3<\varepsilon\}$. We consider the function $\kappa:\RR^2\mapsto \ZZ^2$ given by
$$\kappa(x')=k' \Leftrightarrow x'\in C_1^{k'}.$$
We observe that $\kappa$ is well defined, except for a set of zero measure in
$\RR^2$. In addition, for any $\tau>0$, it holds
$$\kappa\left({x'\over \tau}\right)=k'\Leftrightarrow x'\in
C_\tau^{k'}.$$

We denote $C_{\ep^\ell}(x')$, for a.e. $x'\in \RR^2$, by the square
$C_{\ep^\ell}^{k'}$ such that $x'\in C_{\ep^\ell}^{k'}$.\par

Given $\rho>0$, we take 
\begin{equation}\label{Def:omega_rho}
	\omega_\rho=\{x\in\omega\,:\, {\rm dist}(x,\partial\omega)>\rho\},
	\end{equation}
$$I_{\rho,\ep}=\{k'\in\mathbb{Z}^2\,:\, \omega_\rho\cap C_{\ep^\ell}^{k'} \neq \emptyset\}.$$

 By $\mathcal{V}$ we define the space of functions $\widehat \varphi:
\RR^2\times(0,+\infty)\mapsto \RR$ such that $\widehat \varphi\in H^1_{\#}(\widehat Q_M)$, for every $M>0$, and $\nabla
\widehat \varphi\in L^2_{\#}(\widehat Q)^3$.  We remark that $\mathcal{V}$  is a Hilbert space by considering $\|\cdot\|_{\mathcal{V}}$ given by $$\|\widehat \varphi\|_{\mathcal{V}}^2=\|\nabla \widehat
\varphi\|_{L^2(\widehat Q)^3}^2+\|\widehat \varphi\|_{L^2(Z'\times\{0\})}^2.$$

We observe that when we use  $O_\ep$, we refer to a generic real sequence which is devoted to tend to zero
when  $\ep\to 0$. Moreover, $O_\ep$ is allowed to change change from line to line. By $C$, we denote a generic positive constant, which does not depend on $\varepsilon$ and  it can also change from line to line.

\section{Main results} \label{Sec:MainResults}

As discussed before, different asymptotic behaviours of the flow may be deduced depending on the order of magnitude of the viscosity coefficients. Indeed, if we compare the characteristic number $R_M$ defined by (\ref{non-dim}) and appearing in the equation (\ref{31,1}) with small parameter $\ep$, three different asymptotic situations can be formally identified (see e.g.   \cite{BayadaChambatGamouana, Ja9, SuaBMMS}). The most interesting one is, of course, the one leading to a strong coupling at main order, namely the regime 
\begin{equation}\label{RM-R_c}
R_M=\ep^2 R_c,\qquad R_c=\mathcal{O}(1).
\end{equation}
Hence, we will perform our analysis assuming the above scalings of $R_M$ and $R_c$ with respect to $\ep$. {\color{black} Concerning the other parameters, we recall that $N^2$, $\alpha$ and $\beta$ are of order $\mathcal{O}(1)$.}

Besides, in the case of a squeeze film model, we also assume that the (vertical) velocity of the upper plate $V_\vps$ is of order $\vps$ as $\vps$ tends to zero. Hence, we consider the asymptotic regime
\begin{equation}\label{Regime:Veps}
V_\vps = \vps S,
\end{equation}
where $S$ is a positive constant.

In order to study the asymptotic behaviour of the solution to system \eqref{1,1}--\eqref{bc-5}, we also need to assume a certain regularity on the boundary data $\gb_\ep$, and uniform estimates of relevant norms. A very general way of stating those properties is the following: there exists a sequence of lift functions $\mathbf{J}_\ep\in H^1(\Omega_\varepsilon)^3$ satisfying $\mbox{div}\,\mathbf{J}_\ep = 0$ in $\Omega_\ep$, the boundary  conditions
\begin{equation}\label{BC-lift-Jeps}
\mathbf{J}_\ep = -V_\ep e_3\ \hbox{on}\ \G_\ep^1,\quad \mathbf{J}_\ep = \mathbf{g}_\ep \ \hbox{on}\ \G_\ep^\ell,\quad \mathbf{J}_\ep\cdot {\bf n}_\ep = 0 \  \hbox{on}\ \G_\ep^0,
\end{equation}
and the estimates
\begin{equation}\label{Estimate-lift-Jeps}
\forall \ep>0\qquad 
\|\mathbf{J}_\varepsilon\|_{L^2(\Omega_\varepsilon)^3}\leq C\varepsilon^{1\over 2},\quad \| D \mathbf{J}_\varepsilon\|_{L^2(\Omega_\varepsilon)^{3\times 3}}\leq C\varepsilon^{-{1\over 2}},
\quad \|\Jb_\ep\|_{L^2(\Gamma^0_\ep)^3 } \leq C,
\end{equation}
where $C>0$ is a universal constant.

\begin{rem}
One typical construction of a boundary data $\gb_\ep$ and the associate lift function $\Jb_\ep$ is the following, see \cite{Bayada4}. Consider a regular vector field ${\Jb}\in H^1(\Omega)^3$, satisfying
\[
\mathrm{div}\, \Jb = 0\textrm{ in }\Omega,\quad \Jb =  -Se_3\textrm{ on }\omega\times\{h\},\quad \Jb = 0\textrm{ on }\omega\times\{0\}.
\]
Extending $\Jb=(J',J_3)$ by zero on $\omega\times(-\infty,0)$, we can define $\Jb_\ep\in H^1(\Omeps)^3$ by
\[
\Jb_\ep(x',x_3) = (J'(x',\frac{x_3}{\ep}),\ep J_3(x',\frac{x_3}{\ep}))\quad \forall (x',x_3)\in \Omeps,
\]
and $\gb_\vps:={\Jb_\ep}_{|{\G^{\ell}_\ep}}$ in the sense of traces. By the change of variable $(x',x_3) = (y',\ep y_3)$, there holds
\begin{align*}
\int_{\Omega_\ep}|D \Jb_{\ep}|^2\, dx'dx_3 & = \ep\int_{\Omega}\left(|D_{y'} J'|^2 + \frac{1}{\ep^2}| \partial_{y_3} J'|^2 + \ep^2|\nabla_{y'}J_3|^2 + | \partial_{y_3}J_3|^2 \right)\, dy'dy_3,\\
\int_{\Omega_\ep}|  \Jb_{\ep}|^2\, dx'dx_3 & = \ep \int_{\Omega} \left( |J'|^2 + \ep^2 |J_3|^2 \right)\, dy'dy_3,
\end{align*}
so that $\Jb_\ep$ satisfies all the required properties \eqref{BC-lift-Jeps}-\eqref{Estimate-lift-Jeps}.

{\color{black} Since such vector field $\Jb$ is not unique, the lift function $\Jb_{\ep}$ and the boundary data $\gb_\ep$ are quite arbitrary. In fact, they do not play a significant role in the asymptotic analysis of the problem, provided that conditions \eqref{BC-lift-Jeps}-\eqref{Estimate-lift-Jeps} are satisfied.}
\end{rem}

Let us start with an existence and uniqueness result for the solution of problems \eqref{1,1}--\eqref{bc-5}, whose  proof is given in the Section \ref{sec:proofs}.
\begin{tm}\label{thm_existence}
Assume that the coupling parameter $N^2$ satisfies the condition
\begin{equation}\label{ExtraConditionOnN}
N^2\leq \frac{1}{2}\, ,
\end{equation}
and define the nonnegative parameter $\gamma$ by
\begin{equation}\label{Def:gamma}
\gamma = \frac{1}{\alpha}-N^2-N^2\beta\, .
\end{equation}
Assume that the asymptotic regimes \eqref{RM-R_c} and \eqref{Regime:Veps} hold. Then, for any $\beta$ such that
\begin{equation}\label{Condition-gamma}
\gamma^2 < \frac{R_c(1-2N^2)}{h^2}\, ,
\end{equation}
there exists $\ep_0>0$ such that for any $0<\ep < \ep_0$,  there exists a unique weak solution $(\ub_\ep,\wb_\ep,p_\ep)$ in $H^1(\Omeps)^3\times H^1(\Omeps)^3\times L^2_0(\Omeps)$ to system \eqref{1,1}--\eqref{bc-5} (in the sense of Definition~\ref{Def:WeakSolution-u-w-p}). 
\end{tm}
\begin{rem}
{\color{black}In the case of a flat boundary, Bayada \emph{et al.}\! obtained in~\cite{Bayada4} less restrictive conditions, namely $N^2<1$ and $\gamma^2<\frac{R_c(1-N^2)}{h^2}$. However, we stress that the restriction of parameter $N$ (\ref{ExtraConditionOnN}) that it is necessary to guarantee existence and uniqueness of the weak solution, is in fact in agreement with tribology models, where different considerations lead to the same assumption $N^2\leq 1/2$ (see \cite{Singh1, Singh2}). }
\end{rem}

\subsection{Rescaling}

We also want to describe the asymptotic behaviour of the sequence $({\bf u}_\ep, {\bf w}_\ep, p_\ep)$  of solution of the micropolar system (\ref{1,1})-(\ref{31,1}) supplemented with boundary conditions (\ref{bc-4})-(\ref{bc-5}), as $\ep$ tends to $0$. We start by introducing a change of variables classically used in asymptotic analysis of flows in thin domains: the dilatation
\begin{equation}\label{changevar1}
y'=x',\quad y_3={x_3\over \ep},
\end{equation}
which changes $\Omega_\ep$ to the set $\widetilde \Omega_\ep$ of height of order $h$, defined as follows:
\begin{equation}\label{tilOmep}
\widetilde \Omega_\ep=\left\{(y',y_3)\in\RR^2\times \RR\ :\  y'\in \omega,\quad -\widetilde \Psi_\ep(y')<y_3<h\right\},
\end{equation}
where
$$\widetilde \Psi_\ep(y')={1\over \ep} \Psi_\ep(y')=\ep^{\delta-1}\Psi\left({1\over \ep^\ell}y'\cdot e'_1\right).$$
The lower, upper and lateral boundaries of the rescaled domain $\widetilde \Omega_\ep$  are now defined by
\begin{align*}
\tilde{\Gamma}^0_\ep &=\left\{( y', y_3)\in\RR^2\times \RR\ :\   y'\in   \omega,\quad y_3=- \tilde \Psi_\ep(y') \right\},\\
\tilde\Gamma^{1}_\ep &=\left\{( y', y_3)\in\RR^2\times \RR\ :\   y'\in \omega,\quad y_3= h\right\},\\
  \tilde\Gamma^{\ell}_\ep & =\partial \tilde\Omega_\ep-(\tilde\Gamma^0_\ep \cup \widetilde\Gamma^{1}_\ep).
\end{align*}
Accordingly, we define the functions $\widetilde {\bf u}_\ep$, $\widetilde {\bf w}_\ep\in H^1(\widetilde\Omega_\ep)^3$  and $\widetilde p_\ep\in L^2_0(\widetilde\Omega_\ep)$ by 
\begin{equation}\label{changevar1_fun}
\begin{array}{c}\dis \widetilde {\bf u}_\ep(y)={\bf u}_\ep(y',\ep y_3),\quad  \widetilde {\bf w}_\ep(y)={\bf w}_\ep(y',\ep y_3),\quad \widetilde p_\ep(y)=p_\ep(y',\ep y_3),
\end{array}\quad\hbox{a.e. } y\in \widetilde\Omega_\ep\, .
\end{equation}

Since $\delta>1$, it is clear that the sequence of domains $\widetilde \Omega_\ep$ converges (for instance, in the sense of Hausdorff complementary topology) to the limit domain $\Omega$ defined by
\begin{equation*}
\Omega= \omega\times (0,h)\subset\RR^2\times \RR\, .
\end{equation*}
We denote by $\Gamma:=\omega\times\{0\}$ the lower boundary of $\Omega$.

The next step of the analysis is to identify the effective system satisfied by these rescaled functions.

\subsection{Effective system}
In this subsection, we give  the result concerning the asymptotic behaviour of the rescaled functions $\widetilde {\bf u}_\ep$, $\widetilde {\bf w}_\ep$, $\widetilde p_\ep$. Depending on the relation between the amplitude parameter $\delta$ and the period parameter $\ell$, we obtain three different regimes, that we call critical, sub-critical and super-critical. Here, we state the result in the critical case; the other cases will be discussed in Remark \ref{rem_sub_sup}. The proof of the corresponding results is given in Section \ref{sec:proofs}.

\begin{tm}\label{thm_effective}
Assume that the asymptotic regimes \eqref{RM-R_c} and \eqref{Regime:Veps} and conditions (\ref{ExtraConditionOnN}) and (\ref{Condition-gamma}) hold.   Assume that $\delta, \ell$ satisfy the relation $\delta= {3\over 2}\ell-{1\over 2}$ (critical case). Let $({\bf u}_\ep, {\bf w}_\ep, p_\ep)$ be a sequence of weak solutions of (\ref{1,1})-(\ref{bc-5}). Then, there exist $\widetilde {\bf u}', \widetilde {\bf w}'\in  H^1(0,h;L^2(\omega))^2$ and $ p\in H^1(\omega)\cap L^2_0(\omega)$ such that the rescaled  functions $\widetilde {\bf u}_\ep, \widetilde {\bf w}_\ep, \widetilde p_\ep$ satisfy
\begin{equation}\label{convergencias_solutions}
\begin{array}{c}
\widetilde {\bf u}_\ep\rightharpoonup (\widetilde {\bf u}',0)\hbox{ in } H^1(0,h;L^2(\omega))^3,\quad \ep\widetilde {\bf w}_\ep\rightharpoonup (\widetilde {\bf w}',0)\hbox{ in }  H^1(0,h;L^2(\omega))^3,\\
\noame
\ep^2\widetilde p_\ep\to p \hbox{ in }L^2(\Omega)\,.
\end{array}
\end{equation}
The triplet $(\widetilde \ub',\widetilde \wb',  p)$ is the unique solution of the  following problem
\begin{equation}\label{limit_system_1}\left\{\begin{array}{rl}
\dis
-\partial_{y_3}^2 \widetilde \ub' + \nabla_{y'}   p -2N^2 {\rm rot}_{y_3} \widetilde \wb'=0& \hbox{ in } \Omega\, ,\\
\noame
-R_c\partial_{y_3}^2 \widetilde \wb'+4N^2 \widetilde \wb'-2N^2{\rm rot}_{y_3}\widetilde \ub'=0& \hbox{ in } \Omega\, ,\\
\noame
\dis \mathrm{div}_{y'}\int_0^{h} \widetilde {\bf u}'(y',y_3)\,dy_3=S& \hbox{ in }  \omega\, ,
\end{array}\right.\end{equation}
with  the boundary conditions 
\begin{equation}\label{bc_top_system_1}
\hspace{-0.8cm}\widetilde \ub'=0,\ \widetilde \wb'=0 \hbox{ on }\omega\times\{h\},
\end{equation}
\begin{equation}\label{bc_botom_case2_system_1}
\partial_{y_3} \widetilde \ub'=-{2\over \alpha} [ \widetilde \wb']^\perp + E_\lambda \widetilde (\ub'\cdot e'_1)\, e_1'\ \hbox{ on }\Gamma,\quad 
 R_c\partial_{y_3} \widetilde \wb'=-2N^2\beta[\widetilde \ub']^\perp+R_c F_\lambda\widetilde (\wb'\cdot e_1')\, e_1'\ \hbox{ on }\Gamma\,.
\end{equation}
Coefficients $E_\lambda, F_\lambda \in \mathbb{R}$ appearing in boundary conditions \eqref{bc_botom_case2_system_1} are defined by
\begin{equation}\label{matrixM}
E_\lambda=\int_{\widehat Q}|D_z \widehat \phi^{1,\lambda}|^2\,dz,\quad F_\lambda=\int_{\widehat Q}|D_z \widehat \phi^{2,\lambda}|^2\,dz,
\end{equation}
where
$(\widehat {\phi}^{i,\lambda},\widehat q^{i,\lambda})\in\mathcal{V}^3\times L^2_\sharp(\widehat Q)$, $i=1,2$, are respectively the solutions of 
\begin{equation}\label{system_phi_1}
\left\{\begin{array}{rll}
-\Delta_z\widehat {\phi}^{1,\lambda} + \nabla_z\widehat q^{1,\lambda}=0&\hbox{ in }&\mathbb{R}^2\times \mathbb{R}^+,\\
\noame\dis
{\rm div}_z\widehat{\phi}^{1,\lambda}=0&\hbox{ in }&\mathbb{R}^2\times \mathbb{R}^+,\\
\noame\dis
\widehat{\phi}_3^{1,\lambda}(z',0)=\lambda\partial_{z_1}\Psi(  z'\cdot e'_1)&\hbox{ on }&\mathbb{R}^2\times \{0\},\\
\noame\dis
-\partial_{z_3}\widehat{\phi}^{1,\lambda}_1=0,\ -\partial_{z_3}\widehat{\bf \phi}^{1,\lambda}_2=0&\hbox{ on }&\mathbb{R}^2\times \{0\}\,,
\end{array}\right.
\end{equation}
and 
\begin{equation}\label{system_phi_2}
\left\{\begin{array}{rll}
-\Delta_z\widehat \phi^{2,\lambda} =0&\hbox{ in }&\mathbb{R}^2\times \mathbb{R}^+,\\
\noame\dis
\widehat\phi_3^{2,\lambda}(z',0)=\lambda\partial_{z_1}\Psi( z'\cdot e'_1)&\hbox{ on }&\mathbb{R}^2\times \{0\},\\
\noame\dis
-\partial_{z_3}\widehat\phi^{2,\lambda}_1=0,\ -\partial_{z_3}\widehat\phi^{2,\lambda}_2=0&\hbox{ on }&\mathbb{R}^2\times \{0\}\,.
\end{array}\right.
\end{equation}
\end{tm}

\begin{rem}\label{rem_sub_sup} 
Theorem \ref{thm_effective} can be adapted easily to describe the two other asymptotic regimes:\\

\begin{itemize}
\item[$\bullet$] In the sub-critical case $\delta>{3\over 2}\ell-{1\over 2}$, the riblets are so small that there is no effect of roughness, so we obtain the nonzero boundary conditions on $\Gamma$,
$$\partial_{y_3} \widetilde \ub'=-{2\over \alpha} [ \widetilde \wb']^\perp \ \hbox{ on }\Gamma,\quad 
 R_c\partial_{y_3} \widetilde \wb'=-2N^2\beta[\widetilde \ub']^\perp\ \hbox{ on }\Gamma\,.
 $$
Thus, we deduce that the model obtained in \cite{Bayada4, Bayada5} even holds for a very slightly rough boundary.
\item[$\bullet$] In the super-critical case $1<\delta<{3\over 2}\ell-{1\over 2}$, the effect of the riblets is maximal. Thus, boundary conditions given in (\ref{bc_botom_case2_system_1}) are replaced by
\begin{equation}\label{super_critical_bc}\widetilde \ub'\cdot e'_1=\widetilde \wb'\cdot e_1'=0\quad\hbox{on }\Gamma,\quad \partial_{y_3}\widetilde \ub' \cdot e'_2=\partial_{y_3}\widetilde \wb\cdot e'_2=0\quad\hbox{on }\Gamma\,.
\end{equation}
Thus, we deduce that the roughness is so strong that the fluid adheres to the boundary and fluid elements cannot rotate on the fluid-solid interface in the $x_1$-direction.
\end{itemize}
\end{rem}

\subsection{Generalized micropolar Reynolds equations}\label{Sec:generalmicropolar}
In this subsection, we obtain a generalized Reynolds equation associated to the homogenized micropolar system given in Theorem \ref{thm_effective} (critical case). For the sake of simplicity, we will consider a 2D domain characteristic of the lubrication assumption. Thus, we consider in Theorem \ref{thm_effective} that 
the flow does not depend on the $y_2$-coordinate, and that velocity component $\widetilde u_2$ and micropolar component $\widetilde w_1$ are zero. Hence, we address the following limit problem posed in $\Omega=(0,1)\times (0,h)$:
\begin{equation}\label{limit_system_1_reynolds}\left\{\begin{array}{r}
\dis
-\partial_{y_3}^2 \widetilde u_1 + \partial_{y_1}   p +2N^2 \partial_{y_3} \widetilde w_2=0\ \hbox{ in } \Omega,
\\
\noame\dis
-R_c\partial_{y_3}^2 \widetilde w_2+4N^2 \widetilde w_2-2N^2\partial_{y_3}\widetilde u_1=0\ \hbox{ in } \Omega,\end{array}\right.
\end{equation}
completed with the boundary conditions 
\begin{equation}\label{bc_top_system_1_reynolds}
\hspace{-0.8cm}\widetilde u_1=0,\ \widetilde w_2=0 \hbox{ on }\Gamma^1=(0,1)\times\{h\},
\end{equation}
\begin{equation}\label{bc_botom_case2_system_1_reynolds}
\partial_{y_3} \widetilde u_1={2\over \alpha}  \widetilde w_2 + E_\lambda \widetilde u_1\ \hbox{ on }\Gamma=(0,1)\times\{0\},\quad 
 R_c\partial_{y_3} \widetilde w_2=-2N^2\beta \widetilde u_1\ \hbox{ on }\Gamma\,,
\end{equation}
and the incompressibility condition
\begin{equation}\label{limit_system_1_reynolds2}
\partial_{y_1}\int_0^{h} \widetilde u_1(y_1,y_3)\,dy_3=S\quad \hbox{ in } (0,1)\,.
\end{equation} 
{\color{black} We give in the Appendix the expression of $(\tilde u_1,\tilde w_2)$, solution of system (\ref{limit_system_1_reynolds})--(\ref{bc_botom_case2_system_1_reynolds}), in terms of $p$ (see Lemmas~\ref{lemma_alpha_neq_1} and~\ref{lemma_alpha_igual_1}). }
Putting these expressions in (\ref{limit_system_1_reynolds2}) will lead to the corresponding Reynolds equations that take into account the roughness-induced effects.

\begin{tm}\label{thm_reynolds}
In the critical case $\delta = \frac{3}{2}\ell-\frac{1}{2}$, the pressure $p$ satisfies the following Reynolds equation
\begin{equation}\label{reynolds}
\int_0^1\Theta_{\lambda}\partial_{y_1}  p(y_1)\,\partial_{y_1}\theta(y_1)\,dy_1=\int_0^1 S\theta(y_1)\,dy_1,\quad\forall \theta\in H^1(0,1),
\end{equation}
with  $\Theta_{\lambda}$ defined in the case $\alpha\neq 1$ by
\begin{equation}\label{theta1_alpha_neq_1}
\begin{array}{rl}
\Theta_{\lambda} =&  {h^3\over 3(1-N^2)}-(1-\eta_\lambda){3h^3\over 4(1-N^2)}\\
\noame
& -\left({2N^2\over k}\Big[{ch(kh)-1\over k}-\eta_\lambda h sh(kh)\Big]+{\gamma_\alpha\over 2} h^2(1-2\eta_\lambda)
-(1-\eta_\lambda)\left[\gamma_\alpha h+{2N^2\over k}sh(kh)\right]\right)A \\
\noame
& -\left({2N^2\over k}\left[{sh(kh)\over k}-\eta_\lambda h ch(kh)\right]-(1-\eta_\lambda)(1+ch(kh))h{N^2\over k}\right)B\,,
\end{array}
\end{equation}
and in the case $\alpha= 1$ by
\begin{equation}\label{theta1_alpha_igual_1}
\begin{array}{rl}
\Theta_{\lambda} =& -{1\over 2(1-N^2)}\Big({h^3\over 3}-\mu_\lambda h^3\Big)-(1-\mu_\lambda){h^2\over k(1-N^2)}{1-ch(kh)\over sh(kh)}
\\
\noame & 
-\Big[{1\over 1-N^2}\Big({h^2\over 2}-\mu_\lambda h^2\Big)+(1-\mu_\lambda){h\over k(1-N^2)}{1-ch(kh)\over sh(kh)}\Big]A'
\\
\noame &
 -\Big[{2N^2\over k}\Big({sh(kh)\over k}- \mu_\lambda h ch(kh)\Big)-(1-\mu_\lambda){2N^2\over k}\Big(h+{(1-ch(kh))^2\over k\, sh(kh)}\Big)\Big]B' \,,
\end{array}
\end{equation}
where $A, A', B$ and $B'$  are defined in Lemmas \ref{lemma_alpha_neq_1} and \ref{lemma_alpha_igual_1} in the Appendix.
\end{tm}

\begin{rem} It is worth mentioning that the effective expressions given in Lemmas \ref{lemma_alpha_neq_1} and \ref{lemma_alpha_igual_1}, and Theorem \ref{thm_reynolds} are explicitly corrected by the roughness-induced
 coefficient  $E_\lambda$. Indeed, by putting $E_\lambda = 0$, which implies $\eta_\lambda=\mu_\lambda=1$  (i.e., no roughness introduced), we obtain the same expressions as derived in \cite{Bayada4}, which also corresponds to the sub-critical case $\delta> {3\over 2}\ell-{1\over 2}$.
\end{rem}

{\color{black} Using the explicit expressions from Lemma~\ref{lemma_alpha_neq_1}, Lemma~\ref{lemma_alpha_igual_1} and Theorem~\ref{thm_reynolds}, it is possible to develop $\tilde u_1, \tilde w_2$ and $p$ in powers of $\lambda^2$. This will be useful in the numerical computations from Section~\ref{Sec:num}. However, since the corresponding formulas are rather long, we have gathered them in the Appendix for the sake of clarity (see Corollary~\ref{cor_alpha_neq_1}).}


\medskip

{\color{black} Finally, we give the micropolar Reynolds equation corresponding  to the super-critical case $1<\delta<{3\over 2}\ell-{1\over 2}$. As in the critical case, its derivation is based on explicit expressions of the velocity and microrotation (see Lemma~\ref{lm_super_critical} in the Appendix).}

\begin{tm}\label{thm_reynolds_super}
In the super-critical case $1<\delta<{3\over 2}\ell-{1\over 2}$, the pressure $p$ satisfies the following Reynolds equation 
\begin{equation}\label{reynolds}
\int_0^1\Theta\partial_{y_1}  p(y_1)\,\partial_{y_1}\theta(y_1)\,dy_1=\int_0^1 S\theta(y_1)\,dy_1,\quad\forall \theta\in H^1((0,1)),
\end{equation}
with  $\Theta$ defined by
$$
\Theta={h^3\over 12(1-N^2)}-{2N^2\over k}\Big[{ch(kh)-1\over k}-{h\over 2}sh(kh)\Big]A''-{2N^2\over k}\Big[{sh(kh)\over k}-{h\over 2}(ch(kh)+1)\Big]B'',\\
$$
where $A''$ and $B''$  are defined in Lemma \ref{lm_super_critical}.
\end{tm}

\section{Proofs of the main results}\label{sec:proofs}
We start by proving the existence and uniqueness of solution of problem~\eqref{1,1}--\eqref{bc-5}.\\

\noindent {\bf Proof of Theorem~\ref{thm_existence}. } Let $\Jb_\ep\in H^1(\Omeps)^3$ be a sequence of free divergence lift functions satisfying \eqref{BC-lift-Jeps}-\eqref{Estimate-lift-Jeps}. Replacing $\ubeps$ by $\vbeps+\Jb_\ep$ in the weak formulation~\eqref{WeakForm-u}, we see that $(\ubeps,\wbeps,p_\ep)$ is a weak solution to system~\eqref{1,1}--\eqref{bc-5} if and only if $(\vbeps,\wb_\eps,p_\ep)\in \Veps^0\times\Veps\times L^2_0(\Omeps)$ and satisfies for any $(\varphi,\psi)\in \Veps\times \Veps$
\begin{align}
& \int_{\Omeps} D\vb_\eps  : D\varphi  \,dx  - \int_{\Omeps} p_\ep\, \dv \varphi \,dx + R_M\int_{\Omeps} D\wb_\ep : D\psi \,dx - 2N^2\int_{\Omeps}(\nabla \times \vb_\eps)\cdot \psi \,dx  \nonumber \\
& - 2N^2 \int_{\Omeps}\wb_\ep\cdot (\nabla\times \varphi)  \,dx + 4N^2\int_{\Omeps} \wb_\eps\cdot \psi\,dx  -2\left({1\over\alpha}-N^2\right) \int_{\Geps^0}(\wb_\ep\times {\bf n}_\ep)\cdot \varphi\,d\sigma \nonumber
   \\
& - 2N^2\beta \int_{\Geps^0}(\vb_\ep\times {\bf n}_\ep)\cdot \psi \,d\sigma \nonumber\\
& =  - \int_{\Omeps}D\Jb_\ep : D\varphi \,dx + 2N^2\int_{\Omeps}(\nabla\times \Jb_\eps)\cdot \psi \,dx   +  2N^2\beta \int_{\Geps^0}(\Jb_\ep\times {\bf n}_\ep)\cdot \psi\,d\sigma \, .  \label{WeakForm} 
\end{align}
Equation \eqref{WeakForm} justifies the introduction of the bilinear forms $\mathcal{A}_\vps:(\Veps\times \Veps)^2 \rightarrow \RR$ and  $\mathcal{B}_\vps:(\Veps\times \Veps)\times L^2_0(\Omeps) \rightarrow \RR$  respectively defined by
\begin{equation}\label{Def:BilinearFormA} 
\begin{array}{rl}
\dis \mathcal{A}_\vps((\vb,\wb),(\varphi,\psi))  & \dis =
 \int_{\Omeps} D\vb   : D\varphi \,dx - \int_{\Omeps} p_\ep\, \dv \varphi \,dx + R_M\int_{\Omeps} D\wb  : D\psi \,dx \\
\noame
 &\dis - 2N^2\int_{\Omeps}(\nabla \times \vb )\cdot \psi  \,dx - 2N^2 \int_{\Omeps}\wb \cdot (\nabla\times \varphi) \,dx + 4N^2\int_{\Omeps} \wb \cdot \psi \,dx \\
\noame &\dis -2\left({1\over\alpha}-N^2\right) \int_{\Geps^0}(\wb\times {\bf n}_\ep)\cdot \varphi \,d\sigma 
  - 2N^2\beta \int_{\Geps^0}(\vb\times {\bf n}_\ep)\cdot \psi  \,d\sigma \, ,
\end{array}\end{equation}
and
\begin{equation}\label{Def:BilinearFormB}
\dis \mathcal{B}_\vps((\vb,\wb),q) =  -\int_{\Omeps} q\, \dv \vb \,dx \, ,\hfill
\end{equation}
and of the linear form $\mathcal{L}_\vps : \Veps\times \Veps \rightarrow \RR $ defined by
\begin{equation}\label{Def:LinearFormL}
\begin{array}{rl}
\dis \mathcal{L}_\vps(\varphi,\psi) &\dis =   - \int_{\Omeps}D\Jb_\ep : D\varphi  \,dx + 2N^2\int_{\Omeps}(\nabla\times \Jb_\eps)\cdot \psi \,dx +  2N^2\beta \int_{\Geps^0}(\Jb_\ep\times {\bf n}_\ep)\cdot \psi \,d\sigma\, . 
\end{array}\end{equation}
Hence, $(\ubeps,\wbeps,p_\eps)$ is a weak solution to system~\eqref{1,1}--\eqref{bc-5} if and only if $(\vbeps,\wbeps,p_\eps)\in \Veps\times\Veps\times L^2_0(\Omeps)$ and satisfies the following mixed formulation
\begin{align}
\mathcal{A}_\vps((\vb,\wb),(\varphi,\psi)) +\mathcal{B}_\vps((\varphi,\psi),p_\ep)  &= \mathcal{L}_\vps(\varphi,\psi) &&\forall  (\varphi,\psi)\in \Veps\times\Veps\, , \label{MixedForm1}\\
\mathcal{B}_\vps((\vbeps,\wbeps),q) &= 0  &&\forall q\in L^2(\Omeps) \, . \label{MixedForm2}
\end{align}

The existence and uniqueness of the solution $(\vb_\ep,\wb_\ep,p_\ep)$ to the mixed formulation~\eqref{MixedForm1}-\eqref{MixedForm2} is established in \cite{Bayada4}, in the case where the oscillating boundary $\G^0_\ep$ is replaced by a flat boundary $\omega\times \{0\}$. For the sake of completeness, we recall the main steps of the proof, highlighting the differences that are implied by the oscillations of the lower boundary $\G^0_\vps$.

First, let us state some useful quantitative inequalities.

\paragraph{Trace inequality on $\Geps^0$}

Since the lower boundary $\Geps^0$ is not flat, one needs to take into account the variations of the normal direction ${\bf n}_\ep$ in order to estimate the $L^2$-norm of the trace of a function $\psi\in \Veps$. To this end, we introduce the quantity $\tau_\ep$ defined by
\[
\tau_\ep:= \sup_{x'\in \omega}\sqrt{1+|\nabla_{x'}\Psi_\ep(x')|^2}.
\]
We also denote by $h_\vps$ the height of the domain $\Omeps$, defined by
\[
h_\vps:=\sup_{x'\in \omega} \big(\vps h + \Psi_\ep(x')\big) = \ep \sup_{x'\in \omega} \big(h +\lambda \ep^{\delta-1}\Psi({1\over \ep^\ell}x'\cdot e'_1)\big)\, .
\]
In particular, $\lim_{\ep\rightarrow 0} h_\ep/\ep = h$.
With this notation, there holds the trace inequality
\begin{equation}\label{Ineq-TraceOmeps}
\forall \psi \in \Veps\quad \|\psi\|_{L^2(\G^0_\ep) ^3}\leq \sqrt{\tau_{\ep}h_\ep} \|D \psi\|_{L^2(\Omeps)^{3\times 3}}.
\end{equation}
In fact, this inequality holds true for any vector field $\psi\in H^1(\Omeps)^3$ vanishing on $\Geps^1$. By density, it is enough to prove it for any $\psi\in H^1(\Omeps)^3\cap C^1(\overline\Omeps)^3$ such that $\psi_{|\Geps^1}=0$. Integrating on vertical lines, we obtain
\begin{align*}
\int_{\G^0_\ep} |\psi|^2 & = \int_{\omega} |\psi(x',-\Psi_\ep(x'))|^2\, \sqrt{1+|\nabla_{x'}\Psi_\ep(x')|^2}\ dx' \\
& \leq \tau_{\ep} \int_{\omega} |\psi(x',-\Psi_\ep(x'))|^2 \ dx' \\
& \leq \tau_{\ep} \int_{\omega} \Big|\int_{-\Psi_\ep(x')}^{\ep h} \partial_{x_3}\psi(x',x_3)\ dx_3\Big|^2 \ dx' \\
&\leq  \tau_{\ep} \int_{\omega} (\ep h +  \sup_{\omega}\Psi_\ep ) \Big(\int_{-\Psi_\ep(x')}^{\ep h} |\partial_{x_3}\psi(x',x_3)|^2\ dx_3\Big) dx' \\
&\leq \tau_{\ep} h_\ep   \int_{\Omeps} |D \psi|^2 \,dx\, .
\end{align*}
This proves \eqref{Ineq-TraceOmeps}.

\paragraph{Poincar\'e inequality}

In the same fashion, Poincar\'e inequality in $\Veps$ reads
\begin{equation}\label{Ineq-PoincareOmeps}
\forall \varphi\in \Veps\quad \|\varphi\|_{L^2(\Omeps)^3} \leq {h_\ep} \|D \varphi\|_{L^2(\Omeps)^{3\times 3}}.
\end{equation}

\paragraph{Relation between $\|\nabla  \times \varphi\|_{L^2(\Omeps)^3}$ and $\|D\varphi\|_{L^2(\Omeps)^{3\times 3}}$}
Let us recall that for any vector field $\varphi\in \Veps$,
\begin{equation}\label{Identity:L2norms-gradient-curl}
\int_{\Omeps}\big(|\dv \varphi|^2+|\nabla\times \varphi|^2\big) \,dx = \int_{\Omeps}|D\varphi|^2\,dx + \int_{\Geps^0}((\varphi\cdot\nabla){\bf n}_\ep)\cdot\varphi\,d\sigma \, ,
\end{equation}
(see, for instance, \cite{BoyerFabrie}, formula~(IV.23)). In particular, if the lower boundary $\Geps^0$ was flat, the identity $\|\nabla\times \varphi\|_{L^2(\Omeps)^3}^2 =  \|D\varphi\|_{L^2(\Omeps)^{3\times 3}}^2$ would hold for any $\varphi\in \Veps^0$, since the remaining term $ \int_{\Geps^0}((\varphi\cdot\nabla){\bf n}_\ep)\cdot\varphi\,d\sigma$ would vanish. However, in the present geometric configuration, one cannot expect this term to be zero in general. In fact, a classical estimate reads
\[
\Big| \int_{\Geps^0}((\varphi\cdot\nabla){\bf n}_\ep)\cdot\varphi  \,d\sigma \Big| \leq  \mathrm{Lip}({\bf n}_\ep)\, \|\varphi\|_{L^2(\Geps^0)^3}^2 \, ,
\]
where $\mathrm{Lip}({\bf n}_\ep)$ is the Lipschitz constant of the normal vector field ${\bf n}_\ep$, locally extended in a neighborhood of the surface $\{x_3=\Psi_\ep(x')\}$ (in the sense of~\cite[Section~3.4]{BoyerFabrie}). However, using definition of $\Psi_\ep$ given in \eqref{rugprofile-nondim} and condition~\eqref{relation_parameters} on parameters $\delta,\ell$, it turns out that in the general case where $\|\partial^2_{11}\Psi\|_{\infty}>0$, $\mathrm{Lip}({\bf n}_\ep)$ is of order $\ep^{\delta-2\ell}$, hence diverging since $\ep^{\delta-\ell}$ goes to zero. As a result, we cannot use identity~\eqref{Identity:L2norms-gradient-curl} to estimate the $L^2$ norms of  $\dv \varphi$ and $\nabla\times \varphi$ over $\Omeps$ by the $L^2$ norm of $D\varphi$, as is done in~\cite{Bayada4} in the case of a flat boundary.

Instead, we rely on the following elementary estimate:
\begin{equation}
	\|\nabla\times \varphi\|_{L^2(\Omeps)^3}\leq \sqrt{2}\|D\varphi\|_{L^2(\Omeps)^{3\times 3}} \quad \forall\varphi\in H^1(\Omeps)^3\, .\label{Estimate:Sqrt2}
\end{equation}
The presence of the constant $\sqrt{2}$ in the previous estimate is at the origin of the term $2N^2$ in condition~\eqref{Condition-gamma}, which was simply $N^2$ in the case of a flat surface, as established in~\cite{Bayada4}.

\paragraph{Existence and uniqueness of the solution of the mixed formulation~\eqref{MixedForm1}-\eqref{MixedForm2}}

Using Cauchy-Schwarz inequality and inequalities~\eqref{Ineq-TraceOmeps}, \eqref{Ineq-PoincareOmeps} and \eqref{Estimate:Sqrt2}, it is easy to see that $\mathcal{A}_\ep,\mathcal{B}_\ep$ and $\mathcal{L}_\ep$ are continuous on their respective domains of definition, for any fixed value of parameter $\ep$. Hence, noticing that by definition of $\Veps^0$,
\[
\Veps^0\times \Veps = \left\lbrace (\varphi,\psi)\in \Veps\times \Veps,\ \mathcal{B}_{\ep}((\varphi,\psi),q)=0 \ \text{for any }q\in L^2_0(\Omeps) \right\rbrace\, ,
\]
and denoting by $\|(\cdot,\cdot)\|_{\Veps\times\Veps}$ the norm defined by
\[
\|(\varphi,\psi)\|_{\Veps\times\Veps} = \Big( \|D \varphi\|_{L^2(\Omeps)^{3\times 3}}^2 + \|D \psi\|_{L^2(\Omeps)^{3\times 3}}^2  \Big)^{1/2}
\]
the existence and uniqueness of the solution $(\vbeps,\wbeps,p_\ep)$ to the mixed formulation~\eqref{MixedForm1}-\eqref{MixedForm2} result from the following properties (see~\cite[paragraph 4.1 p.\! 57]{GiraultRaviart}) :
\begin{itemize}
\item[(i) coerciveness of $\mathcal{A}_{\ep}$:] there exists $\eta=\eta(\ep)>0$ such that
\[
\forall (\varphi,\psi)\in \Veps^0\times\Veps\quad
\mathcal{A}_{\ep}((\varphi,\psi),(\varphi,\psi))\geq \eta \big( \|D\varphi\|_{L^2(\Omeps)^{3\times 3}}^2 +  \|D\psi\|_{L^2(\Omeps)^{3\times 3}}^2 \big)\, ,
\]
\item[(ii) inf-sup condition:] there exists $c=c(\ep)>0$ such that
\[
\inf_{q\in L^2_0(\Omeps)} \sup_{(\varphi,\psi)\in \Veps\times \Veps} \frac{\mathcal{B}_\ep((\varphi,\psi),q)}{\|(\varphi,\psi)\|_{\Veps\times \Veps}\|q\|_{L^2_0(\Omeps)}}  \geq c \, .
\]
\end{itemize}

The inf-sup condition (ii) can be proved using the exact same arguments as in the proof of Theorem~2.2 in~\cite{Bayada4}, that relies on the solvability in $H^1_0(\Omeps)^3$ of equation $\dv \varphi = q$, for an arbitrary $q\in L^2_0(\Omeps)$, with natural estimates.

To establish the coerciveness condition (i), we use H\"older inequality, Poincar\'e inequality~\eqref{Ineq-PoincareOmeps}, the trace inequality~\eqref{Ineq-TraceOmeps} and estimate~\eqref{Estimate:Sqrt2} to obtain the lower estimate
\begin{align}
\Ac_\ep((\varphi,\psi),(\varphi,\psi))  =& \int_{\Omeps} |D\varphi|^2 \,dx  + R_M\int_{\Omeps} |D\psi|^2 \,dx 
- 4N^2\int_{\Omeps}(\nabla \times \varphi )\cdot \psi  \,dx  + 4N^2\int_{\Omeps} |\psi|^2  \,dx \nonumber\\
&  -2\left({1\over\alpha}-N^2\right) \int_{\Geps^0}(\psi\times {\bf n}_\ep)\cdot \varphi \,d\sigma 
  - 2N^2\beta \int_{\Geps^0}(\varphi\times {\bf n}_\ep)\cdot \psi  \,d\sigma \nonumber\\
\geq & \|D\varphi\|_{L^2(\Omeps)^{3\times 3}}^2 + R_M \|D\psi\|_{L^2(\Omeps)^{3\times 3}}^2 - 4N^2\sqrt{2} \|D\varphi\|_{L^2(\Omeps)^{3\times 3}} \| \psi\|_{L^2(\Omeps)^3}   \label{CoercivityA-1} \\
& +4N^2\|\psi\|_{L^2(\Omeps)^3}^2- 2\gamma\, \tau_\ep h_\ep\,  \|D \varphi\|_{L^2(\Omeps)^{3\times 3}} \|D \psi\|_{L^2(\Omeps)^{3\times 3}}\, ,\nonumber
\end{align}
where $\gamma$ is the defined by~\eqref{Def:gamma}.

Now,  by condition \eqref{Condition-gamma}, there exists $c_1>0$ satisfying
\begin{equation}\label{Proof-Coercivity-Def:c_1}
\frac{\gamma h}{R_c} < c_1 < \frac{1-2N^2}{\gamma h}\, ,
\end{equation} 
and by Young inequality,
\[
\|D\varphi\|_{L^2(\Omeps)^{3\times 3}} \| D\psi \|_{L^2(\Omeps)^{3\times 3}}  \leq \frac{c_1}{2\ep }\|D\varphi\|_{L^2(\Omeps)^{3\times 3}}^2 + \frac{\ep}{2 c_1}\|D\psi\|_{L^2(\Omeps)^{3\times 3}}^2\, .
\]
By continuity, there exists a real number $c_2$ satisfying $0<c_2<1$, and such that
\begin{equation}\label{Proof-Coercivity-Def:c_2}
c_1<\frac{1-\frac{2N^2}{c_2}}{\gamma h}\, ,
\end{equation}
and we also have
\[
\|D\varphi\|_{L^2(\Omeps)^{3\times 3}} \|\psi\|_{L^2(\Omeps)^3}  \leq \frac{\sqrt{2}}{4c_2}\|D\varphi\|_{L^2(\Omeps)^{3\times 3}}^2 + \frac{c_2}{\sqrt{2}}\|\psi\|_{L^2(\Omeps)^3}^2.
\]
Going back to estimate \eqref{CoercivityA-1}, we obtain
\begin{equation}\label{Ineq:CoercivityAeps}
\Ac_\ep((\varphi,\psi),(\varphi,\psi))  \geq A_\ep \|D\varphi\|_{L^2(\Omeps)^{3\times 3}}^2 + \ep^2 B_\ep \|D\psi\|_{L^2(\Omeps)^{3\times 3}}^2 + 4N^2(1-c_2)\|\psi\|_{L^2(\Omeps)^3}^2,
\end{equation}
where $A_\ep$, $B_\ep$ are defined by
\[
A_\ep = 1-\frac{2N^2}{c_2}- c_1 \gamma \tau_\ep \frac{h_\ep}{\ep},\qquad B_\ep = R_c - \frac{\gamma \tau_\ep} {c_1} \frac{h_\ep}{\ep}.
\]
In particular, there holds
\[
\lim_{\ep \rightarrow 0} A_\ep = 1-\frac{2N^2}{c_2}-c_1 \gamma h\, , \qquad
\lim_{\ep \rightarrow 0} B_\ep =R_c - \frac{\gamma}{c_1} h\, .
\]
Using conditions \eqref{Proof-Coercivity-Def:c_1}-\eqref{Proof-Coercivity-Def:c_2}, we conclude that for $\ep$ small enough, $\Ac_\ep$ is coercive.

\par\hfill$\square$

\subsection{A priori estimates and convergences}\label{Subsec-AprioriEstimates}
In this subsection we give a priori estimates and convergence results for the rescaled functions $\widetilde {\bf u}_\ep, \widetilde {\bf w}_\ep, \widetilde p_\ep$. Also, in order to take into account the effects of the rough boundary, we will introduce the unfolding method before proceeding with the proof of the Theorem \ref{thm_effective}.

\begin{prop}\label{prop_estim_tilde_velocity}
Assume that the asymptotic regimes \eqref{RM-R_c} and \eqref{Regime:Veps} and conditions (\ref{ExtraConditionOnN}) and (\ref{Condition-gamma}) hold. Then,  there exists a constant $C > 0$, independent on $\ep$, such that for any $0<\ep < \ep_0$,  the solution ($\ub_\ep, \wb_\ep, p_\ep)$ satisfies the estimates
\begin{align}
&\|{\bf u}_\ep\|_{L^2(\Omega_\ep)^3}\leq C\ep^{1\over 2},\quad \|D {\bf u}_\ep\|_{L^2(\Omega_\ep)^{3\times 3}}\leq C \ep^{-{1\over 2}},\label{vestim}\\
&\|{\bf w}_\ep\|_{L^2(\Omega_\ep)^3}\leq C \ep^{-{1\over 2}},\quad \|D {\bf w}_\ep\|_{L^2(\Omega_\ep)^{3\times 3}}\leq C \ep^{-{3\over 2}},\label{westim}\\
&\|p_\ep\|_{L^2(\Omega_\ep)}\leq C\ep^{-{3\over 2}},\quad\|\nabla p_\ep\|_{H^{-1}(\Omega_\ep)^3}\leq C\ep^{-{1\over 2}}\label{pestim}.
\end{align}
Moreover, the rescaled solution $(\widetilde {\bf u}_\ep,\widetilde {\bf w}_\ep,\widetilde p_\ep)$  defined by \eqref{changevar1_fun} satisfies the estimates:
\begin{align}
&\|\widetilde {\bf u}_\ep\|_{L^2(\widetilde \Omega_\ep)^3}\leq C,\quad \|\partial_{y_3}\widetilde {\bf u}_\ep\|_{L^2(\widetilde\Omega_\ep)^3}\leq C,\quad \|D_{y'}\widetilde {\bf u}_\ep\|_{L^2(\widetilde\Omega_\ep)^{3\times 2}}\leq C\ep^{-1}, \label{vtilestim}\\
& \|\widetilde {\bf w}_\ep\|_{L^2(\widetilde \Omega_\ep)^3}\leq C\ep^{-1},\quad \|\partial_{y_3}\widetilde {\bf w}_\ep\|_{L^2(\widetilde\Omega_\ep)^3}\leq C\ep^{-1},\quad \|D_{y'}\widetilde {\bf w}_\ep\|_{L^2(\widetilde\Omega_\ep)^{3\times 2}}\leq C\ep^{-2}, \label{wtilestim}\\
&\|\widetilde p_\ep\|_{L^2(\widetilde\Omega_\ep)}\leq C\ep^{-2},\quad \|\partial_{y_3}\widetilde p_\ep\|_{H^{-1}(\widetilde\Omega_\ep)}\leq C \ep^{-1},\quad \|\nabla_{y'}\widetilde p_\ep\|_{H^{-1}(\widetilde\Omega_\ep)^2}\leq C\ep^{-2}. \label{ptilestim}
\end{align}

\end{prop}
\medskip

%

\noindent {\bf Proof of Proposition \ref{prop_estim_tilde_velocity}.} According to Theorem \ref{thm_existence}, there exists $\ep_0>0$ such that for any $0<\ep < \ep_0$,  there exists a unique weak solution $(\ub_\ep,\wb_\ep,p_\ep)\in H^1(\Omeps)^3\times H^1(\Omeps)^3\times L^2_0(\Omeps)$ to system \eqref{1,1}--\eqref{bc-5}. Now, we obtain the estimates on the velocity and microrotation and then, we obtain the estimates for the pressure.

\paragraph{Estimates on $\ub_\ep$, $\wb_\ep$, $\widetilde \ub_\ep$ and $\widetilde \wb_\ep$} Using H\"older inequality, inequalities \eqref{Ineq-TraceOmeps} and \eqref{Ineq-PoincareOmeps}, and estimate~\eqref{Estimate:Sqrt2}, we obtain for any $(\varphi,\psi)\in \Veps\times\Veps^0$ :
\begin{align}
\Lc_\ep(\varphi,\psi)&  \leq \|D\Jb_\ep\|_{L^2(\Omeps)^{3\times 3}} \|D\varphi\|_{L^2(\Omeps)^{3\times 3}} + 2N^2 \|\nabla\times \Jb_\ep \|_{L^2(\Omeps)^{3}} h_\ep \|D\psi \|_{L^2(\Omeps)^{3\times 3}}\nonumber\\
 & \quad + 2N^2\beta\|\Jb_\ep\|_{L^2(\Geps^0)^3} \sqrt{\tau_\ep}\sqrt{h_\ep} \|D\psi\|_{L^2(\Omeps)^{3\times 3}}
\label{Ineq:ContinuityLeps} \\
& \leq \|D\Jb_\ep\|_{L^2(\Omeps)^{3\times 3}} \|D\varphi\|_{L^2(\Omeps)^{3\times 3}}
+ 2N^2\sqrt{h_\ep}\, C_\ep \|D\psi\|_{L^2(\Omeps)^{3\times 3}}\,, \nonumber
\end{align}
where $C_\ep$ is defined by $C_\ep = \sqrt{2}\sqrt{h_\ep} \|D\Jb_\ep\|_{L^2(\Omeps)^{3\times 3}} + \beta\sqrt{\tau_\ep}\|\Jb_\ep\|_{L^2(\Geps^0)^3} $. Notice that, by properties~\eqref{Estimate-lift-Jeps} satisfied by $\Jb_\ep$, there exists a constant $C>0$ such that $C_\ep\leq C$ for any $\ep>0$.

To obtain estimates \eqref{vestim} and \eqref{westim}, we test against $(\varphi,\psi) = (\vb_\ep,\wb_\ep)$ in~\eqref{MixedForm1} and use inequalities \eqref{Ineq:ContinuityLeps} and   \eqref{Ineq:CoercivityAeps}  to obtain
\begin{align*}
&  A_\ep \|D \vbeps\|_{L^2(\Omeps)^{3\times 3}}^2 + \ep^2 B_\ep \| D\wbeps\|_{L^2(\Omeps)^{3\times 3}}^2 \\
& \leq \|D\Jb_\ep\|_{L^2(\Omeps)^{3\times 3}} \|D\vbeps \|_{L^2(\Omeps)^{3\times 3}}
 + 2N^2\sqrt{h_\ep}\, C_\ep \|D\wbeps\|_{L^2(\Omeps)^{3\times 3}}\, .
\end{align*}
Last inequality can be written equivalently as
\begin{align*}
& A_\vps \left( \|D\vbeps\|_{L^2(\Omeps)^{3\times 3}} -\frac{1}{2A_\vps}\|D \Jb_\ep\|_{L^2(\Omeps)^{3\times 3}}  \right)^2 + \vps^2  B_\vps \left(\|D\wbeps\|_{L^2(\Omeps)^{3\times 3}} -\frac{N^2}{\ep^2B_\ep}\sqrt{h_\ep}C_\ep   \right)^2 \\
& \leq \frac{1}{4A_\ep}\|D \Jb_\ep\|_{L^2(\Omeps)^{3\times 3}}^2 + \frac{N^4}{\ep^2B_\ep} h_\ep C_\ep^2\, .
\end{align*}
Since $\tau_\ep, b_\ep/\ep, A_\ep$, $B_\ep$ and $C_\ep$ are uniformly bounded, using assumptions \eqref{Estimate-lift-Jeps}, we deduce than the right hand side of the previous inequality is bounded by $C/\ep$ for a certain constant $C>0$. This implies the following bounds:
\[
\|D\vbeps\|_{L^2(\Omeps)^{3\times 3}} \leq {C} \ep^{-1/2},\quad \|D\wbeps\|_{L^2(\Omeps)^{3\times 3}} \leq {C} \ep^{-3/2}.
\]
Hence, using Poincar\'e inequality \eqref{Ineq-PoincareOmeps}, the relation $\ub_\ep = \vb_\ep+\Jb_\ep$ and properties \eqref{Estimate-lift-Jeps} satisfied by $\Jb_\ep$, we obtain the desired estimates~\eqref{vestim}-\eqref{westim}. Finally, estimates~\eqref{vtilestim}-\eqref{wtilestim} are direct consequences of the rescaling~\eqref{changevar1_fun}.

\paragraph{Estimates on $p_\ep$ and $\widetilde{p}_\ep$}
In order to estimate $\nabla p_\ep$ in $H^{-1}(\Omeps)^3$, we test against $\varphi\in H^1_0(\Omeps)^3$ in~\eqref{WeakForm-u-Step1} :
\begin{align*}
\left\langle \nabla p_\ep,\varphi \right\rangle_{H^{-1}(\Omeps)^3\times H^1_0(\Omeps)^3 } & = -\int_{\Omeps} p_\ep\, \dv \varphi \,dx \\
&  = -\int_{\Omeps} D\ub_\eps  : D\varphi  \,dx    + 2N^2 \int_{\Omeps}\wb_\ep\cdot (\nabla\times \varphi)   \,dx \, .
\end{align*}
Using H\"older inequality and estimate~\eqref{Estimate:Sqrt2}, we deduce
\[
\big| \left\langle \nabla p_\ep,\varphi \right\rangle_{H^{-1}(\Omeps)^3\times H^1_0(\Omeps)^3 }\big| \leq \big( \|D\ubeps\|_{L^2(\Omeps)^{3\times 3}} + 2N^2\sqrt{2} \|\wbeps\|_{L^2(\Omeps)^{3\times 3}} \big)\|D\varphi\|_{L^2(\Omeps)^{3\times 3}}\, .
\]
We conclude from the upper bounds~\eqref{vestim}-\eqref{westim} that  $\|\nabla p_\ep\|_{H^{-1}(\Omeps)^3}\leq C\ep^{-1/2}$, where $C$ depends only on $N$.

Finally, to estimate $p_\ep$ in $L^2_0(\Omeps)$, we apply the following inequality, whose proof is given in~\cite[Corollary 4.2]{CLS-SIMA}:
\[
\|p_\ep\|_{L^2_0(\Omeps)} \leq \frac{C}{\ep} \|\nabla p_\ep\|_{H^{-1}(\Omeps)^3}\, .
\]
This proves~\eqref{pestim}. Finally, estimates~\eqref{ptilestim} are direct consequences of the rescaling~\eqref{changevar1_fun},
which concludes the proof of Proposition~\ref{prop_estim_tilde_velocity}.
\qed

\medskip

As a consequence of the a priori estimates stated in Proposition \ref{prop_estim_tilde_velocity}, and the fact that $\Omega\subset \widetilde\Omega_\ep$ and $|\widetilde \Omega_\ep\setminus\Omega|\to 0$, we have the following convergences for rescaled solutions $\widetilde \ub_\ep$, $\widetilde \wb_\ep$ restricted to the limit domain $\Omega$.

\begin{lm}\label{lema_tilde_v} 
Assume that the asymptotic regimes \eqref{RM-R_c} and \eqref{Regime:Veps} and conditions (\ref{ExtraConditionOnN}) and (\ref{Condition-gamma}) hold. Then, for a subsequence of $\ep$, still denoted by $\ep$, there exist $\widetilde {\bf u}', \widetilde {\bf w}'\in H^1(0,h;L^2(\omega))^2$
with $\widetilde {\bf u}'(x',h)=\widetilde {\bf w}'(x',h)=0$ for a.e.\! $ x'\in \omega$, and 
\begin{equation}\label{div_cond_u}\dis \mathrm{div}_{y'}\int_0^{h} \widetilde {\bf u}'(y',y_3)\,dy_3=S\hbox{ in }  H^{-1}(\omega),\end{equation}
 such that
\begin{align}
\widetilde \ub_{\ep|\Omega}&\rightharpoonup (\widetilde {\bf u}',0) \quad  \hbox{ in }H^1(0,h;L^2(\omega))^3, \label{convergence_u}\\ 
\ep\widetilde {\bf w}_{\ep|\Omega}&\rightharpoonup   (\widetilde {\bf w}',0)\quad \hbox{ in }H^1(0,h;L^2(\omega))^3. \label{conv_w_ep}
\end{align} 
%
\end{lm}

\noindent {\bf Proof. } The space $H^1(0,h;L^2(\omega))$ is a Hilbert space for the norm 
\[\|v\|_{H^1(0,h;L^2(\omega))}=\big( \|v\|_{L^2(\Omega)}^2 + \|\partial_{y_3}v\|_{L^2(\Omega)}^2  \big)^{1/2}\, .\]
By estimates~\eqref{vtilestim}-\eqref{wtilestim}, $\widetilde {\bf u}_{\ep|\Omega}$ and $\ep\widetilde {\bf w}_{\ep|\Omega}$ are bounded in $H^1(0,h;L^2(\omega))^3$, so there exist $\widetilde {\bf u}$ and $\widetilde {\bf w}$  such that for a subsequence of $\ep$, still denoted by $\ep$, we have
\begin{equation}\label{WeakConv-tilde_u-tilde_w}
\widetilde {\bf u}_{\ep|\Omega} \rightharpoonup \widetilde {\bf u}\quad \textrm{and}\quad \ep \widetilde {\bf w}_{\ep|\Omega} \rightharpoonup \widetilde {\bf w}\, \quad \textrm{in }  H^1(0,h;L^2(\omega))^3\, .
\end{equation}
By continuity of the trace operator from $H^1(0,h;L^2(\omega))$ into $L^2(\omega\times \{h\})$, the conditions   $\widetilde \ub_\ep(x',h)=-\ep Se_3$ and $\ep\widetilde \wb_\ep(x',h)=0$ for a.e.\! $x'\in \omega$ pass to the limit, yielding $\widetilde \ub(x',h)=\widetilde \wb(x',h)=0 $ for a.e.\! $x'\in \omega$.

Now, we prove that $\widetilde u_3 = 0$. Since $\ub_\ep$ is divergence free, using definition~\eqref{changevar1_fun}, the rescaled function $\widetilde{\ub}_\ep$ satisfies
\begin{equation}\label{div_tilde}
\mathrm{div}_{y'}\,  \widetilde{\ub}'_\ep + \frac{1}{\ep}\partial_{y_3}\widetilde{u}_{\ep,3}=0\quad  \textrm{a.e.\! in }\widetilde \Omega_\ep,
\end{equation}
so for any $\phi\in C_c^{\infty}(\Omega)$,
\begin{align*}
0 & = \int_{\Omega} \mathrm{div}_{y'}\,  \widetilde{\ub}'_\ep\, \phi  \,dy + \int_{\Omega} \frac{1}{\ep}\partial_{y_3}\widetilde{u}_{\ep,3}\, \phi \,dy 
 = -  \int_{\Omega} \widetilde{\ub}'_\ep  \cdot \nabla_{y'} \phi \,dy + \int_{\Omega} \frac{1}{\ep}\partial_{y_3}\widetilde{u}_{\ep,3}\, \phi  \,dy \, .
\end{align*}
Hence, $\int_{\Omega} \partial_{y_3}\widetilde{u}_{\ep,3}\, \phi \,dx  = \ep  \int_{\Omega} \widetilde{\ub}'_\ep  \cdot \nabla_{y'} \phi  \,dx $ and using~\eqref{WeakConv-tilde_u-tilde_w}, we deduce that $\int_{\Omega} \partial_{y_3}\widetilde{u}_3\, \phi  \,dx = 0$. As a result, $\widetilde{u}_3$ does not depend on $y_3$, and since it vanishes on $y_3=h$, it is identically null.

Next, we prove the divergence equation (\ref{div_cond_u}).  Using condition~\eqref{21,1}, integration by parts, boundary conditions~\eqref{bc-1},~\eqref{bc-3} and the change of variables~\eqref{changevar1}, we have for any $\phi\in C_c^\infty(\omega)$
\begin{align*}
	0&=\int_{\Omeps} (\mbox{div}\, \ub_\ep)\phi(x')\, dx\\
	& = -\int_{\Omeps} \ub_\ep'\cdot \nabla_{x'}\phi  \, dx +\int_{\G_\ep^1}u_{\ep,3}\, \phi\, d\sigma\\
	 & = - \eps\int_{\widetilde \Omega_\ep} \widetilde{\ub}_{\ep}'\cdot \nabla_{y'}\phi(y')\, dy - \eps S\int_{\omega}\phi(y')\, dy'\,.
	\end{align*}
Noticing that, by the bound~\eqref{vtilestim} and H\"older inequality, there holds
$\lim_{\ep\rightarrow 0}\int_{\widetilde{\Omega}_\eps \setminus \Omega} | \widetilde{\ub}_{\ep}|^2dy = 0$, we deduce that
\[
-\int_{\Omega} \widetilde{\ub}_{\ep}'\cdot \nabla_{y'}\phi(y')\, dy =  S\int_{\omega}\phi(y')\, dy' + O_\ep \, .
\]
Using the weak convergence~\eqref{WeakConv-tilde_u-tilde_w}, we can pass to the limit in the previous equality and obtain
\begin{align*}
	S\int_\omega\phi(y')\,dy'& = -\int_{\Omega} \widetilde{\ub}'  \cdot \nabla_{y'} \phi(y') \,dy \\
	& = \int_{\omega} {\rm div}_{y'}\left(\int_0^h\widetilde{\ub}'(y',y_3)\, dy_3\right)\, \phi(y') \,dy'\, ,
	\end{align*}
which proves  \eqref{div_cond_u}.

Finally, it remains to prove that $\widetilde w_3=0$. To do this, for any $\psi\in C_c^\infty(\Omega)$, we consider  $\psi_\ep=\ep \psi(x',x_3/\ep) e_3$ as test function in the variational formulation (\ref{WeakForm-u-Step2}). 
Applying the change of variables (\ref{changevar1}) and extending the integrals to $\Omega$,  we get
$$ R_c\int_{\Omega}\ep\partial_{y_3}\widetilde w_{\ep,3} \partial_{y_3}\psi_3\,dy+4N^2\int_{\Omega}\ep \widetilde w_{\ep,3}\psi\,dy=2N^2\int_{\Omega}\ep{\rm Rot}_{x'}\widetilde \ub'_\ep\psi_3\,dy + O_\ep.
$$
Integrating by parts the right-hand side, we get
$$ R_c\int_{\Omega}\ep\partial_{y_3}\widetilde w_{\ep,3} \partial_{y_3}\psi_3\,dy+4N^2\int_{\Omega}\ep \widetilde w_{\ep,3}\psi\,dy=2N^2\int_{\Omega}\ep[\widetilde \ub'_\ep]^\perp \nabla_{y'}\psi_3\,dy + O_\ep.
$$
Using convergences \eqref{WeakConv-tilde_u-tilde_w}, when $\ep$ tends to zero, we get
\begin{equation}\label{w3limit}
R_c\int_{\Omega}\partial_{y_3}\widetilde w_{3} \, \partial_{y_3}\psi\,dy+4N^2\int_{\Omega}  \widetilde w_{3}\psi\,dy=0.
\end{equation}
Next, we prove that $\widetilde w_3(y',0)=0$ for a.e.\! $x'\in \omega$. The condition ${\bf w}_\varepsilon\cdot{\bf n}_\varepsilon=0$ on $\Gamma_\varepsilon^0$ can be rewritten as follows
$$\varepsilon \tilde w_{\varepsilon,1}(y',-\tilde \Psi_\varepsilon(y'))\lambda \varepsilon^{\delta-\ell}\partial_{1}\Psi({1\over \varepsilon^\ell}y'\cdot e'_1) + \varepsilon \tilde w_{\varepsilon,3}(y',-\tilde \Psi_\varepsilon(y'))=0\quad\hbox{a.e.}\ y'\in\omega.$$
Multiplying this equality by $\psi\in C^\infty_c(\omega)$ and integrating on $\omega$, we get
\begin{equation}\label{w3zero}
\int_{\omega}\varepsilon \tilde w_{\varepsilon,1}(y',-\tilde \Psi_\varepsilon(y'))\lambda e^{\delta-\ell}\partial_{1}\Psi({1\over \varepsilon^\ell}y'\cdot e'_1)\psi(y')dy' +\int_{\omega} \varepsilon \tilde w_{\varepsilon,3}(y',-\tilde \Psi_\varepsilon(y'))\psi(y')dy'=0.
\end{equation}
We can write the second term of (\ref{w3zero}) as follows
$$\int_{\omega} \varepsilon \tilde w_{\varepsilon,3}(y',-\tilde \Psi_\varepsilon(y'))\psi(y')dy'=\int_\omega 
\varepsilon \tilde w_{\varepsilon,3}(y',0)\psi(y')dy'-\int_\omega\left(\int_{-\tilde\Psi_\varepsilon(y')}^0\varepsilon\partial_{y_3}\tilde w_{\varepsilon,3}(y',s)ds\right)\psi(y')dy'.$$
Then, since $\varepsilon\tilde{\bf w}_\varepsilon|_{\Omega}$ is  bounded in $H^1(0,h;L^2(\omega))^3$, by continuity of the trace operator from $H^1(0,h;L^2(\omega))^3$ into $L^2(\omega\times \{0\})^3$ and  convergence (\ref{WeakConv-tilde_u-tilde_w}), we have that
$$\int_\omega 
\varepsilon \tilde w_{\varepsilon,3}(y',0)dy'=\int_\omega 
\varepsilon \tilde w_{ 3}(y',0)dy'+O_\varepsilon.$$
Moreover, from the {\color{black} Cauchy-Schwarz} inequality, estimate (\ref{wtilestim})$_2$ and $|\widetilde\Omega_\varepsilon^{-}|\to 0 $ (recall that $\|\widetilde\Psi_\varepsilon\|_{L^\infty}\leq C\varepsilon^{\delta-1}$), then
$$\left|\int_\omega\int_{-\tilde\Psi_{\varepsilon}(y')}^0\varepsilon\partial_{y_3}(y',s)ds\varphi(y')dy'\right|\leq \left(\int_{\widetilde\Omega_\varepsilon}|\varepsilon\partial_{y_3}\tilde w_{\varepsilon,3}|^2\right)^{1\over 2}
\left(\int_{\widetilde\Omega_\varepsilon^-}|\varphi(y')|^2\right)^{1\over 2}=O_\varepsilon,$$
and so, we get that
$$\int_{\omega} \varepsilon \tilde w_{\varepsilon,3}(y',-\tilde \Psi_\varepsilon(y'))dy'=\int_\omega 
\varepsilon \tilde w_{3}(y',0)dy'+O_\varepsilon.$$
A similar argument works for the first term of (\ref{w3zero}) works replacing $\varphi(y')$ by $\lambda \varepsilon^{\delta-\ell}\partial_1\Psi({1\over \varepsilon^\ell}y'\cdot e'_1)\varphi(y')$, which goes to $0$ in $L^\infty(\omega)$. \\

Then, from the above, passing to the limit in (\ref{w3zero}), we get 
$$\int_{\omega}\widetilde w_3(y',0)\varphi(y')dy'=0,$$
which is equivalent to $\widetilde w_3(y',0)=0$ for a.e.\! $y'\in \omega$.\\

Finally, from (\ref{w3limit})  and taking into account that $\widetilde w_3(y',h) =\widetilde w_3(y',0)= 0$ for a.e.\! $y'\in \omega$, it is easily deduced that $\widetilde w_3 = 0$, which ends the proof.
\qed 

\medskip

In order to give the convergence of the rescaled pressure $\widetilde p_\ep$, let us give a more accurate estimate for pressure $p_\ep$. For this, we need to recall a decomposition result for $p_\ep$ whose proof can be found in  \cite[Corollary 4.2]{CLS-SIMA}.

\begin{prop}\label{Prop:PressureDecomposition} The following decomposition for $p_\varepsilon\in L^2_0(\Omega_\ep)$ holds
\begin{equation}\label{decomposition} p_\ep=p_\ep^0+ p_\ep^1,\end{equation} where $p_\ep^0\in H^1(\omega)$, which is independent of  $x_3$, and $p_\ep^1\in L^2(\Omega_\ep)$. Moreover, the following estimates hold 
\begin{equation}\label{decomposition-estim}
\|p_\ep^0\|_{H^1(\omega)}\leq C\ep^{-{3\over 2}}\|\nabla p_\ep\|_{H^{-1}(\Omega_\ep)^3},\qquad
\|p_\ep^1\|_{L^2(\Omega_\ep)}\leq C\|\nabla p_\ep\|_{H^{-1}(\Omega_\ep)^3}.
\end{equation}
\end{prop}
From this result, we are able to give the convergence result for $\widetilde p_\ep$. We denote by $\widetilde p_{\ep}^1$ the rescaled function associated with $p_\ep^1$, defined by $\widetilde p_\ep^1(y)=p_\ep^1(y',\ep y_3)$ for a.e. $y\in \widetilde{\Omega}_\ep$.
\begin{cor}  \label{lm_conv_p}Previous result implies the existence of $ p\in H^1(\omega)$  and $\widetilde p^1\in L^2(\Omega)$ satisfying
\begin{equation}\label{conv_p0_p1}
\ep^2  p_\ep^0\rightharpoonup  p\hbox{  in }H^1(\omega),\qquad \ep\widetilde p_{\ep|\Omega}^1\rightharpoonup \widetilde p^1\hbox{  in }L^2(\Omega),
\end{equation}
and moreover 
\begin{equation}\label{converg_p}\ep^2\widetilde p_{\ep|\Omega} \to p\hbox{  in }L^2(\Omega).
\end{equation}
\end{cor}

\noindent {\bf Proof. } From (\ref{decomposition-estim}) and (\ref{pestim})$_2$, we get 
\begin{equation}\label{estimates_p0_p1}\|p_\ep^0\|_{H^1(\omega)}\leq C\ep^{-2},\quad \|p_\ep^1\|_{L^2(\Omega_\ep)}\leq C\ep^{-{1\over 2}},
\end{equation}
and after rescaling $p_\ep^1$, last inequality  becomes $$\|\widetilde p_\ep^1\|_{L^2(\widetilde \Omega_\ep)}\leq C\ep^{-1}.$$
Previous estimates and the fact that $\Omega\subset \widetilde\Omega_\ep$ and $|\widetilde \Omega_\ep\setminus\Omega|\to 0$ imply (\ref{conv_p0_p1}). The strong convergence (\ref{converg_p}) for the complete pressure $\widetilde p_\ep$ is a direct consequence of (\ref{conv_p0_p1}) and the decomposition (\ref{decomposition}).
\par \hfill $\square$

\subsection{Unfolding method}
In order to capture the behaviour of ${\bf u}_\ep$, ${\bf w}_\ep$ and $p_\ep^1$ (introduced in Proposition~\ref{Prop:PressureDecomposition}) near the rough boundary  $\Gamma_\ep$, we need to introduce a new change of variables, which is adapted from the unfolding method (see \cite{ArDoHo,CLS-SIMA,CDG}). To do this, for ${\bf u}_\ep, {\bf w}_\ep\in H^1(\Omega_\ep)^3$ satisfying boundary conditions~\eqref{bc-1}--\eqref{bc-3},  $p_\ep^1\in L^2(\Omega_\ep)$,  and $\rho>0$,  we set $\widehat {\bf u}_\ep$, $\widehat {\bf w}_\ep$ and $\widehat p_\ep^1$ by 
\begin{equation}\label{uhat}
\widehat {\bf u}_\ep(x',z)={\bf u}_\ep\left(\ep^\ell\kappa\left({x'\over \ep^\ell}\right)+\ep^\ell z',\ep^\ell z_3\right),
\end{equation}
\begin{equation}\label{what}
\widehat {\bf w}_\ep(x',z)={\bf w}_\ep\left(\ep^\ell\kappa\left({x'\over \ep^\ell}\right)+\ep^\ell z',\ep^\ell z_3\right),\end{equation}
\begin{equation}\label{defp1}
\widehat p_\ep^1(x',z)=p_\ep^1\left(\ep^\ell\kappa\left({x'\over \ep^\ell}\right)+\ep^\ell z', \ep^\ell z_3\right),
\end{equation}
for a.e. $(x',z)\in\omega_\rho\times \widehat Z_\ep$, where $\omega_\rho$ is defined by~\eqref{Def:omega_rho} and 
$$\widehat Z_\ep=\left\{z\in Z'\times\mathbb{R}\,:\, -\ep^{\delta-\ell}\Psi(z'\cdot e_1')<z_3<\ep^{1-\ell}h\right\}.$$
\begin{rem}
For  every $k'\in I_{\rho,\eps}$, the functions $\widehat \ub_\ep$, $\widehat \wb_\ep$ and $\widehat p_\ep^1$ restricted to $C_{\ep^\ell}^{k'}\times \widehat Z_\ep$ are independent of  $x'$. However,  as functions depending  on $z$, they are obtained from their original functions by means of 
\begin{equation}\label{changevar2}
z'={x'-\ep^\ell k'\over \ep^\ell},\quad z_3={x_3\over \ep^\ell},
\end{equation}
that converts $Q_{\ep^\ell}^{k'}$ in $\widehat Z_\ep$. 
\end{rem}

Thanks to the estimates satisfied by ${\bf u}_\ep$, ${\bf w}_\ep$ and $p_\ep^1$ given in (\ref{vestim}), (\ref{pestim}) and (\ref{estimates_p0_p1})$_2$, respectively, we have the following compactness results. 

\begin{lm} \label{lem_conv_hat_uv} Consider two sequences ${\bf u}_\ep, {\bf w}_\ep\in \Vb_\vps$ satisfying (\ref{vestim})$_2$ and (\ref{westim})$_2$, respectively. Define $\widetilde {\bf u}_\ep, \widetilde {\bf w}_\ep\in H^1(\widetilde \Omega_\ep)^3$ by (\ref{changevar1_fun}), so that (\ref{convergence_u}) and (\ref{conv_w_ep}) hold. Set $\delta\leq {3\over 2}\ell-{1\over 2}$. Then,
\\

\begin{itemize}
\item[(i)] If $\delta< {3\over 2}\ell-{1\over 2}$, it holds 
\begin{equation}\label{unablapsi}
\widetilde u_1(x',0) \partial_{y_1}\Psi(z'\cdot e_1')=0,\quad\hbox{ a.e. }(x',z')\in\omega\times Z',
\end{equation}
\begin{equation}\label{wnablapsi}
\widetilde w_1(x',0) \partial_{y_1}\Psi(z'\cdot e_1')=0,\quad\hbox{ a.e. }(x',z')\in\omega\times Z'.
\end{equation}

\item[(ii)] If $\delta= {3\over 2}\ell-{1\over 2}$, there exist $\widehat {\bf u}, \widehat {\bf w}\in L^2(\omega,\mathcal{V}^3)$, {\color{black} where $\mathcal{V}$ is the space of functions $\widehat \varphi:
\RR^2\times(0,+\infty)\mapsto \RR$ such that $\widehat \varphi\in H^1_{\#}(\widehat Q_M)$, for every $M>0$, and $\nabla
\widehat \varphi\in L^2_{\#}(\widehat Q)^3$}, satisfying
\begin{equation}\label{propertyhatu3}\widehat u_3(x',z',0)=-\lambda \partial_{y_1}\Psi(z'\cdot e_1')\widetilde u_1(x',0), \quad\hbox{a.e. }(x',z')\in\omega\times Z',\end{equation}
\begin{equation}\label{propertyhatw3}\widehat w_3(x',z',0)=-\lambda \partial_{y_1}\Psi(z'\cdot e_1')\widetilde w_1(x',0), \quad\hbox{a.e. }(x',z')\in\omega\times Z',\end{equation}
and such that, for any $\rho, M>0$, the sequences $\widehat {\bf u}_\ep$ and $\widehat {\bf w}_\ep$, respectively given by (\ref{uhat}) and (\ref{what}), satisfy
\begin{equation}\label{Dhatuw}
\ep^{1-\ell\over 2}D_z\widehat {\bf u}_\ep\rightharpoonup D_z \widehat {\bf u}\  \hbox{in }L^2(\omega_\rho\times \widehat Q_M)^{3\times 3},\quad \ep^{3-\ell\over 2}D_z\widehat {\bf w}_\ep\rightharpoonup D_z \widehat {\bf w}\  \hbox{in }L^2(\omega_\rho\times \widehat Q_M)^{3\times 3}.
\end{equation}

Moreover, if one assumes ${\rm div}\, {\bf u}_\ep=0$ in $\Omega_\ep$, then $\widehat {\bf u}$ satisfies
\begin{equation}\label{divhatu}{\rm div}_z \widehat {\bf u}=0\ \hbox{ in }\ \omega\times \widehat Q.
\end{equation}
\end{itemize}
\end{lm}

{\bf Proof. } This result is a direct consequence of \cite[Lemma 5.4]{CLS-SIMA}, applied to the sequences $\ep^2\ub_\ep$ and $\ep^3\wb_\ep$. 
\par\hfill $\square$

\begin{lm} 
We consider $p_\ep^1\in L^2(\Omega_\ep)$ such that (\ref{estimates_p0_p1})$_2$ holds. 
Then, there exists $\widehat p^1\in L^2(\omega\times \widehat Q)$ satisfying, up to a subsequence, the convergence
\begin{equation}\label{convp1}
\ep^{1+\ell \over 2}\widehat p_\ep^1\rightharpoonup \widehat p^1\ \hbox{ in }\ L^2\left(\omega_\rho\times \widehat Q_M\right)\quad\forall \rho, M>0.
\end{equation}
\end{lm}
%
{\bf Proof. } It is a direct consequence of~\cite[Lemma 5.5]{CLS-SIMA}, applied to the sequence $\ep^2 p^1_\ep$. 
\par\hfill$\square$

We are now in position to prove Theorem~\ref{thm_effective} in the critical case $\delta= {3\over 2}\ell-{1\over 2}$, which is the most relevant from the mechanical point of view since it describes the coupling effects between the riblets and nonzero boundary conditions.\\

\noindent {\bf Proof of Theorem \ref{thm_effective}. }  Let us consider $\delta= {3\over 2}\ell-{1\over 2}$ with $\ell>1$.

First of all,   Lemma \ref{lema_tilde_v} and Corollary \ref{lm_conv_p} implies the existence of $\widetilde {\bf u}', \widetilde {\bf w}'\in  H^1(0,h;L^2(\omega))^2$ such that (\ref{bc_top_system_1}), and $ p\in H^1(\omega)$  so that convergences of $\widetilde \ub_\ep$, $\widetilde \wb_\ep$ and $\widetilde p_\ep$ given in (\ref{convergencias_solutions}) hold. Also, we have that the divergence condition (\ref{limit_system_1})$_3$ holds. From Corollary \ref{lm_conv_p}, the sequences $p_\ep^0$ and $\tilde p_\ep^1$ satisfy convergences given in (\ref{conv_p0_p1}).

%
%
%

We recall the variational formulation given by (\ref{WeakForm-u-Step1}) and (\ref{WeakForm-u-Step2}). For $\varphi, \psi \in  \Veps$, $({\bf u}_\ep,{\bf w}_\ep,p_\ep)$ satisfies
\begin{equation}\label{vf_vel_1}
\begin{array}{l}
\dis \int_{\Omeps} D\ub_\eps  : D\varphi \,dx -\int_{\Omeps} p_\ep\, \dv \varphi  \,dx - 2N^2 \int_{\Omeps}\wb_\ep\cdot (\nabla\times \varphi) \,dx\\
\noame
\dis -2\left({1\over\alpha}-N^2\right) \int_{\Geps^0}(\wb_\ep\times {\bf n}_\ep)\cdot \varphi  \,d\sigma = 0\,,
\end{array}
\end{equation}
\begin{equation}\label{vf_mic_1}
\begin{array}{l}
\dis
\ep^2R_c\int_{\Omeps} D\wb_\ep : D\psi  \,dx - 2N^2\beta \int_{\Geps^0}(\ub_\ep\times {\bf n}_\ep)\cdot \psi  \,d\sigma + 4N^2\int_{\Omeps} \wb_\eps\cdot \psi  \,dx  \\
\noame
\dis - 2N^2\int_{\Omeps}(\nabla \times \ub_\eps)\cdot \psi  \,dx = 0 \, .
\end{array}
\end{equation}

%
%

Now, we want to pass to the limit in the above variational formulations. To do this, we will use  appropriate test functions $\varphi_\ep$, $\psi_\ep$. We divide the proof in four steps. \\

{\it Step 1.  Definition of the test functions. }  Lemma \ref{lem_conv_hat_uv} gives the existence of $\widehat {\bf u},\,\widehat {\bf w}\in L^2(\omega;\mathcal{V}^3)$ satisfying (\ref{propertyhatu3}), (\ref{propertyhatw3}) and (\ref{divhatu}). Thanks to this, we consider the following test functions. For any $\widetilde \varphi, \widetilde\psi \in C_c^1(\omega\times(-h,h))^3$, with $\varphi_3=\psi_3=0$, $\widehat \varphi,\widehat \psi\in C_c^1(\omega;C_{\#}^1(\widehat Q))^3$ satisfying

\begin{equation}\label{testfunc}
\left\{\begin{array}{l}
D_z\widehat \varphi(x',z)=0\hbox{ a.e. in }\{z_3>M\}\hbox{ for some  }M>0,\\
\noame\dis
\widetilde\varphi'(y',y_3)=\widetilde \varphi'(y',0) \hbox{ when  }y_3\leq 0,\\
\noame\dis
\widehat\varphi(x',z',z_3)=\widehat \varphi(x',z',0) \hbox{ when   }z_3\leq 0,\\
\noame\dis
\lambda\, \partial_{z_1}\Psi(z'\cdot e_1')\widetilde\varphi_1(y',0)+\widehat \varphi_3(y',z',0)=0,
\end{array}\right.
\end{equation}
\begin{equation}\label{testfunc_psi}
\left\{\begin{array}{l}
D_z\widehat \psi(x',z)=0\hbox{ a.e. in }\{z_3>M\}\hbox{ for some  }M>0,\\
\noame\dis
\widetilde\psi'(y',y_3)=\widetilde \psi'(y',0) \hbox{ when }y_3\leq 0,\\
\noame\dis
\widehat\psi(x',z',z_3)=\widehat \psi(x',z',0) \hbox{ when }z_3\leq 0,\\
\noame\dis
 \lambda\,\partial_{z_1}\Psi(z'\cdot e_1')\widetilde\psi_1(y',0)+\widehat \psi_3(y',z',0)=0,
\end{array}\right.
\end{equation}
and a function $\zeta\in C^\infty(\mathbb{R})$ such that
\begin{equation}\label{zeta}
\zeta(s)=\left\{\begin{array}{l} \dis1\hbox{ when }s<{1\over 3},\\
\noame\dis
0\hbox{ when }s>{2\over 3},
\end{array}\right.
\end{equation}
we set $\varphi_\ep,\psi_\ep\in H^1(\Omega_\ep)^3$ as follows
$$\left\{\begin{array}{l}\dis
\varphi_\ep'(x)={\ep}\widetilde\varphi'\left(x',{x_3\over \ep}\right)+\ep^{1+\ell\over 2}\widehat \varphi'\left(x',{x\over \ep^\ell}\right)\zeta\left({x_3\over \ep}\right),\\
\noame\dis
\varphi_{\ep,3}=\ep^{1+\ell\over 2}\widehat \varphi_3\left(x',{x\over \ep^\ell}\right)\zeta\left({x_3\over \ep}\right) -\ep^{\ell}    \widehat \varphi_1\left(x',{x\over \ep^\ell}\right)  \lambda\,  \partial_{x_1}\Psi\left({1\over \ep^\ell}x'\cdot e_1'\right)\zeta\left({x_3\over \ep^\ell}\right),
\end{array}\right.$$

$$\left\{\begin{array}{l}\dis
\psi_\ep'(x)=\widetilde\psi'\left(x',{x_3\over \ep}\right)+\ep^{\ell-1\over 2}\widehat \psi'\left(x',{x\over \ep^\ell}\right)\zeta\left({x_3\over \ep}\right),\\
\noame\dis
\psi_{\ep,3}=\ep^{\ell-1\over 2}\widehat \psi_3\left(x',{x\over \ep^\ell}\right)\zeta\left({x_3\over \ep}\right) -\ep^{\ell-1}  \widehat \psi_1\left(x',{x\over \ep^\ell}\right) \lambda\, \partial_{x_1}\Psi\left({1\over \ep^\ell}x'\cdot e_1'\right)\zeta\left({x_3\over \ep^\ell}\right).
\end{array}\right.$$
%
Since $\widetilde\varphi'(x)$, $\widetilde\psi'(x)$, $\widehat\varphi'(x',z)$ and $\widehat\psi(x',z)$ are zero when $x'$  is out of a compact subset of $\omega$,  (\ref{testfunc}) and (\ref{testfunc_psi}), then  $\varphi_\ep$, $\psi_\ep$ are such that 
$$\varphi_\ep=\psi_\ep=0\hbox{ on }\partial\Omega_\ep\setminus\Gamma_\ep^0,\quad \varphi_\ep\cdot {\bf n}_\ep=\psi_\ep\cdot {\bf n}_\ep=0\hbox{ on }\Gamma_\ep^0.$$
So, we are able to consider $\varphi_\ep$ and $\psi_\ep$, respectively, as test functions in (\ref{vf_vel_1}) and (\ref{vf_mic_1}).  The difficulty now is to obtain the limit of every terms of (\ref{vf_vel_1}) and (\ref{vf_mic_1}). For this, we observe that from the conditions $D_z\widehat \varphi=D_z\widehat \psi=0$ a.e. in $\{z_3>M\}$ and (\ref{zeta}), it follows

\begin{equation}\label{phieptest}
\varphi_\ep(x)=\ep\left(\widetilde\varphi'\left(x',{x_3\over \ep}\right),0\right) + g_\ep\quad\hbox{in } \Omega_\ep,
\end{equation}
\begin{equation}\label{Dphieptest}
D\varphi_\ep(x)=\sum_{i=1}^2\partial_{y_3}\widetilde\varphi_i\left(x',{x_3\over \ep}\right)e_i\otimes e_3+ \ep^{1-\ell\over 2}D_z\widehat \varphi\left(x',{x\over \ep^\ell}\right) + h_\ep(x)\quad\hbox{in }\Omega_\ep,
\end{equation}
\begin{equation}\label{psieptest}
\psi_\ep(x)=\left(\widetilde\psi'\left(x',{x_3\over \ep}\right),0\right) + \breve g_\ep\quad\hbox{in }\overline\Omega_\ep,
\end{equation}
\begin{equation}\label{Dpsieptest}
D\psi_\ep(x)={ \ep^{-1}}\sum_{i=1}^2\partial_{y_3}\widetilde\psi_i\left(x',{x_3\over \ep}\right)e_i\otimes e_3+ \ep^{-{1+\ell\over 2}}D_z\widehat \psi\left(x',{x\over r_\ep}\right) + \breve h_\ep(x)\quad\hbox{in }\Omega_\ep,
\end{equation}
where $g_\ep, \breve g_\ep\in C^0(\overline\Omega_\ep)^3$, $h_\ep, \breve h_\ep\in C^0(\overline\Omega_\ep)^{3\times 3}$ (thanks to $\ell>1$) are such that

\begin{equation}\label{gbounded1}
\eps^{-3}  \int_{\Omega_\ep}|g_\ep|^2dx\leq C \left(\ep^{\ell+1} + \ep^{3(\ell-1)}\right)   =O_\ep,
\end{equation}

\begin{equation}\label{gbounded2}
\ep^{-2}\int_{\Gamma_\ep^0}|g_\ep|^2\,d\sigma \leq C\ep^{\ell-1}=O_\ep,
\end{equation}

\begin{equation}\label{Dgbounded}
\ep^{-1}\int_{\Omega_\ep}|h_\ep|^2dx\leq C\ep^{3}\left(\ep^{\ell-2} +\ep^{\ell-4}+\frac{1}{\ep} \right)=O_\ep,
\end{equation}

\begin{equation}\label{g_2bounded1}
\ep^{-1}\, \int_{\Omega_\ep}|\breve g_\ep|^2dx\leq C \left(\ep^{\ell-1} + \ep^{3(\ell-1) }\right)   =O_\ep,
\end{equation}

\begin{equation}\label{g_2bounded2}
\int_{\Gamma_\ep^0}|\breve g_\ep|^2\,d\sigma \leq C\ep^{\ell-1}  =O_\ep,
\end{equation}

\begin{equation}\label{Dg_2bounded}
	\ep\int_{\Omega_\ep}|\breve h_\ep|^2dx\leq C \ep^{3}\left(\ep^{\ell-2} +\ep^{\ell-4}+\frac{1}{\ep} \right)=O_\ep.
\end{equation}

{\color{black} We remark that functions $g_\ep, \breve g_\ep, h_\ep$ and $\breve h_\ep$ and previous estimates are devoted to identify terms of the variational formulation that are negligible in the asymptotic analysis.}\\


{\it Step 2. Passing to the limit in variational formulation (\ref{vf_vel_1}).} We can pass to the limit in every term of (\ref{vf_vel_1}).\\

\noindent $\bullet$  {\it 1st term of (\ref{vf_vel_1})}: From (\ref{vestim}), (\ref{Dphieptest}) and (\ref{Dgbounded}), we deduce
\begin{equation}\label{F1-v}
\begin{array}{ll}\dis
\int_{\Omeps} D\ub_\eps  : D\varphi_\ep \,dx=&\dis \int_{\Omega_\ep^+}\partial_{x_3} {\bf u}_\ep'(x)\cdot \partial_{y_3}\widetilde\varphi'\left(x',{x_3\over \ep}\right)\,dx\\
\noame &\dis+\ep^{1-\ell\over 2}\int_{\Omega_\ep^+}D {\bf u}_\ep(x): D_z \widehat \varphi\left(x',{x\over \ep^\ell}\right)\,dx+O_\ep.
\end{array}
\end{equation}
{\color{black} Observe that main order terms are defined in $\Omega_\varepsilon^+= \omega\times (0,\varepsilon h)$. To obtain this, here we have used that, in $\Omega_\varepsilon^-=\omega\times (-\Psi_\varepsilon(x'),0)$,}
	\[
	\ep^{1-\ell\over 2}\int_{\Omega_\ep^-}D {\bf u}_\ep(x): D_z \widehat \varphi\left(x',{x\over \ep^\ell}\right)\,dx=O_\ep,
	\] 
	Last estimate results from Cauchy-Schwarz inequality, the estimate of $D{\bf u}_\varepsilon$ given in (\ref{vestim}) and the estimate
	$$\int_{\Omega_\varepsilon^-}|\ep^{1-\ell\over 2}D_z\widehat \varphi|^2\,dx\leq C \varepsilon^{1-\ell }|\Omega_\varepsilon^-|=C\varepsilon^{{1-\ell}+\delta}=C\varepsilon^{1+\ell\over 2}.$$	
{\color{black} Also, we have used that
$$	\int_{\Omeps} D\ub_\eps: h_\varepsilon\,dx\leq \| D\ub_\eps\|_{L^2(\Omega_\varepsilon)^{3\times 3}}\|h_\varepsilon\|_{L^2(\Omega_\varepsilon)^{3\times 3}}\leq C\varepsilon^{-{1\over 2}}\|h_\varepsilon\|_{L^2(\Omega_\varepsilon)^{3\times 3}}=O_\varepsilon.$$
	
}
	
 Now, using  the dilatation (\ref{changevar1}) and (\ref{convergencias_solutions})$_1$, we get 
$$ \int_{\Omega_\ep^+}\partial_{x_3} {\bf u}_\ep'(x)\cdot \partial_{y_3}\widetilde\varphi\left(x',{x_3\over \ep}\right)\,dx=\int_{\Omega} \partial_{y_3} \widetilde {\bf u}_\ep'(y)\cdot \partial_{y_3}\widetilde\varphi'(y)\,dy=\int_{\Omega} \partial_{y_3}\widetilde {\ub}'(y)\cdot \partial_{y_3}\widetilde \varphi'(y)\,dy+O_\ep.
$$
Next, from the {\it unfolding} (\ref{changevar2}),  the hypothesis of the support of $D_z\widehat \varphi$ and (\ref{Dhatuw})$_1$, we deduce
\begin{equation}\label{Du_Dhatu}\begin{array}{rl}\dis\ep^{1-\ell\over 2}\int_{\Omega_\ep^+}D {\bf u}_\ep(x): D_z \varphi\left(x',{x\over \ep^\ell}\right)\,dx
&\dis=\int_{\omega\times\widehat Q_M} D_z \left(\ep^{1-\ell\over 2}\widehat {\bf u}_\ep(x',z)\right) : D_z\widehat \varphi(x',z)\,dx'dz +O_\ep\\
\noame &\dis
=\int_{\omega\times\widehat Q}D_z  \widehat {\bf u}(x',z) : D_z  \widehat \varphi(x',z)\,dx'dz+ O_\ep,
\end{array}\end{equation}
{\color{black} where $\widehat Q_M=Z'\times (0,M)$ and $\widehat Q=Z'\times
(0,+\infty)$. For more details of the unfolding change, we refer to \cite{CLS-SIMA}.}
Then, we have that (\ref{F1-v}) is given by
\begin{equation}\label{F1}
\int_{\Omeps} D\ub_\eps  : D\varphi \,dx=
\int_{\Omega} \partial_{y_3}\widetilde {\ub}'(y)\cdot \partial_{y_3}\widetilde \varphi'(y)\,dy+\int_{\omega\times\widehat Q}D_z\widehat {\bf u}(x',z):D_z\widehat\varphi(x',z)\,dx'dz+ O_\ep.
\end{equation}

\noindent $\bullet$ {\it 2nd term of (\ref{vf_vel_1})}. Using the decomposition (\ref{decomposition}),  (\ref{phieptest}) and (\ref{gbounded1}), we have that 
$$-\int_{\Omega_\ep}p_\ep\,{\rm div}\,\varphi_\ep\,dx=\int_{\Omega_\ep}\nabla_{x'}p_\ep^0(x')\cdot\varphi'_\ep(x)\,dx-\int_{\Omega_\ep}p_\ep^1(x)\,{\rm div}\,\varphi_\ep(x)\,dx\,.$$
Applying the change of variables (\ref{changevar1})  and (\ref{conv_p0_p1})$_1$ to the first integral, we get
$$\begin{array}{l}
\dis\int_{\Omega_\ep}\nabla_{x'}p_\ep^0(x')\cdot \varphi_\ep'(x)\,dx=\ep\int_{\Omega_\ep^+}\nabla_{x'}  p_\ep^0(x')\cdot \widetilde\varphi'\left(x',{x_3\over \ep}\right)\,dx+O_\ep\\
\noame\dis
=\int_{\Omega}\ep^2\nabla_{y'}  p_\ep^0(y')\cdot \widetilde\varphi'(y)\,dy+O_\ep
=\int_{\Omega}\nabla_{y'} p(y')\cdot \widetilde \varphi'(y)\,dy + O_\ep.
\end{array}
$$
{\color{black} Here, we have used that 
$$\int_{\Omega_\varepsilon}\nabla_{x'}p_\varepsilon^0(x')g_\varepsilon(x) dx\leq \|\nabla_{x'}p_\varepsilon^0(x')\|_{L^2(\Omega_\varepsilon)}\|g_\varepsilon\|_{L^2(\Omega_\varepsilon)}\leq C \varepsilon^{-{3\over 2}}\|g_\varepsilon\|_{L^2(\Omega_\varepsilon)}=O_\varepsilon.$$
}

For the second integral, using  the change of variables (\ref{changevar2})  and (\ref{convp1}) {\color{black} (see \cite{CLS-SIMA} for more details)}, we obtain 
$$
\begin{array}{l}\dis
\int_{\Omega_\ep}p^1_\ep(x){\rm div}\,\varphi_\ep(x)\,dx=\ep^{1-\ell\over 2}\int_{\Omega_\ep^+}p_\ep^1(x)\,{\rm div}_z\,\widehat \varphi\left(x',{x\over \ep^\ell}\right)\,dx+O_\ep\\
\noame\dis
=\int_{\omega\times \widehat Q_M}
 \ep^{1+\ell \over 2}\widehat  p_\ep^1(x',z) \,{\rm div}_z\widehat \varphi(x',z)\,dx'dz+O_\ep=\int_{\omega\times \widehat Q}\widehat p^1(x',z)\,{\rm div}_z\widehat \varphi(x',z)\,dx'dz+O_\ep.
\end{array}$$
Then, we get 
\begin{equation}\label{p0}-\int_{\Omega_\ep}p_\ep\,{\rm div}\,\varphi_\ep\,dx=\int_{\Omega}\nabla_{y'}  p(y')\cdot \widetilde \varphi'(y)\,dy-\int_{\omega\times \widehat Q}\widehat p^1(x',z)\,{\rm div}_z\widehat \varphi(x',z)\,dx'dz+O_\ep\,.
\end{equation}

\noindent $\bullet$ {\it 3rd term of (\ref{vf_vel_1})}. Using  (\ref{Dphieptest}), (\ref{Dgbounded}), the change of variables (\ref{changevar1}) and (\ref{conv_w_ep}), we have
\begin{equation}\label{FourthT_1}\begin{array}{ll}\dis
-2N^2\int_{\Omega_\ep}{\bf w}_\ep\cdot (\nabla\times \varphi_\ep) \,dx=&\dis -2N^2 \int_{\Omega_\ep^+}   {\bf w}_\ep(x)\cdot\left({\rm rot}_{y_3} \widetilde \varphi\left(x',{x_3\over \ep}\right)\right)\,dx + O_\ep\\
\noame &\dis 
=-2N^2\int_{\Omega}\ep \widetilde {\bf w}_\ep'(y)\cdot {\rm rot}_{y_3}\widetilde \varphi'(y) \,dy + O_\ep\\
\noame &\dis
=-2N^2\int_{\Omega}\widetilde {\bf w}'(y)\cdot {\rm rot}_{y_3}\widetilde \varphi' (y)\,dy + O_\ep\,.
\end{array}\end{equation}
{\color{black}
Among others, here we have used that
$$\int_{\Omega_\ep}{\bf w}_\ep\cdot (\nabla\times g_\varepsilon)(x)\,dx\leq \|{\bf w}_\ep\|_{L^2(\Omega_\varepsilon)^3}\|h_\varepsilon\|_{L^2(\Omega_\varepsilon)^{3\times 3}}\leq C\varepsilon^{-{1\over 2}}\|h_\varepsilon\|_{L^2(\Omega_\varepsilon)^{3\times 3}}=O_\varepsilon.$$

}
Integrating by parts by using  the formula (\ref{bypartsform}),  we then get
\begin{equation}\label{FourthT}
-2N^2\int_{\Omega_\ep} {\bf w}_\ep\cdot(\nabla\times \varphi_\ep)\,dx
=-2N^2\int_{\Omega}{\rm rot}_{y_3}\widetilde {\bf w} '(y)\cdot \widetilde \varphi'(y)\,dy -2N^2\int_{\Gamma}[\widetilde {\bf w}']^\perp(y)\cdot\widetilde\varphi'(y)\,d\sigma+ O_\ep.
\end{equation}

\noindent $\bullet$ {\it  4th term of (\ref{vf_vel_1})}. 
From ${\bf w}_\ep=0$ on $\Gamma_\ep^1$ and  estimate (\ref{westim}), we obtain
$$\int_{\Gamma_\ep^0}|{\bf w}_\ep|^2d\sigma\leq C \ep\int_{\Omega_\ep}|D{\bf w}_\ep|^2\,dx\leq C\ep^{-2}.$$
Then, from (\ref{phieptest}), (\ref{gbounded2}) and (\ref{testfunc})$_2$, we deduce
$$\begin{array}{l}\dis
-2\left({1\over\alpha}-N^2\right)\int_{\Gamma_\ep^0}({\bf w}_\ep\times {\bf n}_\varepsilon)\cdot\varphi_\ep\,d\sigma\\
\noame\dis= \dis -2\left({1\over\alpha}-N^2\right)\ep\int_{ \omega}\left({\bf w}_\ep\left(x',-\lambda\ep^{3\ell-1\over 2}\Psi\left({x_1\over \ep^\ell}\right)\right)\times {\bf n}_\ep\right)\cdot\widetilde \varphi(x',0)\,\\
\noame\dis\qquad\qquad\qquad\qquad\qquad \sqrt{1+ \lambda^2\left(\ep^{\ell-1\over 2}\right)^2\left|\partial_{x_1}\Psi\left({1\over \ep^\ell}x'\cdot e_1'\right)\right|^2}\,dx'+ O_\ep
\\
\noame\dis= \dis -2\left({1\over\alpha}-N^2\right)\ep\int_{ \omega}\left({\bf w}_\ep\left(x',-\lambda\ep^{3\ell-1\over 2}\Psi\left({x'\over \ep^\ell}\right)\right)\times {\bf n}\right)\cdot\widetilde \varphi(x',0)\,dx'+ O_\ep\,,
\end{array}$$
where ${\bf n}=(0,0,-1)$.  {\color{black} Here we have used that 
$$\int_{\Gamma_\ep^0}({\bf w}_\ep\times {\bf n}_\varepsilon)\cdot\varphi_\ep\,d\sigma\leq \|{\bf w}_\ep\|_{L^2(\Gamma_\ep^0)^3}\|g_\varepsilon\|_{L^2(\Gamma_\ep^0)^3}\leq C\varepsilon^{-1}\|g_\varepsilon\|_{L^2(\Gamma_\ep^0)^3}=O_\varepsilon.$$
}
By means of integration in the variable $x_3$, we get 
\begin{equation}\label{comparacion}
\int_\omega \left|\ep {\bf w}_\ep\left(x',-\lambda\ep^{3\ell-1\over 2}\Psi\left({x_1\over \ep^\ell}\right)\right)-\ep {\bf w}_\ep(x',0)\right|^2\,dx'\leq C\ep^{3\ell+3\over 2}\int_{\Omega_\ep}|D{\bf w}_\ep|^2\,dx\leq C\ep^{3\ell-1\over 2}.
\end{equation}
Then, from ${\bf w}_\ep(x',0)=\widetilde {\bf w}_\ep(x',0)$ and (\ref{conv_w_ep}), we obtain
\begin{equation}\label{intboundary}\begin{array}{rl}
\dis-2\left({1\over\alpha}-N^2\right)\int_{\Gamma_\ep^0}({\bf w}_\ep\times {\bf n}_\ep)\cdot\varphi_\ep\,d\sigma
&\dis =-2\left({1\over\alpha}-N^2\right)\int_{ \omega}(\ep \widetilde {\bf w}_\ep(y',0)\times {\bf n})\cdot\widetilde \varphi(y',0)\,dy'+O_\ep
\\
\noame& \dis
=-2\left({1\over\alpha}-N^2\right)\int_{ \omega}[\widetilde {\bf w}']^\perp(y',0)\cdot\widetilde \varphi'(y',0)\,dy'+O_\ep.
\end{array}\end{equation}
By considering (\ref{F1}), (\ref{p0}), (\ref{FourthT}) and (\ref{intboundary}), then we get that  $\widetilde {\bf u}'$, $\widetilde {\bf w}'$, $p$, $\widehat {\bf u}$ and $\widehat p^1$ satisfy
\begin{equation}\label{limitvf-1}
\begin{array}{l}
\dis \int_{\Omega} \partial_{y_3}\widetilde {\ub}'(y)\cdot \partial_{y_3}\widetilde \varphi'(y)\,dy+\int_{\omega\times\widehat Q}D_z\widehat {\bf u}(x',z):D_z\widehat\varphi(x',z)\,dx'dz +\int_{\Omega}\nabla_{y'} p(y')\cdot \widetilde \varphi'(y)\,dy \\
\noame\dis 
-\int_{\omega\times \widehat Q}\widehat p^1(x',z)\,{\rm div}_z\widehat \varphi(x',z)\,dx'dz -2N^2\int_{\Omega}{\rm rot}_{y_3}\widetilde {\bf w} '(y)\cdot \widetilde \varphi'(y)\,dy\\
\noame\dis -{2\over\alpha}\int_{ \omega}[\widetilde {\bf w}']^\perp(y',0)\cdot\widetilde \varphi'(y',0)\,dy'=0
\end{array}
\end{equation}
for any $\widetilde\varphi'\in C_c^1(\omega\times(-h,h))^2$, $\widehat \varphi\in C_c^1(\omega;C_{\#}^1(\widehat Q))^3$ satisfying (\ref{testfunc}), and then, by  arguments of density, for any $\widetilde\varphi'\in H^1(0,h;L^2(\omega))^2$ and any $\widehat \varphi\in L^2(\omega;\mathcal{V})^3$ satisfying
$$\widetilde \varphi(x',h)=0\hbox{ a.e. }x'\in\omega,\quad \lambda\partial_{z_1}\Psi(z'\cdot e_1')\widetilde\varphi_1(x',0)+\widehat \varphi_3(x',z',0)=0\hbox{ a.e. }(x',z')\in \omega\times Z'.$$
By considering $\widetilde \varphi'=0$ in (\ref{limitvf-1}), we get that  $(\widehat {\bf u}, \widehat p^1)\in \mathcal{V}^3\times L^2_{\#}(\widehat Q)$ solves

\begin{equation}\label{phiproblem}
\left\{\begin{array}{rl}\dis
-\Delta_z \widehat {\bf u}+\nabla_z \widehat p^1=0&\hbox{ in }\mathbb{R}^2\times \mathbb{R}^+,\\
\noame\dis
{\rm div}_z\widehat {\bf u}=0&\hbox{ in }\mathbb{R}^2\times \mathbb{R}^+,
\\
\noame\dis
\widehat u_3(x',z',0)=-\lambda\partial_{z_1}\Psi(z'\cdot e_1')\widetilde u_1(x',0)&\hbox{ on }\mathbb{R}^2\times\{0\},\\
\noame\dis
\partial_{z_3}\widehat {\bf u}'=0&\hbox{ on }\mathbb{R}^2\times\{0\},
\end{array}\right.
\end{equation}
a.e. $x'$ in $\omega$. Setting $(\widehat \phi^{1,\lambda}, \widehat q^{1,\lambda} )$ by (\ref{system_phi_1}), we derive 
\begin{equation}\label{Dzphip1}
\begin{array}{l}\dis
D_z\widehat {\bf u}(x',z)=-\widetilde u_1(x',0) D_z\widehat \phi^{1,\lambda}(z), \quad \hbox{ a.e. in }\mathbb{R}^2\times\mathbb{R}^+,\\
\noame\dis
\widehat p^1(x',z)=-\widetilde u_1(x',0) \widehat q^{1,\lambda}(z), \quad \hbox{ a.e. in }\mathbb{R}^2\times\mathbb{R}^+.
\end{array}
\end{equation}
Next, by considering $\widetilde\varphi'\in H^1(0,h;L^2(\omega))^2$, $\widetilde\varphi'(x',h)=0$, a.e. $x'\in \omega$, as test function in  (\ref{limitvf-1}), setting $\widehat \varphi$ by
$$\widehat \varphi(x',z)=-\widetilde \varphi_1(x',0)\widehat \phi^{1,\lambda}(z),$$
and taking into account (\ref{Dzphip1}), we obtain
\begin{equation}\label{limitvf-2}
\begin{array}{l}
\dis \int_{\Omega} \partial_{y_3}\widetilde {\ub}'(y)\cdot \partial_{y_3}\widetilde \varphi'(y)\,dy+\int_{\Omega}\nabla_{y'} p(y')\cdot \widetilde \varphi'(y)\,dy
-2N^2\int_{\Omega}{\rm rot}_{y_3}\widetilde {\bf w} '(y)\cdot \widetilde \varphi'(y)\,dy\\
\noame\dis 
+ E_\lambda\int_{\omega} \widetilde u_1(y',0)\,\widetilde \varphi_1(y',0)\,dy' -{2\over\alpha}\int_{ \omega}[\widetilde {\bf w}']^\perp(y',0)\cdot\widetilde \varphi'(y',0)\,dy'=0,\end{array}
\end{equation}
where $E_\lambda\in \mathbb{R}$ is given by (\ref{matrixM})$_1$. \\

{\it Step 3. Passing to the limit in the variational formulation (\ref{vf_mic_1}). } This step is similar to the previos step, so we will only give some details.\\

\noindent $\bullet$ {\it  1st term of (\ref{vf_mic_1})}.  Analogously to the first step of the previous variational formulation, we deduce
\begin{equation}\label{FTw} 
\begin{array}{rl}
\dis \ep^2R_c\int_{\Omega_\ep}D{\bf w}_\ep:D\psi_\ep\,dx =& \dis \ep R_c\int_{\Omega_\ep^+}\partial_{x_3} {\bf w}'_\ep(x)\cdot \partial_{y_3}\widetilde\psi'\left(x',{x_3\over \ep}\right)\,dx\\
&\dis 
+\ep^{3-\ell\over 2}R_c\int_{\Omega_\ep^+} D {\bf w}_\ep(x):D_z\widehat\psi\left(x',{x\over \ep^\ell}\right)\,dx+O_\ep,
\end{array}
\end{equation}
and by using the changes of variables, we can pass to the limit in every terms by obtaining 
\begin{equation}\label{FTw_2} \ep^2 R_c\int_{\Omega_\ep}D{\bf w}_\ep:D\psi_\ep\,dx=R_c\int_{\Omega}\partial_{y_3}\widetilde {\bf w}'(y)\cdot\partial_{y_3}\widetilde\varphi'(y)\,dy + R_c\int_{\omega\times \widehat Q}D_z\widehat {\bf w}(x',z):D_z\widehat\psi(x',z)\,dx'dz+ O_\ep.
\end{equation}

\noindent $\bullet$ {\it  2nd term of (\ref{vf_mic_1})}. 
From ${\bf u}_\ep=0$ on $\Gamma_\ep^1$ and the estimate (\ref{vestim}), we get
$$\int_{\Gamma_\ep^0}|{\bf u}_\ep|^2d\sigma\leq C \ep\int_{\Omega_\ep}|D{\bf u}_\ep|^2\,dx\leq C.$$
Using this, and proceeding analogously to the development of the fourth term in (\ref{vf_vel_1}), we get 

\begin{equation}\label{fifthterm_w}\begin{array}{l}\dis
\dis -2N^2\beta\int_{\Gamma_\ep^0} ({\bf u}_\ep\times {\bf n}_\ep)\cdot\psi_\ep\,d\sigma=
-2N^2\beta \int_{\omega}[\widetilde {\bf u}' (y',0)]^\perp\cdot\widetilde \psi'(y',0)\,dy'+ O_\ep.
\end{array}
\end{equation}

\noindent $\bullet$ {\it  3rd term of (\ref{vf_mic_1})}. Applying (\ref{g_2bounded1}), the change of variables (\ref{changevar1}) and convergence (\ref{conv_w_ep}), we have
\begin{equation}\label{ST}\begin{array}{l}
\dis 4N^2\int_{\Omega_\ep}{\bf w}_\ep\cdot\psi_\ep\,dx= 4N^2\int_{\Omega} {\bf w}_\ep'(x)\cdot \widetilde\psi'\left(x',{x_3\over \ep}\right)dx+ O_\ep\\
\noame \dis  =4N^2\int_{\Omega} \ep \widetilde {\bf w}_\ep'(y)\cdot \widetilde\psi'(y)\,dy+ O_\ep= 4N^2\int_{\Omega}\widetilde {\bf w}'(y)\cdot \widetilde\psi'(y)\,dy+O_\ep.
\end{array}\end{equation}

\noindent $\bullet$ {\it  4th term of (\ref{vf_mic_1})}.
Similarly to (\ref{FourthT_1}) and taking into account convergence  (\ref{convergence_u}), we have 
\begin{equation}\label{ST_2}
-2N^2\int_{\Omega_\ep}(\nabla\times{\bf u}_\ep)\cdot \psi_\ep\,dx=-2N^2\int_{\Omega} {\rm rot}_{y_3}\widetilde {\bf u}'(y)\cdot \widetilde\psi'(y)\,dy+ O_\ep.
\end{equation}
Finally, from (\ref{FTw_2}) to (\ref{ST_2}), we have obtained
\begin{equation}\label{varform_limit_w}\begin{array}{l}
\dis
R_c\int_{\Omega}\partial_{y_3}\widetilde {\bf w}'(y)\cdot\partial_{y_3}\widetilde\psi'(y)\,dy + R_c\int_{\omega\times \widehat Q}D_z\widehat {\bf w}(x',z):D_z\widehat\psi(x',z)\,dx'dz+4N^2\int_{\Omega}\widetilde {\bf w}'(y)\cdot \widetilde\psi'(y)\,dy\\
\noame
\dis
-2N^2\int_{\Omega} {\rm rot}_{y_3}\widetilde {\bf u}'(y)\cdot \widetilde\psi'(y)\,dy-2N^2\beta \int_{\omega}[\widetilde {\bf u}' (y',0)]^\perp\cdot\widetilde \psi'(y',0)\,dy'=0,
\end{array}\end{equation}
for any $\widetilde\psi'\in C_c^1(\omega\times(-1,1))^2$, $\widehat \psi\in C_c^1(\omega;C_{\#}^1(\widehat Q))^3$ satisfying  (\ref{testfunc}), and by argument of density, for any $\widetilde\psi'\in H^1(0,1;L^2(\omega))^2$, $\widehat \psi\in L^2(\omega;\mathcal{V})^3$ satisfying
$$\widetilde \psi(x',1)=0,\hbox{ a.e. }x'\in\omega,\quad \lambda\partial_1\Psi(z_1)\widetilde\psi_1(x',0)+\widehat \psi_3(x',z',0)=0,\hbox{ a.e. }(x',z')\in \omega\times Z'.$$
Similarly tot he previous step, we eliminate  $\widehat {\bf w}$ from (\ref{varform_limit_w}). To do this, we consider $\widetilde \psi'=0$ in (\ref{varform_limit_w}), which gives that $\widehat {\bf w}\in \mathcal{V}^3\times L^2_{\#}(\widehat Q)^3$ solves
\begin{equation}\label{phiproblem-2}
\left\{\begin{array}{rl}\dis
-\Delta_z \widehat {\bf w}=0&\hbox{ in }\mathbb{R}^2\times \mathbb{R}^+,
\\
\noame\dis
\widehat w_3(x',z',0)=-\lambda\partial_1\Psi(z_1)\widetilde w_1(x',0)&\hbox{ on }\mathbb{R}^2\times\{0\},\\
\noame\dis
\partial_{z_3}\widehat {\bf w}'=0&\hbox{ on }\mathbb{R}^2\times\{0\},
\end{array}\right.
\end{equation}
a.e. $x'$ in $\omega$. Setting $\widehat \phi^{2,\lambda}$ by (\ref{system_phi_2}), we derive
\begin{equation}\label{Dzphip2}
\begin{array}{l}\dis
D_z\widehat {\bf w}(x',z)=-\widetilde w_1(x',0) D_z\widehat \phi^{2,\lambda}(z), \quad \hbox{ a.e. in }\mathbb{R}^2\times\mathbb{R}^+,\end{array}
\end{equation}
Next, by considering $\widetilde\psi'\in H^1(0,h;L^2(\omega))^2$, $\widetilde\psi'(x',h)=0$, a.e. $x'\in \omega$ as test function in  (\ref{varform_limit_w}), defining  $\widehat \psi$ by
$$\widehat \psi(x',z)=-\widetilde \psi_1(x',0)\widehat\phi^{2,\lambda}(z),$$
and from (\ref{Dzphip2}), we deduce
\begin{equation}\label{varform_limit_w-2}\begin{array}{l}
\dis
R_c\int_{\Omega}\partial_{y_3}\widetilde {\bf w}'(y)\cdot\partial_{y_3}\widetilde\psi'(y)\,dy+4N^2\int_{\Omega}\widetilde {\bf w}'(y)\cdot \widetilde\psi'(y)\,dy
-2N^2\int_{\Omega} {\rm rot}_{y_3}\widetilde {\bf u}'(y)\cdot \widetilde\psi'(y)\,dy\\
\noame\dis
+ R_c  F_\lambda\int_{\omega} \widetilde w_1(y',0)\widetilde \psi_1(y',0)\,dy'
-2N^2\beta \int_{\omega}[\widetilde {\bf u}' (y',0)]^\perp\cdot\widetilde \psi'(y',0)\,dy'=0\,,\end{array}\end{equation}
where $F_\lambda$ is given by (\ref{matrixM})$_2$.\\

{\it Step 4. Conclusion. }  Since  $\varphi'$ and $\widetilde \psi'$ are arbitrary, we derive from (\ref{limitvf-2}) and (\ref{varform_limit_w-2}) that $(\widetilde \ub',\widetilde \wb', p)$  satisfies  the system (\ref{limit_system_1})$_{1,2}$ with boundary conditions (\ref{bc_botom_case2_system_1}). To ensure that the whole sequences $\widetilde {\bf u}_\ep, \widetilde {\bf w}_\ep, \ep \widetilde p_\ep$ converge, it remains to prove the existence and uniqueness of weak solution of the effective system  (\ref{limit_system_1})--(\ref{bc_botom_case2_system_1}). 

Note that we can always reduce the non vanishing divergence problem  (\ref{limit_system_1})--(\ref{bc_botom_case2_system_1}) to a free divergence problem by considering the lift function 
${\bf J}\in H^1(\Omega)^3$ such that
$$\left\{\begin{array}{l}
{\rm div}_{y} {\bf J}=0\hbox{ in }\Omega,\\
{\bf J}\cdot n=0 \hbox{ on }   \Gamma^{\ell},\\
{\bf J}=-Se_3\hbox{ on }\Gamma^1,\\
{\bf J}=0\hbox{ on }\Gamma,
\end{array}\right.$$
and hence using the change of unknowns ${\bf U}' ={\bf u}' -{\bf J}' $. Therefore, it is sufficient to study the existence and uniqueness of the weak solution of the problem
$$\left\{\begin{array}{rl}
\dis
-\partial_{y_3}^2 \widetilde \ub' + \nabla_{y'}   p -2N^2 {\rm rot}_{y_3} \widetilde \wb'=\partial_{y_3}^2{\bf J}'& \hbox{ in } \Omega,\\
\noame
-R_c\partial_{y_3}^2 \widetilde \wb'+4N^2 \widetilde \wb'-2N^2{\rm rot}_{y_3}\widetilde \ub'=2N^2{\rm rot}_{y_3}{\bf J}'& \hbox{ in } \Omega,\\
\noame
\dis \mathrm{div}_{y'}\int_0^{h} \widetilde {\bf u}'(y',y_3)\,dy_3=0& \hbox{ in }  \omega.
\end{array}\right.$$
with  the boundary conditions 
$$
\hspace{-0.8cm}\widetilde \ub'=0,\ \widetilde \wb'=0 \hbox{ on }\omega\times\{h\},
$$
$$
\partial_{y_3} \widetilde \ub'=-{2\over \alpha} [ \widetilde \wb']^\perp + E_\lambda \widetilde (\ub'\cdot e'_1)\, e_1'\ \hbox{ on }\Gamma,\quad 
 R_c\partial_{y_3} \widetilde \wb'=-2N^2\beta[\widetilde \ub']^\perp+R_c F_\lambda\widetilde (\wb'\cdot e_1')\, e_1'\ \hbox{ on }\Gamma\,.
$$
The existence and uniqueness of this problem follows the lines of the proof of Theorem \ref{thm_existence} with a flat bottom, taking into account suitable spaces and the new boundary terms $E_\lambda \widetilde \ub'\cdot e'_1$ and $R_c F_\lambda\widetilde \wb'\cdot e_1'$ and the new source terms $\partial_{y_3}^2{\bf J}'$ and $2N^2{\rm rot}_{y_3}{\bf J}'$ in the variational formulation.

\par \hfill $\square$
\\

The proofs of Lemmas \ref{lemma_alpha_neq_1}, \ref{lemma_alpha_igual_1} and \ref{lm_super_critical} and Corollary \ref{cor_alpha_neq_1} are given in the Appendix. We finish this section by providing the proofs of Theorems \ref{thm_reynolds} and \ref{thm_reynolds_super}.\\

\noindent {\bf Proof of Theorem \ref{thm_reynolds}. } Using (\ref{limit_system_1_reynolds2}), we obtain 
$$-\int_0^1\left(\int_0^{h}\widetilde u_1(y_1,y_3)\,dy_3\right)\partial_{y_1}\theta(y_1)\,dy_1=\int_0^1S\theta(y_1)\,dy_1, \quad \forall\,\theta\in H^1((0,1)).$$
From Lemmas \ref{lemma_alpha_neq_1} and \ref{lemma_alpha_igual_1}, by averaging (\ref{tilde_u_1}) and (\ref{tilde_u_1_alpha_igual_1}) we obtain 
$$\int_0^{h}\widetilde u_1(y_1,y_3)\,dy_3=-\Theta_{\lambda} \partial_{y_1} p(y_1).$$
\par\hfill$\square$\\

\noindent {\bf Proof of Theorem \ref{thm_reynolds_super}. } Using (\ref{limit_system_1_reynolds2}), we obtain 
$$-\int_0^1\left(\int_0^{h}\widetilde u_1(y_1,y_3)\,dy_3\right)\partial_{y_1}\theta(y_1)\,dy_1=\int_0^1S\theta(y_1)\,dy_1, \quad \forall\,\theta\in H^1((0,1)).$$
From Lemma \ref{lm_super_critical}, by averaging (\ref{lm_super_critical_u1})$_1$, we obtain 
$$\int_0^{h}\widetilde u_1(y_1,y_3)\,dy_3=-\Theta \partial_{y_1} p(y_1).$$
\par\hfill$\square$\\

\section{Application to squeeze-film bearing}\label{Sec:num}

As a an application of the results presented in Section~\ref{Sec:MainResults}, we consider in this section a squeeze-film bearing composed of two parallel plates separated by a micropolar fluid film. The lower surface is at rest and composed of a rough material, while the smooth upper surface is under normal squeeze motion. Hence, the distance between the plates is a decreasing function of time, which is expected to go to zero as the fluid is squeezed out of the gap. We assume that the motion of the upper plate is slow, so that the inertial effects can be neglected and the behaviour of the bearing can be captured by a quasi-static model.

\subsection{Derivation of the model}\label{Subsec:DerivationModel}

Let $T_0>0$ be the characteristic time of the motion. Following the notation from Section~\ref{Sec:PositionOfPb}, we denote by $L$ the horizontal dimension of the bearing and by $hc(T_0t)$ the distance between the plates at time $t$, where $t$ stands for the dimensionless time variable and $h>0$ is an adimensional constant. 
We introduce the small parameter $\eps(t)=c(T_0t)/L$ and assume that, at each time $t$, the fluid flow is described by the system~\eqref{1,1,1}--\eqref{bc5} with $\eps=\eps(t)$, where all the constant $\nu,\nu_r,c_a,c_d,\alpha,\beta$ are fixed. However, the velocity $\overline{V}$ of the upper plate is related to the load $\overline W$ applied on the bearing, which is assumed to be independant on $t$, through the implicit relation
\begin{equation}\label{Num:prescribed_load}
	\overline W = \int_{\overline\Gamma_{\ep(t)}} \overline p_{\ep(t)}\, .
\end{equation}

{
\color{black}
We assume that the rough bottom is composed of periodically distributed riblets, described by a given function $\overline \Psi$ of the form
\begin{equation}\label{Num:defPsi}
\overline{\Psi}(\overline{x}_1) = L\Lambda \Psi\left( \frac{1}{LM} \overline x_1 \right) 
\end{equation}
where $\Lambda,M>0$ correspond respectively to the amplitude and period of the ribbed surface, divided by $L$, and where the function $\Psi$ is $1$-periodic and regular, and satisfies~\eqref{Psi-normalized}.

Using the notation introduced in Section 2, for a given time $t$, the fluid domain $\overline \Omega_{\ep(t)}$ is given by
\[
\overline \Omega_{\ep(t)}
=\left\{(\overline x',\overline x_3)\in (L\omega)\times \mathbb{R}\, ,\   - \overline \Psi_{\ep(t)}(\overline x_1)<\overline x_3<h\, L\, \ep(t)\right\}\, ,
\]
where according to definition~\eqref{rugprofile}, the function $\overline\Psi_{\ep(t)}$ takes the form
\begin{equation}\label{Num:psi_eps_bar}
\overline \Psi_{\ep(t)}(\overline x_1) = L  \lambda(t) \ep(t)^{\delta(t)}\Psi\left( {1 \over L \ep(t)^{\ell(t)}}  \overline x_1\right)\, .
\end{equation}
Since the geometry of the lower plate is, in fact, independent on time, $\overline\Psi_{\ep(t)}$ coincides with the fixed profile $\overline \Psi$ given by~\eqref{Num:defPsi}. Hence, parameters $\lambda(t),\delta(t)$ satisfy
\begin{equation}\label{Num:def_Lambda-M}
\lambda(t)\, \eps(t)^{\delta(t)} = \Lambda\quad \text{and}\quad \ep(t)^{\ell(t)}=M\, .
\end{equation}
Considering the critical regime $\delta(t)=\frac{3}{2}\ell(t)-\frac{1}{2}$, the parameters $\ell(t),\delta(t),\lambda(t)$ can thus be expressed as functions of $\ep(t),M, L$, and in particular the parameter $\lambda(t)$ is given by
$
\lambda(t) = \Lambda\, \ep(t)^{\frac{1}{2}-\frac{3}{2}\frac{\ln M}{\ln\ep(t)}}
$, which simplifies to
\begin{equation}\label{Num:def_lambda}
\lambda(t) = \frac{\Lambda}{M^{3/2}}\, \ep(t)^{1/2}\, .
\end{equation}
}


Analogously to~\eqref{non-dim}, we introduce the dimensionless pressure $p_{\ep(t)}$ defined by $p_{\ep(t)}= \frac{L}{V_0(\nu+\nu_r)}\overline p_{\ep(t)}$, and accordingly, the adimensional load $W=\frac{\overline W}{L V_0(\nu+\nu_r)}$. Defining also the rescaled pressure $\widetilde p_{\ep(t)}$, the constraint~\eqref{Num:prescribed_load} can be rephrased as
\begin{equation}\label{Num:dimensionless_load}
	\int_\omega \widetilde p_{\ep(t)}(y',h)\, dy' = W \,  .	
\end{equation}
Since the system is independent on the $x_2$-direction, we take $\omega=(0,1)$ and apply Theorems~\ref{thm_effective} and~\ref{thm_reynolds} to approximate $\widetilde p_{\ep(t)}$ by $p^{\lambda(t),S(t)}/\ep(t)^2$, where $p^{\lambda(t),S(t)}$ is the solution of Reynolds equation~\eqref{reynolds}, with $\lambda=\lambda(t)$ and $S=S(t)$, $S(t)$ being fixed consistently with relation~\eqref{Num:dimensionless_load}.

{\color{black} Indeed, for a given time $t$, if one assumes that $\eps(t)$ has been computed, all roughness parameters $\lambda(t),\ell(t),\delta(t)$ are then prescribed by relations~\eqref{Num:def_Lambda-M} and by the critical relation $\delta(t)=\frac{3}{2}\ell(t)-\frac{1}{2}$. Setting $\lambda=\lambda(t),\ell=\ell(t),\delta=\delta(t)$, the state of the system at time $t$ is thus described by~\eqref{1,1}--\eqref{bc-5}, where $\eps$ takes the value $\eps(t)$. All other parameters being fixed (equal to their values at instant $t$), we can define a sequence of problems~\eqref{1,1}--\eqref{bc-5}, where $\eps\in (0,\eps(t))$ is the only parameter that is let to zero. Since $\delta$ and $\ell$ are bound by relation $\delta=\frac{3}{2}\ell-\frac{1}{2}$, and $\eps(t)$ is small, we can apply the asymptotic analysis result from Theorems~\ref{thm_effective} and~\ref{thm_reynolds} to replace the solution of problem~\eqref{1,1}--\eqref{bc-5}, with $\eps=\eps(t)$, by the solution of the effective problem~\eqref{limit_system_1}--\eqref{bc_botom_case2_system_1}.}

{\color{black}
Let us point that this strategy aims at approximating the configuration of the system at each instant $t$, from a computational point of view, by the limit provided by Theorems~\ref{thm_effective} and~\ref{thm_reynolds}. The study of the physical limit of the system when $\eps(t)$ goes to zero (which implies that $\lambda(t)$ also goes to zero by relation~\eqref{Num:def_lambda}) is beyond the scope of the paper. In particular, the value of the physical parameter $\eps(t)$, albeit small, remains greater than the positive quantity $\eps_0/2$ throughout the simulations.
}

\medskip

As is usual in the lubrication field, we impose Dirichlet boundary conditions for the pressure on $x_1\in \{0,1\}$, instead of the Neumann boundary conditions implicitly contained in the weak formulation~\eqref{reynolds}. We obtain the relation
\begin{equation}
	\int_0^1 p^{\lambda(t),S(t)}(y_1)\, dy_1 = \ep(t)^2\, W\, .
\end{equation}
Since Reynolds equation~\eqref{reynolds} is linear with respect to $S$, there holds  $p^{\lambda(t),S(t)}=S(t)\, p^{\lambda(t),1}$ (where $p^{\lambda(t),1}$ satisfies~\eqref{reynolds} with $\lambda=\lambda(t)$ and $S=1$) so that $S(t)$ is given by
\begin{equation}\label{Num:S(t)}
	S(t) = \frac{\ep(t)^2\, W}{ \int_0^1 p^{\lambda(t),1}(y_1)\, dy_1 }\, .
\end{equation}
Using relations~\eqref{Regime:Veps} and~\eqref{non-dim}$_5$, the normal velocity of the upper plate at time $t$, which is given by $-\frac{Lh}{T_0}\, \ep'(t)$, can thus be expressed by $-\frac{Lh}{T_0} \ep'(t) = V_0\, \ep(t)\, S(t)$, yielding the differential equation
\begin{equation}\label{Num:ODE-eps-1}
	\ep'(t) = -\frac{T_0V_0W}{Lh}\frac{ \ep(t)^3}{\int_0^1 p^{\lambda(t),1}(y_1)\, dy_1}\, ,
\end{equation}
where $\lambda(t)$ depends on $\ep(t)$ through relation~\eqref{Num:def_lambda}.

Since $p^{\lambda(t),1}$ satisfies
\[
-\Theta_{\lambda(t)}\partial^2_{y_1} p^{\lambda(t),1} = 1\quad \textrm{in  }(0,1)\, ,\quad  p^{\lambda(t),1}(0)=p^{\lambda(t),1}(1)=0\, ,
\]
is can be expressed by the explicit formula
$
p^{\lambda(t),1}(y_1) = \frac{y_1(1-y_1)}{2\Theta_{\lambda(t)}}
$
so the ODE~\eqref{Num:ODE-eps-1} can be rewritten
\begin{equation}\label{Num:ODE-eps}
	\ep'(t) = -\frac{12T_0V_0W}{Lh}\Theta_{\lambda(t)}\, \ep(t)^3\, .
\end{equation}

\paragraph{Equations of motion}

Denoting respectively by $\kappa$ and $\chi$ the dimensionless constants 
\[
\kappa = \frac{\Lambda^2}{M^3}\, , \quad \chi = \frac{12T_0V_0W}{Lh}
\]
the model can be summarized by the system of equations
\begin{equation}\label{Num:SystemOfEq}
	 \left\lbrace 
	\begin{array}{l}
		\eps'(t)  = -\chi\, \Theta_{\lambda(t)}\, \eps(t)^3\\
		\lambda(t)  = \sqrt{\kappa \eps(t)}
	\end{array}
	\right. 
\end{equation}
which needs to be completed with the definition of an initial state $\eps_0>0$ such that $\eps(0)=\eps_0$.

To facilitate the comparison with the numerical results presented in \cite{Bessonov2, Bayada4, BonSGPaz}, we introduce the dimensionless parameter 
\[
\delta=\frac{R_c}{2N^2\beta}
\]
and define the relative viscosity coefficient $\bar{\nu}_b$ by $\bar{\nu}_b = \nu_b/\nu$, where $\nu_b$ is the boundary viscosity and $\nu$ is the classical viscosity of the fluid. By definition of $\alpha$ (see~\eqref{Def:alpha}), this coefficient can be expressed as 
\[
\bar{\nu}_b = \frac{1-\alpha N^2}{1-N^2}.
\]
Hence, the solution $t\mapsto \eps(t)$ to the system~\eqref{Num:SystemOfEq} depends on the following set of parameters: 
\begin{itemize}
	\item parameters $N, R_c$ characterizing the physical properties of the micropolar fluid,
	\item $\bar\nu_b, \delta$ characterizing the interaction between the micropolar fluid and the upper wall, 
	\item $E, \kappa$ related to the geometry of the rough pattern,
	\item $h, \chi$ associated with characteristic dimensions of the model,
	\item and the initial datum $\eps_0$.
\end{itemize}
Since we are mostly interested in the combined effects of micropolarity and roughness, we will simply put $h=1$ and $\chi=1$ in the sequel.

In the aim of estimating the influence of the parameters on the performance of the squeeze-film bearing, for each set of parameters, we solve numerically the associated system~\eqref{Num:SystemOfEq} (with $\chi=1)$ by a second-order Runge-Kutta method, and compute the "half-life time" $\Thalf$, \emph{i.e.}\! the first instant $t$ such that $\eps(t)<\eps_0/2$, which means that the width of the bearing has been divided by two. In our analysis, the configurations giving rise to the highest values of $\Thalf$ will be considered the most efficient from a mechanical perspective.

\subsection{Numerical results}

\subsubsection{Determination of the initial width $\eps_0$}

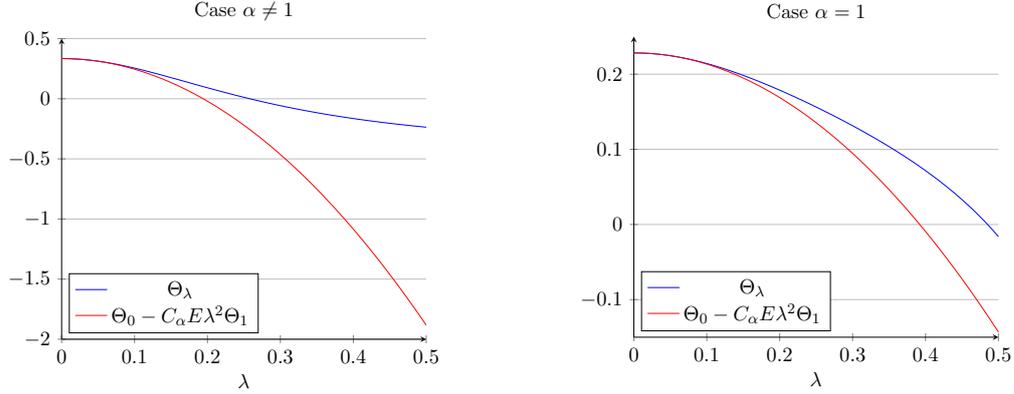
\begin{figure}
		\begin{center}
	\begin{minipage}[c]{0.45\linewidth}
		\begin{tikzpicture}[scale=0.7]
	\begin{axis}[ ymajorgrids, axis x line=bottom, axis y line = left,
		xlabel={$\lambda$}, 
		title={Case $\alpha\neq 1$},
		ymin=-2, ymax=0.5,
		legend entries={$\Theta_\lambda$, $\Theta_0 - C_{\alpha}E\lambda^2\Theta_1$}, legend style={at={(0.02,0.02)}, anchor=south west}
		]
		\addplot+[mark=none] table[x index=0, y index=1]{Data-Plot_Theta-case_2-results.txt};
		\addplot+[mark=none] table[x index =0, y index=2]{Data-Plot_Theta-case_2-results.txt};
	\end{axis}
\end{tikzpicture}
	\end{minipage} 
	\begin{minipage}[c]{0.45\linewidth}
\begin{tikzpicture}[scale=0.7]
	\begin{axis}[ ymajorgrids, axis x line=bottom, axis y line = left,
		xlabel={$\lambda$}, 
		title={Case $\alpha = 1$},
		ymin=-0.15, ymax=0.25,
		legend entries={$\Theta_\lambda$, $\Theta_0 - C_{\alpha}E\lambda^2\Theta_1$}, legend style={at={(0.02,0.02)}, anchor=south west}
		]
		\addplot+[mark=none] table[x index=0, y index=1]{Data-Plot_Theta-case_4-results.txt};
		\addplot+[mark=none] table[x index =0, y index=2]{Data-Plot_Theta-case_4-results.txt};
	\end{axis}
\end{tikzpicture}
	\end{minipage}
\end{center}
	\caption{Example of function $\Theta_{\lambda}$ and its asymptotic development $\Theta_0 - C_{j} E\lambda^2\Theta_1$ plotted against $\lambda$, for the set of parameters $N=0.3, R_c = 0.1, \delta=1, E=10$, and with $\bar\nu_b = 0.1$ (case $\alpha\neq 1$, left) and $\bar \nu_b = 1$ (case $\alpha=1$, right).}
\label{Fig:Theta_lambda}
	\end{figure}

In order to ensure the stability of our numerical method, we choose to determine a small initial value $\eps_0$ that guarantees that, for each tested set of parameters, the computed value of $\eps(t)$ decreases during the simulation, which means by the ODE $\eps'(t)=-\Theta_{\lambda(t)}\, \eps(t)^3$  that the function $\Theta_{\lambda(t)}$ should remain positive. To this aim, we take advantage of the asymptotic development of function $\Theta_\lambda$ as $\lambda$ goes to zero, given in Corollary~\ref{cor_alpha_neq_1} :
\[
\Theta_\lambda=\Theta_0 - C_{j}E\lambda^2 \Theta_1+O(\lambda^4)\, ,
\]
with $j=\alpha$ if $\alpha\neq 1$, $j=N$ if $\alpha = 1$.
As a validation of the above formula, we have plotted in Fig.~\ref{Fig:Theta_lambda} the exact quantity $\Theta_{\lambda}$ (given by~\eqref{theta1_alpha_neq_1}) and its approximation $\Theta_0-C_{j}E\lambda^2 \Theta_1$ as functions of $\lambda$, using 2 different sets of parameters corresponding respectively to a case where $\alpha\neq 1$ and $\alpha=1$. We have also computed numerically the second derivative $\Big(\frac{\partial^2\Theta_{\lambda}}{\partial \lambda^2}\Big)_{|\lambda=0}$ and, in each case, given the expected value $-2C_{j}E\Theta_1$ (see Table~\ref{Num:SecDer}).

\begin{table}
	\caption{\label{Num:SecDer} Second derivative of $\lambda\mapsto\Theta_{\lambda}$ at $\lambda = 0$, computed using the exact formula of $\Theta_\lambda$ or the asymptotic development $\Theta_0 - C_j E \lambda^2 \Theta_1$ , for the set of parameters given in Fig.~\ref{Fig:Theta_lambda}.}
	\begin{center}
\begin{tabular}{|c||c|c|}
	\hline
	$\bar\nu_b$ & $\Big(\frac{\partial^2\Theta_{\lambda}}{\partial \lambda^2}\Big)_{|\lambda=0}$& $-2C_jE\Theta_1$\\
	\hline\hline
	$0.1$ & -17.7344 & -17.7346 \\
	\hline
	$1$ & -2.97108 &-2.97109\\
	\hline
	\end{tabular}
\end{center}
\end{table}

For each set of parameters, we use the approximation $\Theta_\lambda \approx\Theta_0 - C_{\alpha} E\lambda^2\Theta_1 $ to estimate an initial value $\eps_0$ such that the associated value of $\lambda$, given by $\lambda_0=(\kappa\eps_0)^{1/2}$ (see~\eqref{Num:def_lambda}), satisfies $\Theta_{\lambda_0}>0$. Keeping in mind that $\eps_0$ should be small so that the Reynolds equation gives a good approximation of the pressure, we introduce a threshold $\eps_{max}$ and set
\[
\eps_0 = \min\left(\frac{\Theta_0}{\kappa C_j E\, \Theta_1},\eps_{max}\right)\, .
\]
We check a posteriori that the solution $\eps(t)$ to the ODE~\eqref{Num:ODE-eps-1} is indeed a decreasing function of time.

\subsubsection{Influence of parameters $N,R_c,\bar\nu_b,\delta,E$ }

Since the model under study depends on many parameters, in order to perform comparisons, we have chosen to unify the presentation of the numerical results by plotting the half-life time $\Thalf$ as a function of $N\in [0,0.7]$ (which ensures that condition $N^2\leq 1/2$ is fulfilled), after normalization by its value for $N=0$, for different values of $R_c\in\{0.025,0.05,0.1,0.2\}$, using various sets of parameters $\bar\nu_b, \delta, E$.

\paragraph{Influence of $\bar \nu_b$}

\begin{figure}
	
	\begin{center}
		\begin{minipage}[c]{0.45\linewidth}
			\begin{tikzpicture}[scale=0.7]
				\begin{axis}[ axis x line=bottom, axis y line = left,
					xlabel={$N$}, ylabel={$\Thalf$},
					title={$\bar \nu_b = 0.05$},
					legend entries={$R_c = 0.025$, $R_c = 0.05$, $R_c = 0.1$, $R_c = 0.2$}, legend style={at={(0.02,0.02)}, anchor=south west}
					]
					\addplot+[mark=none] table[x index=0, y index=1]{Data-Plot_N-E=0-nub=0-05-Rc=0-025-Results.txt};
					\addplot+[mark=none] table[x index =0, y index=1]{Data-Plot_N-E=0-nub=0-05-Rc=0-05-Results.txt};
					\addplot+[mark=none] table[x index =0, y index=1]{Data-Plot_N-E=0-nub=0-05-Rc=0-1-Results.txt};
					\addplot+[mark=none] table[x index =0, y index=1]{Data-Plot_N-E=0-nub=0-05-Rc=0-2-Results.txt};				
				\end{axis}
			\end{tikzpicture}
		\end{minipage}
		\hfill
		\begin{minipage}[c]{0.45\linewidth}
			\begin{tikzpicture}[scale=0.7]
				\begin{axis}[ axis x line=bottom, axis y line = left,
					xlabel={$N$}, ylabel={$\Thalf$},
					title={$\bar \nu_b = 0.1$},
					legend entries={$R_c = 0.025$, $R_c = 0.05$, $R_c = 0.1$, $R_c = 0.2$}, legend style={at={(0.02,0.02)}, anchor=south west}
					]
					\addplot+[mark=none] table[x index=0, y index=1]{Data-Plot_N-E=0-nub=0-1-Rc=0-025-Results.txt};
					\addplot+[mark=none] table[x index =0, y index=1]{Data-Plot_N-E=0-nub=0-1-Rc=0-05-Results.txt};
					\addplot+[mark=none] table[x index =0, y index=1]{Data-Plot_N-E=0-nub=0-1-Rc=0-1-Results.txt};
					\addplot+[mark=none] table[x index =0, y index=1]{Data-Plot_N-E=0-nub=0-1-Rc=0-2-Results.txt};				
				\end{axis}
			\end{tikzpicture}
		\end{minipage} 
	\end{center}

	\begin{center}
		\begin{minipage}[c]{0.45\linewidth}
			\begin{tikzpicture}[scale=0.7]
				\begin{axis}[ axis x line=bottom, axis y line = left,
					xlabel={$N$}, ylabel={$\Thalf$},
					title={$\bar \nu_b = 0.2$},
					legend entries={$R_c = 0.025$, $R_c = 0.05$, $R_c = 0.1$, $R_c = 0.2$}, legend style={at={(0.02,0.02)}, anchor=south west}
					]
					\addplot+[mark=none] table[x index=0, y index=1]{Data-Plot_N-E=0-nub=0-2-Rc=0-025-Results.txt};
					\addplot+[mark=none] table[x index =0, y index=1]{Data-Plot_N-E=0-nub=0-2-Rc=0-05-Results.txt};
					\addplot+[mark=none] table[x index =0, y index=1]{Data-Plot_N-E=0-nub=0-2-Rc=0-1-Results.txt};
					\addplot+[mark=none] table[x index =0, y index=1]{Data-Plot_N-E=0-nub=0-2-Rc=0-2-Results.txt};				
				\end{axis}
			\end{tikzpicture}
		\end{minipage}
		\hfill
		\begin{minipage}[c]{0.45\linewidth}
			\begin{tikzpicture}[scale=0.7]
				\begin{axis}[ axis x line=bottom, axis y line = left,
					xlabel={$N$}, ylabel={$\Thalf$},
					title={$\bar \nu_b = 0.4$},
					legend entries={$R_c = 0.025$, $R_c = 0.05$, $R_c = 0.1$, $R_c = 0.2$}, legend style={at={(0.02,0.02)}, anchor=south west}
					]
					\addplot+[mark=none] table[x index=0, y index=1]{Data-Plot_N-E=0-nub=0-4-Rc=0-025-Results.txt};
					\addplot+[mark=none] table[x index =0, y index=1]{Data-Plot_N-E=0-nub=0-4-Rc=0-05-Results.txt};
					\addplot+[mark=none] table[x index =0, y index=1]{Data-Plot_N-E=0-nub=0-4-Rc=0-1-Results.txt};
					\addplot+[mark=none] table[x index =0, y index=1]{Data-Plot_N-E=0-nub=0-4-Rc=0-2-Results.txt};				
				\end{axis}
			\end{tikzpicture}
		\end{minipage} 
	\end{center}
	\caption{ $\Thalf$ plotted against $N$, for $\delta=1$, $E=0$ and different values of $R_c\in \{0.025, 0.05, 0.1, 0.2\}$, with $\bar\nu_b\in \{0.05,0.1,0.2,0.4\}$. }
	\label{Fig:Timeratio_over_N-global-E=0}
\end{figure}
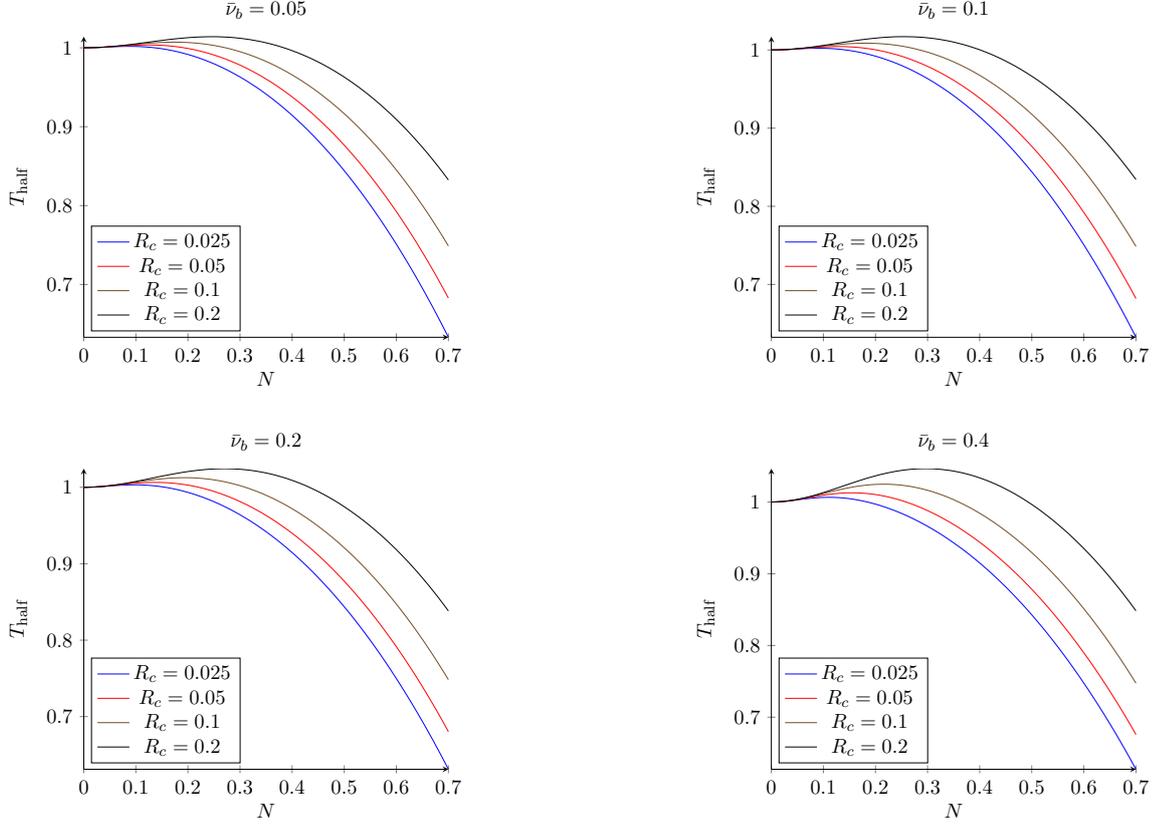

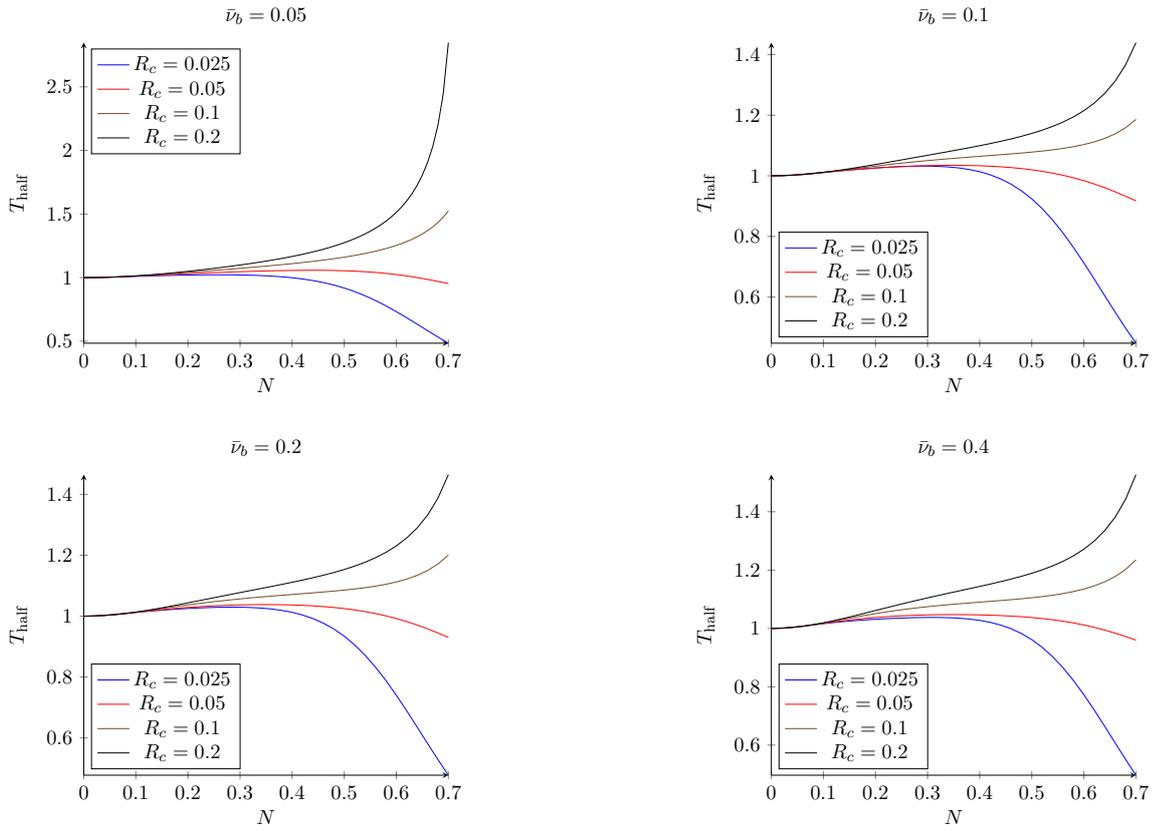
\begin{figure}

	\begin{center}
	\begin{minipage}[c]{0.45\linewidth}
		\begin{tikzpicture}[scale=0.7]
			\begin{axis}[ axis x line=bottom, axis y line = left,
				xlabel={$N$}, ylabel={$\Thalf$},
				title={$\bar \nu_b = 0.05$},
				legend entries={$R_c = 0.025$, $R_c = 0.05$, $R_c = 0.1$, $R_c = 0.2$}, legend style={at={(0.02,0.98)}, anchor=north west}
				]
				\addplot+[mark=none] table[x index=0, y index=1]{Data-Plot_N-E=10-nub=0-05-Rc=0-025-Results.txt};
				\addplot+[mark=none] table[x index =0, y index=1]{Data-Plot_N-E=10-nub=0-05-Rc=0-05-Results.txt};
				\addplot+[mark=none] table[x index =0, y index=1]{Data-Plot_N-E=10-nub=0-05-Rc=0-1-Results.txt};
				\addplot+[mark=none] table[x index =0, y index=1]{Data-Plot_N-E=10-nub=0-05-Rc=0-2-Results.txt};				
			\end{axis}
		\end{tikzpicture}
	\end{minipage}
	\hfill
	\begin{minipage}[c]{0.45\linewidth}
		\begin{tikzpicture}[scale=0.7]
			\begin{axis}[ axis x line=bottom, axis y line = left,
				xlabel={$N$}, ylabel={$\Thalf$},
				title={$\bar \nu_b = 0.1$},
				legend entries={$R_c = 0.025$, $R_c = 0.05$, $R_c = 0.1$, $R_c = 0.2$}, legend style={at={(0.02,0.02)}, anchor=south west}
				]
				\addplot+[mark=none] table[x index=0, y index=1]{Data-Plot_N-E=10-nub=0-1-Rc=0-025-Results.txt};
				\addplot+[mark=none] table[x index =0, y index=1]{Data-Plot_N-E=10-nub=0-1-Rc=0-05-Results.txt};
				\addplot+[mark=none] table[x index =0, y index=1]{Data-Plot_N-E=10-nub=0-1-Rc=0-1-Results.txt};
				\addplot+[mark=none] table[x index =0, y index=1]{Data-Plot_N-E=10-nub=0-1-Rc=0-2-Results.txt};				
			\end{axis}
		\end{tikzpicture}
	\end{minipage} 
\end{center}

	\begin{center}
		\begin{minipage}[c]{0.45\linewidth}
		\begin{tikzpicture}[scale=0.7]
	\begin{axis}[ axis x line=bottom, axis y line = left,
		xlabel={$N$}, ylabel={$\Thalf$},
		title={$\bar \nu_b = 0.2$},
		legend entries={$R_c = 0.025$, $R_c = 0.05$, $R_c = 0.1$, $R_c = 0.2$}, legend style={at={(0.02,0.02)}, anchor=south west}
		]
		\addplot+[mark=none] table[x index=0, y index=1]{Data-Plot_N-E=10-nub=0-2-Rc=0-025-Results.txt};
		\addplot+[mark=none] table[x index =0, y index=1]{Data-Plot_N-E=10-nub=0-2-Rc=0-05-Results.txt};
		\addplot+[mark=none] table[x index =0, y index=1]{Data-Plot_N-E=10-nub=0-2-Rc=0-1-Results.txt};
		\addplot+[mark=none] table[x index =0, y index=1]{Data-Plot_N-E=10-nub=0-2-Rc=0-2-Results.txt};				
	\end{axis}
\end{tikzpicture}
		\end{minipage}
		\hfill
		\begin{minipage}[c]{0.45\linewidth}
		\begin{tikzpicture}[scale=0.7]
	\begin{axis}[ axis x line=bottom, axis y line = left,
		xlabel={$N$}, ylabel={$\Thalf$},
		title={$\bar \nu_b = 0.4$},
		legend entries={$R_c = 0.025$, $R_c = 0.05$, $R_c = 0.1$, $R_c = 0.2$}, legend style={at={(0.02,0.02)}, anchor=south west}
		]
		\addplot+[mark=none] table[x index=0, y index=1]{Data-Plot_N-E=10-nub=0-4-Rc=0-025-Results.txt};
		\addplot+[mark=none] table[x index =0, y index=1]{Data-Plot_N-E=10-nub=0-4-Rc=0-05-Results.txt};
		\addplot+[mark=none] table[x index =0, y index=1]{Data-Plot_N-E=10-nub=0-4-Rc=0-1-Results.txt};
		\addplot+[mark=none] table[x index =0, y index=1]{Data-Plot_N-E=10-nub=0-4-Rc=0-2-Results.txt};				
	\end{axis}
\end{tikzpicture}
		\end{minipage} 
	\end{center}
	\caption{ $\Thalf$ plotted against $N$, for $\delta=1$, $E=10$ and different values of $R_c\in \{0.025, 0.05, 0.1, 0.2\}$, with $\bar\nu_b\in \{0.05,0.1,0.2,0.4\}$. }
	\label{Fig:Timeratio_over_N-global-E=10}
\end{figure}

We have plotted in Fig.~\ref{Fig:Timeratio_over_N-global-E=0} and Fig.~\ref{Fig:Timeratio_over_N-global-E=10} the results obtained with $E=0$ and $E=10$ respectively, considering different values of parameter $\bar\nu_b\in\{0.05,0.1,0.2,0.4\}$ and for a fixed $\delta=1$. It appears that modifying $\bar\nu_b$ does not have much of an impact of the computed value of $\Thalf$, at least qualititatively. Consequently, we have decided to fix this parameter and impose $\bar\nu_b=0.1$ in the rest of the simulations.

\paragraph{Influence of $N,R_c$}

On the opposite, $\Thalf$ is very sensitive to the couple of parameters $(N,R_c)$. It appears that increasing $N$ from the initial value $N=0$ leads, at first, to a slow increasing of $\Thalf$, up to a certain value of $N$. Then, the dependence of $\Thalf$ on $N$ is strongly affected by the value of $R_c$. For instance, in the case $E=10$ (Fig.~\ref{Fig:Timeratio_over_N-global-E=10}), the maximal value $R_c=0.2$ leads to a gradual increasing of the slope of the curve, up to $N=0.7$, whereas the minimal value $R_c=0.025$ corresponds to a reduction of $\Thalf$ as $N$ increases, ending up with a division by a factor $2$ with respect to the initial value for $N=0$.

\paragraph{Influence of $E$}

\begin{figure}
	\begin{center}
		\begin{minipage}[c]{0.32\linewidth}
			\begin{tikzpicture}[scale=0.5]
				\begin{axis}[ axis x line=bottom, axis y line = left,
					xlabel={$y_1$}, ylabel={$\Psi(y_1)$},
					title={$V$-shape},
					ymin=0, ymax=0.6,
					]
					\addplot+[mark=none] table[x index=0, y index=1]{Data-Plot_riblet-sawtooth.txt};
				\end{axis}
			\end{tikzpicture}
		\end{minipage}
		\begin{minipage}[c]{0.32\linewidth}
	\begin{tikzpicture}[scale=0.5]
		\begin{axis}[ axis x line=bottom, axis y line = left,
			xlabel={$y_1$}, ylabel={$\Psi(y_1)$},
			title={$U$-shape},
			ymin=0, ymax=0.35,
			]
			\addplot+[mark=none] table[x index=0, y index=1]{Data-Plot_riblet-scallop.txt};
		\end{axis}
	\end{tikzpicture}
		\end{minipage} 
		\begin{minipage}[c]{0.32\linewidth}
	\begin{tikzpicture}[scale=0.5]
		\begin{axis}[ axis x line=bottom, axis y line = left,
			xlabel={$y_1$}, ylabel={$\Psi(y_1)$},
			title={blade},
			ymin=0, ymax=0.18,
			]
			\addplot+[mark=none] table[x index=0, y index=1]{Data-Plot_riblet-blade.txt};
		\end{axis}
	\end{tikzpicture}
\end{minipage} 
	\end{center}
	\caption{Examples of function $\Psi$ associated, from left to right, with $V$-shape, $U$-shape and blade riblets, normalized by~\eqref{Psi-normalized}.}
	\label{Fig:riblet_psi}
\end{figure}
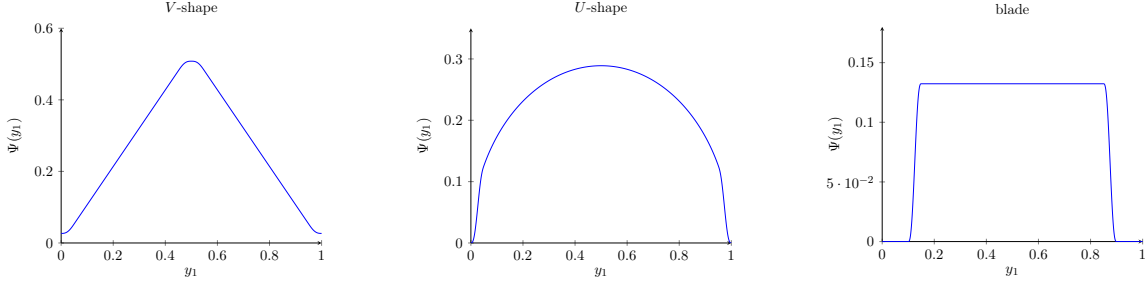

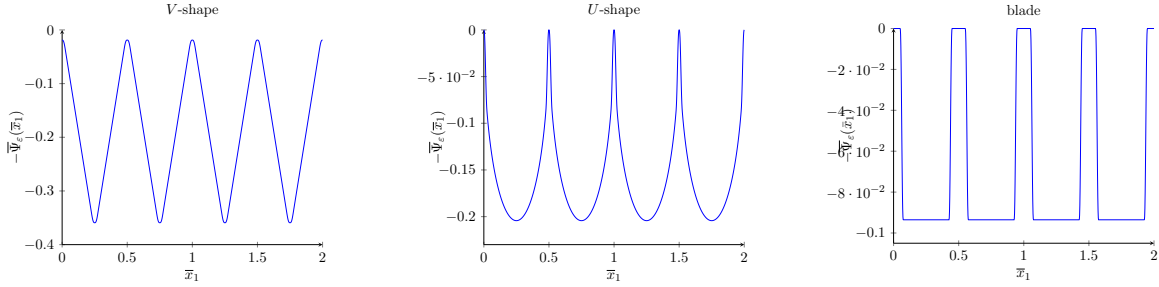
\begin{figure}
	\begin{center}
		\begin{minipage}[c]{0.32\linewidth}
			\begin{tikzpicture}[scale=0.5]
				\begin{axis}[  axis x line=bottom, axis y line = left,
					xlabel={$\overline x_1$}, ylabel={$-\overline{\Psi}_{\eps}(\overline x_1)$},
					title={$V$-shape},
					ymin=-0.4, ymax=0,
					]
					\addplot+[mark=none] table[x index=0, y index=1]{Data-Plot_riblet-sawtooth_eps.txt};
				\end{axis}
			\end{tikzpicture}
		\end{minipage}
		\begin{minipage}[c]{0.32\linewidth}
			\begin{tikzpicture}[scale=0.5]
				\begin{axis}[ axis x line=bottom, axis y line = left,
					xlabel={$\overline x_1$}, ylabel={$-\overline{\Psi}_{\eps}(\overline x_1)$},
					title={$U$-shape},
					ymin=-0.23, ymax=0,
					]
					\addplot+[mark=none] table[x index=0, y index=1]{Data-Plot_riblet-scallop_eps.txt};
				\end{axis}
			\end{tikzpicture}
		\end{minipage} 
		\begin{minipage}[c]{0.32\linewidth}
			\begin{tikzpicture}[scale=0.5]
				\begin{axis}[ axis x line=bottom, axis y line = left,
					xlabel={$\overline x_1$}, ylabel={$-\overline\Psi_\eps(\bar x_1)$},
					title={blade},
					ymin=-0.105, ymax=0,
					]
					\addplot+[mark=none] table[x index=0, y index=1]{Data-Plot_riblet-blade_eps.txt};
				\end{axis}
			\end{tikzpicture}
		\end{minipage} 
	\end{center}
	\caption{Examples of rough geometries associated with functions $\overline{\Psi}_{\eps}(\overline x_1)$ defined by~\eqref{Num:psi_eps_bar} and~\eqref{Num:def_Lambda-M}, with $L=2, M=0.5$ and $\Lambda = M^{3/2}$. Each riblet profile is described by the corresponding function $\Psi$ plotted in Fig.~\ref{Fig:riblet_psi}.}
	\label{Fig:riblet_psi_eps}
\end{figure}

In order to estimate the possible values of $E$ encountered in practical applications, we consider the example of three riblet profiles that are often used in the engineering literature, namely, the $V$-shape, the $U$-shape and the blade riblets. We have plotted in Fig.~\ref{Fig:riblet_psi} the corresponding $\Psi$ functions, normalized by~\eqref{Psi-normalized}, and in Fig.~\ref{Fig:riblet_psi_eps}, a possible rough geometry defined as the graph of the function $-\overline\Psi_{\eps}$, defined by~\eqref{Num:psi_eps_bar} and~\eqref{Num:def_Lambda-M}.
To compute the energy $E$ associated with each of these riblet profiles, we have solved system~\eqref{system_phi_1} (in case $\lambda=1$) by a finite element method using FreeFem++ software \cite{FREEFEM}, implementing a Taylor-Hood approximation for the velocity-pressure pair, \emph{i.e.}, $P_2$ elements for the velocity field and $P_1$ elements for the pressure. The results are summarized in Table~\ref{Num:computation-E}.

Even though the computed values of $E$ associated with the previous examples of riblets are of order $10$ to $100$, it turns out that the results of our simulations are very stable when $E$ exceeds $10$, so we have chosen to represent in Fig.~\ref{Fig:Dependence-on-E-delta=1-custom} the range of values of $E$ that produces the most significant changes in the behaviour of the model, which is $E\in [0,10]$ for $\delta=1$ and $\bar\nu_b=0.1$. In the case $E=0$, \emph{i.e.} in the absence of roughness, raising $N$ from $N=0$ results at first in a slight increase of $\Thalf$. Then, this behaviour reverses and for larger values of $N$, we observe a significant reduction of $\Thalf$, for all the tested values of $R_c$. Raising the value of $E$ produces a visible change in the model, whose behaviour gets more and more dependent on parameter $R_c$. As $E$ increases, the function $N\mapsto \Thalf$ becomes monotonic for the highest values of $R_c$ ($R_c=0.1$ and $R_c=0.2$). This means that for certain micropolar fluids, the roughness of the upper plate may contribute to enhance the performance of the squeeze-film bearing, in the sense that $\Thalf$ becomes larger.

\paragraph{Influence of $\delta$}

We conclude this numerical study by investigating the impact of the slip length $\delta$ on the behaviour of the bearing. We have used values from $\delta=0.7$ to $\delta=10$, a range for which a significant change occurs in the simulations (see Fig.~\ref{Fig:Dependence-on-delta-E=5}). We observe that small values of $\delta$ favor an enhancement of the performance of the bearing, for $R_c\in \{0.1,0.2\}$. On the opposite, large values of $\delta$ such as $\delta=10$ will typically reduce the performance of the micropolar fluid lubrication, with respect to lubrication with a Newtonian fluid, with the exception of the largest value or parameter $R_c=0.2$, where a slight increase of $\Thalf$ occurs.

\begin{figure}
	
	\begin{center}
		\begin{minipage}[c]{0.45\linewidth}
			\begin{tikzpicture}[scale=0.7]
				\begin{axis}[ axis x line=bottom, axis y line = left,
					xlabel={$N$}, ylabel={$\Thalf$},
					title={$E = 0$},
					legend entries={$R_c = 0.025$, $R_c = 0.05$, $R_c = 0.1$, $R_c = 0.2$}, legend style={at={(0.02,0.02)}, anchor=south west}
					]
					\addplot+[mark=none] table[x index=0, y index=1]{Data-Plot_N-E_custom-E=0-Res-global.txt};
					\addplot+[mark=none] table[x index =0, y index=2]{Data-Plot_N-E_custom-E=0-Res-global.txt};
					\addplot+[mark=none] table[x index =0, y index=3]{Data-Plot_N-E_custom-E=0-Res-global.txt};
					\addplot+[mark=none] table[x index =0, y index=4]{Data-Plot_N-E_custom-E=0-Res-global.txt};				
				\end{axis}
			\end{tikzpicture}
		\end{minipage}
		\hfill
		\begin{minipage}[c]{0.45\linewidth}
			\begin{tikzpicture}[scale=0.7]
	\begin{axis}[ axis x line=bottom, axis y line = left,
		xlabel={$N$}, ylabel={$\Thalf$},
		title={$E = 1$},
		legend entries={$R_c = 0.025$, $R_c = 0.05$, $R_c = 0.1$, $R_c = 0.2$}, legend style={at={(0.02,0.02)}, anchor=south west}
		]
		\addplot+[mark=none] table[x index=0, y index=1]{Data-Plot_N-E_custom-E=1-Res-global.txt};
		\addplot+[mark=none] table[x index =0, y index=2]{Data-Plot_N-E_custom-E=1-Res-global.txt};
		\addplot+[mark=none] table[x index =0, y index=3]{Data-Plot_N-E_custom-E=1-Res-global.txt};
		\addplot+[mark=none] table[x index =0, y index=4]{Data-Plot_N-E_custom-E=1-Res-global.txt};				
	\end{axis}
\end{tikzpicture}
		\end{minipage} 
	\end{center}

		\begin{center}
		\begin{minipage}[c]{0.45\linewidth}
			\begin{tikzpicture}[scale=0.7]
				\begin{axis}[ axis x line=bottom, axis y line = left,
					xlabel={$N$}, ylabel={$\Thalf$},
					title={$E = 3$},
					legend entries={$R_c = 0.025$, $R_c = 0.05$, $R_c = 0.1$, $R_c = 0.2$}, legend style={at={(0.02,0.02)}, anchor=south west}
					]
					\addplot+[mark=none] table[x index=0, y index=1]{Data-Plot_N-E_custom-E=3-Res-global.txt};
					\addplot+[mark=none] table[x index =0, y index=2]{Data-Plot_N-E_custom-E=3-Res-global.txt};
					\addplot+[mark=none] table[x index =0, y index=3]{Data-Plot_N-E_custom-E=3-Res-global.txt};
					\addplot+[mark=none] table[x index =0, y index=4]{Data-Plot_N-E_custom-E=3-Res-global.txt};				
				\end{axis}
			\end{tikzpicture}
		\end{minipage}
		\hfill
		\begin{minipage}[c]{0.45\linewidth}
			\begin{tikzpicture}[scale=0.7]
				\begin{axis}[ axis x line=bottom, axis y line = left,
					xlabel={$N$}, ylabel={$\Thalf$},
					title={$E = 5$},
					legend entries={$R_c = 0.025$, $R_c = 0.05$, $R_c = 0.1$, $R_c = 0.2$}, legend style={at={(0.02,0.02)}, anchor=south west}
					]
					\addplot+[mark=none] table[x index=0, y index=1]{Data-Plot_N-E_custom-E=5-Res-global.txt};
					\addplot+[mark=none] table[x index =0, y index=2]{Data-Plot_N-E_custom-E=5-Res-global.txt};
					\addplot+[mark=none] table[x index =0, y index=3]{Data-Plot_N-E_custom-E=5-Res-global.txt};
					\addplot+[mark=none] table[x index =0, y index=4]{Data-Plot_N-E_custom-E=5-Res-global.txt};				
				\end{axis}
			\end{tikzpicture}
		\end{minipage} 
	\end{center}

			\begin{center}
		\begin{minipage}[c]{0.45\linewidth}
			\begin{tikzpicture}[scale=0.7]
				\begin{axis}[ axis x line=bottom, axis y line = left,
					xlabel={$N$}, ylabel={$\Thalf$},
					title={$E = 7$},
					legend entries={$R_c = 0.025$, $R_c = 0.05$, $R_c = 0.1$, $R_c = 0.2$}, legend style={at={(0.02,0.02)}, anchor=south west}
					]
					\addplot+[mark=none] table[x index=0, y index=1]{Data-Plot_N-E_custom-E=7-Res-global.txt};
					\addplot+[mark=none] table[x index =0, y index=2]{Data-Plot_N-E_custom-E=7-Res-global.txt};
					\addplot+[mark=none] table[x index =0, y index=3]{Data-Plot_N-E_custom-E=7-Res-global.txt};
					\addplot+[mark=none] table[x index =0, y index=4]{Data-Plot_N-E_custom-E=7-Res-global.txt};				
				\end{axis}
			\end{tikzpicture}
		\end{minipage}
		\hfill
		\begin{minipage}[c]{0.45\linewidth}
			\begin{tikzpicture}[scale=0.7]
				\begin{axis}[ axis x line=bottom, axis y line = left,
					xlabel={$N$}, ylabel={$\Thalf$},
					title={$E = 10$},
					legend entries={$R_c = 0.025$, $R_c = 0.05$, $R_c = 0.1$, $R_c = 0.2$}, legend style={at={(0.02,0.98)}, anchor=north west}
					]
					\addplot+[mark=none] table[x index=0, y index=1]{Data-Plot_N-E_custom-E=10-Res-global.txt};
					\addplot+[mark=none] table[x index =0, y index=2]{Data-Plot_N-E_custom-E=10-Res-global.txt};
					\addplot+[mark=none] table[x index =0, y index=3]{Data-Plot_N-E_custom-E=10-Res-global.txt};
					\addplot+[mark=none] table[x index =0, y index=4]{Data-Plot_N-E_custom-E=10-Res-global.txt};				
				\end{axis}
			\end{tikzpicture}
		\end{minipage} 
	\end{center}

	\caption{ $\Thalf$ plotted against $N$, for $\bar\nu_b=0.1$, $R_c\in \{0.025, 0.05, 0.1, 0.2\}$ and $\delta=1$. From left to right and top to bottom: $E=0$, $E=1$, $E=3$, $E=5$, $E=7$ and $E=10$. }
	\label{Fig:Dependence-on-E-delta=1-custom}
\end{figure}
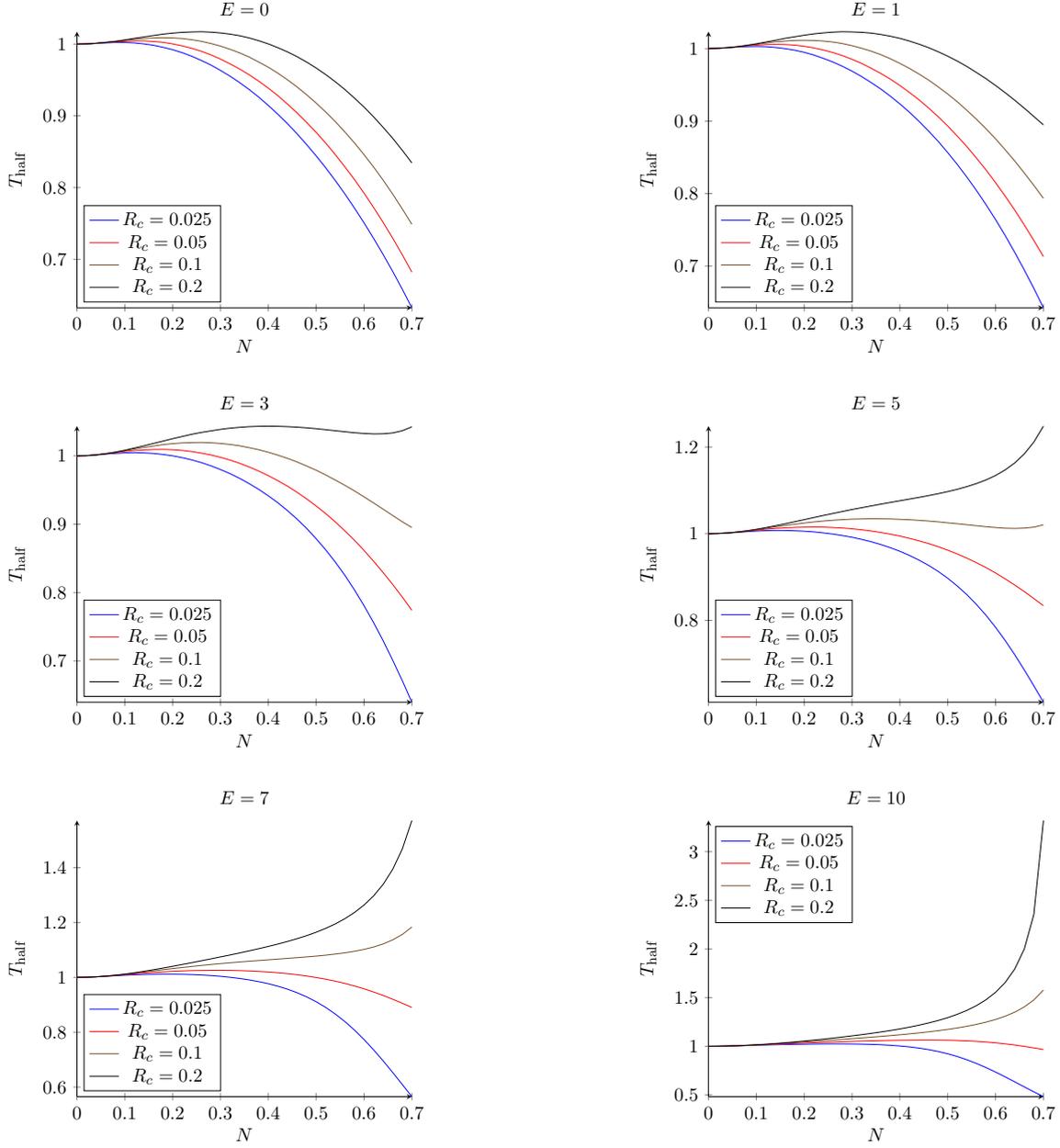


\begin{figure}
	
	\begin{center}
		\begin{minipage}[c]{0.45\linewidth}
			\begin{tikzpicture}[scale=0.7]
				\begin{axis}[ axis x line=bottom, axis y line = left,
					xlabel={$N$}, ylabel={$\Thalf$},
					title={$\delta=0.7$},
					ymin=0.9, ymax=1.4,
					legend entries={$R_c = 0.025$, $R_c = 0.05$, $R_c = 0.1$, $R_c = 0.2$}, legend style={at={(0.02,0.98)}, anchor=north west}
					]
					\addplot+[mark=none] table[x index=0, y index=1]{Data-Plot_N-E=5-nub=0-1-delta_custom-delta=0-7-Res-global.txt};
					\addplot+[mark=none] table[x index=0, y index=2]{Data-Plot_N-E=5-nub=0-1-delta_custom-delta=0-7-Res-global.txt};
					\addplot+[mark=none] table[x index=0, y index=3]{Data-Plot_N-E=5-nub=0-1-delta_custom-delta=0-7-Res-global.txt};
					\addplot+[mark=none] table[x index=0, y index=4]{Data-Plot_N-E=5-nub=0-1-delta_custom-delta=0-7-Res-global.txt};
	
				\end{axis}
			\end{tikzpicture}
		\end{minipage}
		\hfill
		\begin{minipage}[c]{0.45\linewidth}
				\begin{tikzpicture}[scale=0.7]
		\begin{axis}[ axis x line=bottom, axis y line = left,
			xlabel={$N$}, ylabel={$\Thalf$},
			title={$\delta=0.8$},
			legend entries={$R_c = 0.025$, $R_c = 0.05$, $R_c = 0.1$, $R_c = 0.2$}, legend style={at={(0.02,0.02)}, anchor=south west}
			]
			\addplot+[mark=none] table[x index=0, y index=1]{Data-Plot_N-E=5-nub=0-1-delta_custom-delta=0-8-Res-global.txt};
			\addplot+[mark=none] table[x index=0, y index=2]{Data-Plot_N-E=5-nub=0-1-delta_custom-delta=0-8-Res-global.txt};
			\addplot+[mark=none] table[x index=0, y index=3]{Data-Plot_N-E=5-nub=0-1-delta_custom-delta=0-8-Res-global.txt};
			\addplot+[mark=none] table[x index=0, y index=4]{Data-Plot_N-E=5-nub=0-1-delta_custom-delta=0-8-Res-global.txt};
			
		\end{axis}
	\end{tikzpicture}
		\end{minipage} 
	\end{center}

	\begin{center}
	\begin{minipage}[c]{0.45\linewidth}
		\begin{tikzpicture}[scale=0.7]
			\begin{axis}[ axis x line=bottom, axis y line = left,
				xlabel={$N$}, ylabel={$\Thalf$},
				title={$\delta=1$},
				legend entries={$R_c = 0.025$, $R_c = 0.05$, $R_c = 0.1$, $R_c = 0.2$}, legend style={at={(0.02,0.02)}, anchor=south west}
				]
				\addplot+[mark=none] table[x index=0, y index=1]{Data-Plot_N-E=5-nub=0-1-delta_custom-delta=1-Res-global.txt};
				\addplot+[mark=none] table[x index=0, y index=2]{Data-Plot_N-E=5-nub=0-1-delta_custom-delta=1-Res-global.txt};
				\addplot+[mark=none] table[x index=0, y index=3]{Data-Plot_N-E=5-nub=0-1-delta_custom-delta=1-Res-global.txt};
				\addplot+[mark=none] table[x index=0, y index=4]{Data-Plot_N-E=5-nub=0-1-delta_custom-delta=1-Res-global.txt};
				
			\end{axis}
		\end{tikzpicture}
	\end{minipage}
	\hfill
	\begin{minipage}[c]{0.45\linewidth}
		\begin{tikzpicture}[scale=0.7]
			\begin{axis}[ axis x line=bottom, axis y line = left,
				xlabel={$N$}, ylabel={$\Thalf$},
				title={$\delta=1.2$},
				legend entries={$R_c = 0.025$, $R_c = 0.05$, $R_c = 0.1$, $R_c = 0.2$}, legend style={at={(0.02,0.02)}, anchor=south west}
				]
				\addplot+[mark=none] table[x index=0, y index=1]{Data-Plot_N-E=5-nub=0-1-delta_custom-delta=1-2-Res-global.txt};
				\addplot+[mark=none] table[x index=0, y index=2]{Data-Plot_N-E=5-nub=0-1-delta_custom-delta=1-2-Res-global.txt};
				\addplot+[mark=none] table[x index=0, y index=3]{Data-Plot_N-E=5-nub=0-1-delta_custom-delta=1-2-Res-global.txt};
				\addplot+[mark=none] table[x index=0, y index=4]{Data-Plot_N-E=5-nub=0-1-delta_custom-delta=1-2-Res-global.txt};
				
			\end{axis}
		\end{tikzpicture}
	\end{minipage} 
\end{center}

	\begin{center}
	\begin{minipage}[c]{0.45\linewidth}
		\begin{tikzpicture}[scale=0.7]
			\begin{axis}[ axis x line=bottom, axis y line = left,
				xlabel={$N$}, ylabel={$\Thalf$},
				title={$\delta=2$},
				legend entries={$R_c = 0.025$, $R_c = 0.05$, $R_c = 0.1$, $R_c = 0.2$}, legend style={at={(0.02,0.02)}, anchor=south west}
				]
				\addplot+[mark=none] table[x index=0, y index=1]{Data-Plot_N-E=5-nub=0-1-delta_custom-delta=2-Res-global.txt};
				\addplot+[mark=none] table[x index=0, y index=2]{Data-Plot_N-E=5-nub=0-1-delta_custom-delta=2-Res-global.txt};
				\addplot+[mark=none] table[x index=0, y index=3]{Data-Plot_N-E=5-nub=0-1-delta_custom-delta=2-Res-global.txt};
				\addplot+[mark=none] table[x index=0, y index=4]{Data-Plot_N-E=5-nub=0-1-delta_custom-delta=2-Res-global.txt};
				
			\end{axis}
		\end{tikzpicture}
	\end{minipage}
	\hfill
	\begin{minipage}[c]{0.45\linewidth}
		\begin{tikzpicture}[scale=0.7]
			\begin{axis}[ axis x line=bottom, axis y line = left,
				xlabel={$N$}, ylabel={$\Thalf$},
				title={$\delta=10$},
				legend entries={$R_c = 0.025$, $R_c = 0.05$, $R_c = 0.1$, $R_c = 0.2$}, legend style={at={(0.02,0.02)}, anchor=south west}
				]
				\addplot+[mark=none] table[x index=0, y index=1]{Data-Plot_N-E=5-nub=0-1-delta_custom-delta=10-Res-global.txt};
				\addplot+[mark=none] table[x index=0, y index=2]{Data-Plot_N-E=5-nub=0-1-delta_custom-delta=10-Res-global.txt};
				\addplot+[mark=none] table[x index=0, y index=3]{Data-Plot_N-E=5-nub=0-1-delta_custom-delta=10-Res-global.txt};
				\addplot+[mark=none] table[x index=0, y index=4]{Data-Plot_N-E=5-nub=0-1-delta_custom-delta=10-Res-global.txt};
				
			\end{axis}
		\end{tikzpicture}
	\end{minipage} 
\end{center}

	\caption{ $\Thalf$ plotted against $N$, for $\bar\nu_b=0.1$, $R_c\in \{0.025, 0.05, 0.1, 0.2\}$ and $E=5$. From left to right and top to bottom: $\delta=0.7$, $\delta=0.8$, $\delta=1$,$\delta=1.2$, $\delta=2$ and $\delta=10$. }
	\label{Fig:Dependence-on-delta-E=5}
\end{figure}
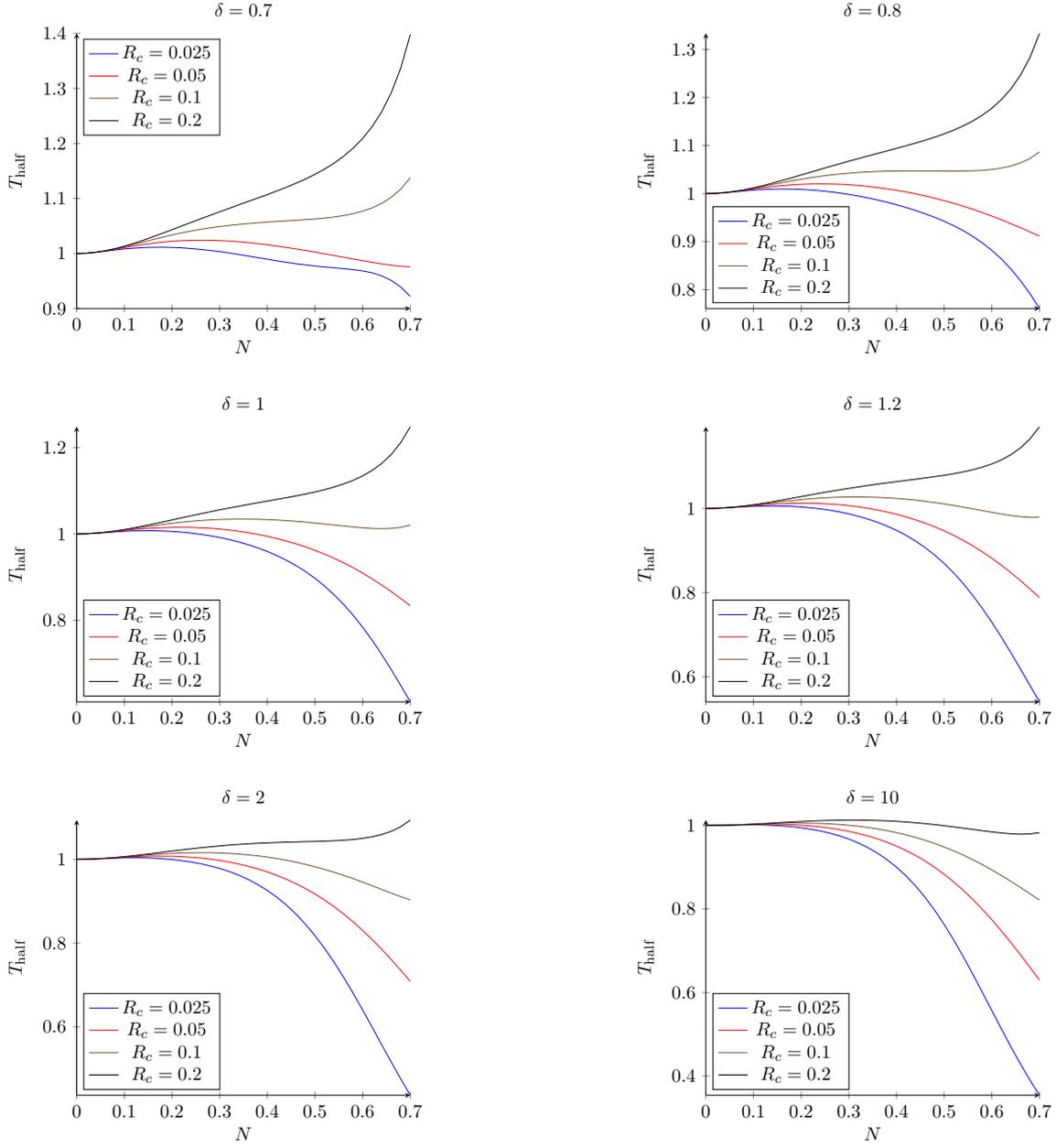


\begin{table}
	\begin{center}
		\begin{tabular}{|c||c|c|c|}
			\hline
			Riblet profile & $V$-shape & $U$-shape& Blade\\
			\hline
			Value of $E$ &12.12 & 62.85 &93.24\\
			\hline
		\end{tabular}
	\end{center}
\caption{\label{Num:computation-E} Value of $E$ associated with three riblet profiles given in Fig.~\ref{Fig:riblet_psi}: the $V$-shape, $U$-shape and blade riblets.}
\end{table}

\subsubsection{Interpretation in terms of pressure using Corollary~\ref{cor_alpha_neq_1}}

As presented in Subsection~\ref{Subsec:DerivationModel}, we have assumed in this study that $\lambda$ and $\eps$ satisfy a relation of the form $\lambda=(\kappa \eps)^{1/2}$, with $\kappa = \Lambda^2/M^3$ (see~\eqref{Num:def_lambda}), where $\eps$ is small. Consequently, $\lambda$ itself is a small parameter, so one can apply the asymptotic developments obtained in Corollary~\ref{cor_alpha_neq_1} to shed a light on the observed behaviour of the model. Indeed, for a small value of $\lambda$, equation~\eqref{Num:S(t)} shows that the vertical velocity of the upper plate is inversely proportional to $\int_{0}^1 p^{\lambda,1}(y_1)\, dy_1$, which can be approximated by
\begin{align*}
\int_{0}^1 p^{\lambda,1}(y_1)\, dy_1 & \approx \int_{0}^1 \Big(p_0(y_1)\, dy_1 + C_\alpha E\lambda^2  p_1(y_1)\Big)\, dy_1  \\
& \approx \int_0^1 (1+C_\alpha E\frac{\Theta_1}{\Theta_0}\lambda^2)p_0(y_1)dy_1\, ,
\end{align*}
where $p_0$, $p_1$ are the respective solutions of~\eqref{reynoldsp0_neq1} and \eqref{reynoldsp1_neq1} with $S=1$.
In particular, at first order in $\lambda$, there holds (in case $\alpha\neq 1$)
\[
\frac{\int_{0}^1 p^{\lambda,1}(y_1)\, dy_1 - \int_{0}^1 p_0(y_1)\, dy_1 }{\int_{0}^1 p_0(y_1)\, dy_1 } \approx C_\alpha E\frac{\Theta_1}{\Theta_0}\lambda^2\, .
\] 

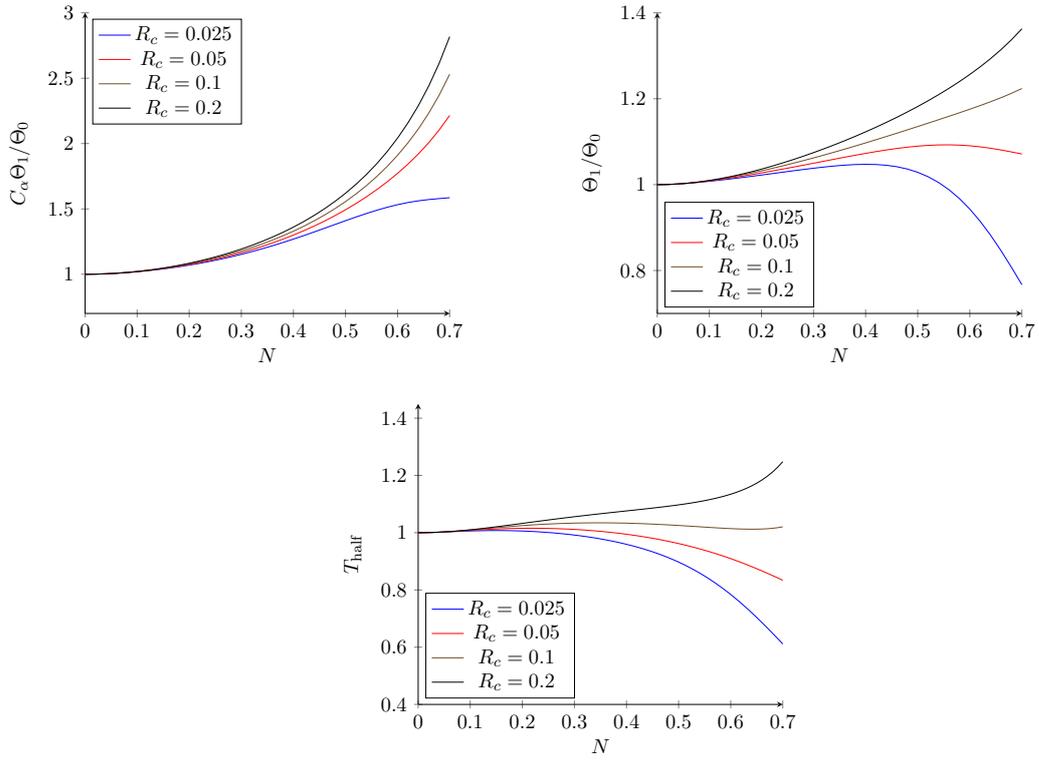
\begin{figure}
	\begin{center}
		\begin{minipage}[c]{0.45\linewidth}
			\begin{tikzpicture}[scale=0.7]
				\begin{axis}[ axis x line=bottom, axis y line = left,
					xlabel={$N$}, ylabel={$C_{\alpha}{\Theta_1}/{\Theta_0}$},
					ymin=0.7, ymax=3,
					legend entries={$R_c = 0.025$, $R_c = 0.05$, $R_c = 0.1$, $R_c = 0.2$}, legend style={at={(0.02,0.98)}, anchor=north west}
					]
					\addplot+[mark=none] table[x index=0, y index=1]{Data-Interp_Plot_N-E=5-delta=1-Rc=0-025-Results.txt};
					\addplot+[mark=none] table[x index =0, y index=1]{Data-Interp_Plot_N-E=5-delta=1-Rc=0-05-Results.txt};
					\addplot+[mark=none] table[x index =0, y index=1]{Data-Interp_Plot_N-E=5-delta=1-Rc=0-1-Results.txt};
					\addplot+[mark=none] table[x index =0, y index=1]{Data-Interp_Plot_N-E=5-delta=1-Rc=0-2-Results.txt};				
				\end{axis}
			\end{tikzpicture}
		\end{minipage} 
		\begin{minipage}[c]{0.45\linewidth}
			\begin{tikzpicture}[scale=0.7]
				\begin{axis}[ axis x line=bottom, axis y line = left,
					xlabel={$N$}, ylabel={${\Theta_1}/{\Theta_0}$},
					ymin=0.7, ymax=1.4,
					legend entries={$R_c = 0.025$, $R_c = 0.05$, $R_c = 0.1$, $R_c = 0.2$}, legend style={at={(0.02,0.02)}, anchor=south west}
					]
					\addplot+[mark=none] table[x index=0, y index=6]{Data-Interp_Plot_N-E=5-delta=1-Rc=0-025-Results.txt};
					\addplot+[mark=none] table[x index =0, y index=6]{Data-Interp_Plot_N-E=5-delta=1-Rc=0-05-Results.txt};
					\addplot+[mark=none] table[x index =0, y index=6]{Data-Interp_Plot_N-E=5-delta=1-Rc=0-1-Results.txt};
					\addplot+[mark=none] table[x index =0, y index=6]{Data-Interp_Plot_N-E=5-delta=1-Rc=0-2-Results.txt};				
				\end{axis}
			\end{tikzpicture}
		\end{minipage}
	\end{center}
	\begin{center}
		\begin{tikzpicture}[scale=0.7]
			\begin{axis}[ axis x line=bottom, axis y line = left,
				xlabel={$N$}, ylabel={$\Thalf$},
				ymin=0.4, ymax=1.45,
				legend entries={$R_c = 0.025$, $R_c = 0.05$, $R_c = 0.1$, $R_c = 0.2$}, legend style={at={(0.02,0.02)}, anchor=south west}
				]
				\addplot+[mark=none] table[x index=0, y index=1]{Data-Plot_N-E=5-nub=0-1-Rc=0-025-Results.txt};
				\addplot+[mark=none] table[x index =0, y index=1]{Data-Plot_N-E=5-nub=0-1-Rc=0-05-Results.txt};
				\addplot+[mark=none] table[x index =0, y index=1]{Data-Plot_N-E=5-nub=0-1-Rc=0-1-Results.txt};
				\addplot+[mark=none] table[x index =0, y index=1]{Data-Plot_N-E=5-nub=0-1-Rc=0-2-Results.txt};				
			\end{axis}
		\end{tikzpicture}
	\end{center}
	
	\caption{$C_{\alpha}{\Theta_1}/{\Theta_0}$ (top, left) and ${\Theta_1}/{\Theta_0}$ (top, right) plotted against $N$, normalized by the value at $N=0$, for $\bar{\nu}_b = 0.1$, $\delta=1$, $E=5$ and different values of $R_c\in \{0.025, 0.05, 0.1, 0.2\}$. Bottom: $\Thalf$ plotted against $N$ for the same sets of parameters.}
	\label{Fig:Interp-N-global-delta=1}
\end{figure}

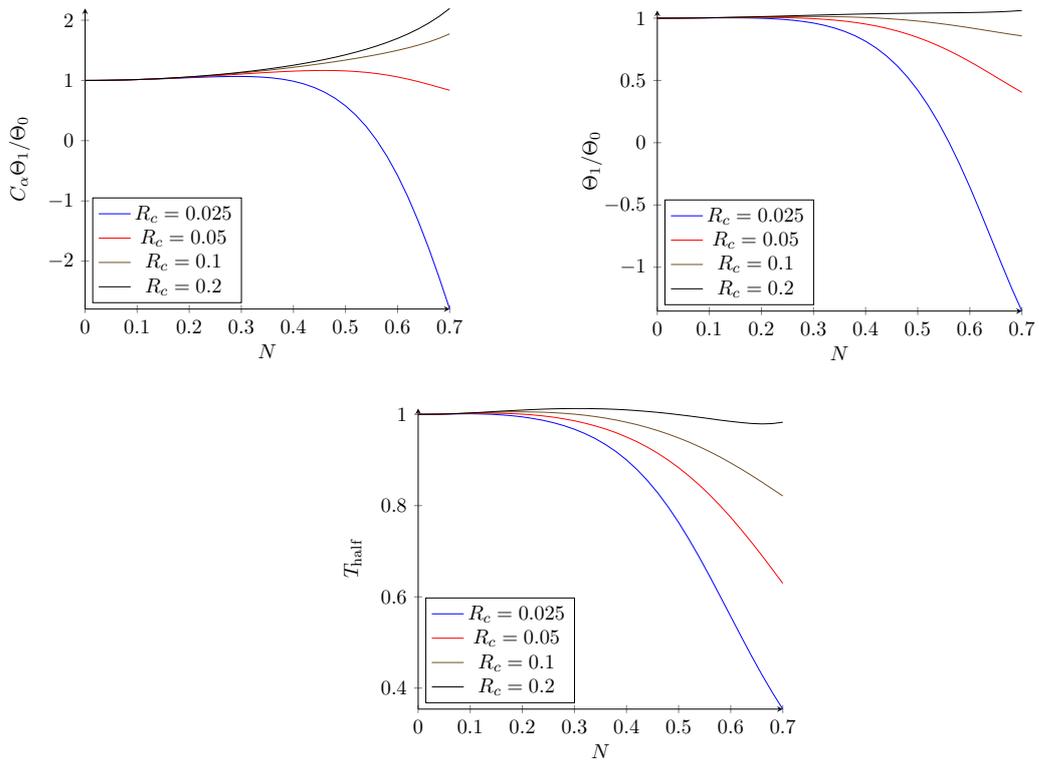
\begin{figure}
	\begin{center}
		\begin{minipage}[c]{0.45\linewidth}
			\begin{tikzpicture}[scale=0.7]
				\begin{axis}[ axis x line=bottom, axis y line = left,
					xlabel={$N$}, ylabel={$C_{\alpha}{\Theta_1}/{\Theta_0}$},
					legend entries={$R_c = 0.025$, $R_c = 0.05$, $R_c = 0.1$, $R_c = 0.2$}, legend style={at={(0.02,0.02)}, anchor=south west}
					]
					\addplot+[mark=none] table[x index=0, y index=1]{Data-Interp_Plot_N-E=5-delta=10-Rc=0-025-Results.txt};
					\addplot+[mark=none] table[x index =0, y index=1]{Data-Interp_Plot_N-E=5-delta=10-Rc=0-05-Results.txt};
					\addplot+[mark=none] table[x index =0, y index=1]{Data-Interp_Plot_N-E=5-delta=10-Rc=0-1-Results.txt};
					\addplot+[mark=none] table[x index =0, y index=1]{Data-Interp_Plot_N-E=5-delta=10-Rc=0-2-Results.txt};				
				\end{axis}
			\end{tikzpicture}
		\end{minipage} 
		\begin{minipage}[c]{0.45\linewidth}
			\begin{tikzpicture}[scale=0.7]
				\begin{axis}[ axis x line=bottom, axis y line = left,
					xlabel={$N$}, ylabel={${\Theta_1}/{\Theta_0}$},
					legend entries={$R_c = 0.025$, $R_c = 0.05$, $R_c = 0.1$, $R_c = 0.2$}, legend style={at={(0.02,0.02)}, anchor=south west}
					]
					\addplot+[mark=none] table[x index=0, y index=6]{Data-Interp_Plot_N-E=5-delta=10-Rc=0-025-Results.txt};
					\addplot+[mark=none] table[x index =0, y index=6]{Data-Interp_Plot_N-E=5-delta=10-Rc=0-05-Results.txt};
					\addplot+[mark=none] table[x index =0, y index=6]{Data-Interp_Plot_N-E=5-delta=10-Rc=0-1-Results.txt};
					\addplot+[mark=none] table[x index =0, y index=6]{Data-Interp_Plot_N-E=5-delta=10-Rc=0-2-Results.txt};				
				\end{axis}
			\end{tikzpicture}
		\end{minipage}
	\end{center}
	\begin{center}
			\begin{tikzpicture}[scale=0.7]
	\begin{axis}[ axis x line=bottom, axis y line = left,
		xlabel={$N$}, ylabel={$\Thalf$},
		legend entries={$R_c = 0.025$, $R_c = 0.05$, $R_c = 0.1$, $R_c = 0.2$}, legend style={at={(0.02,0.02)}, anchor=south west}
		]
		\addplot+[mark=none] table[x index=0, y index=1]{Data-Plot_N-E=5-nub=0-1-delta_custom-delta=10-Res-global.txt};
		\addplot+[mark=none] table[x index=0, y index=2]{Data-Plot_N-E=5-nub=0-1-delta_custom-delta=10-Res-global.txt};
		\addplot+[mark=none] table[x index=0, y index=3]{Data-Plot_N-E=5-nub=0-1-delta_custom-delta=10-Res-global.txt};
		\addplot+[mark=none] table[x index=0, y index=4]{Data-Plot_N-E=5-nub=0-1-delta_custom-delta=10-Res-global.txt};
		
	\end{axis}
\end{tikzpicture}
\end{center}

	\caption{$C_{\alpha}{\Theta_1}/{\Theta_0}$ (top, left) and ${\Theta_1}/{\Theta_0}$ (top, right) plotted against $N$, normalized by the value at $N=0$, for $\bar{\nu}_b = 0.1$, $\delta=10$, $E=5$ and different values of $R_c\in \{0.025, 0.05, 0.1, 0.2\}$. Bottom: normalized value of $\Thalf$ plotted against $N$ for the same sets of parameters.}
	\label{Fig:Interp-N-global-delta=10}
\end{figure}

As a result, the behaviour of the system appears to be driven by the factor $C_\alpha E\frac{\Theta_1}{\Theta_0}$, in the sense that increasing this quantity should result in reducing the velocity of the upper plate subject to a given load, thus increasing $\Thalf$ and enhancing the performance of the bearing. However, as shown in Figs.~\ref{Fig:Interp-N-global-delta=1} and~\ref{Fig:Interp-N-global-delta=10}, the behaviour of $\Thalf$ as a function of $N$ and $R_c$ seems more appropriately described by the ratio $\frac{\Theta_1}{\Theta_0}$ than by the quantity $C_{\alpha}E\frac{\Theta_1}{\Theta_0}$ itself. This suggests that computing the quotient $\frac{\Theta_1}{\Theta_0}$ using the explicit formulas provided by Corollary~\ref{cor_alpha_neq_1} may be a valuable tool to predict and compare the relative performance of squeeze-film bearings, depending on the physical parameters associated with the system.


\subsubsection{Conclusion}

We may conclude from these simulations that, at least for a moderate range of $N$ values and for small values of slip length $\delta$, the model predicts an enhancement of the performance of the bearing lubricated by a micropolar fluid, with respect to a bearing lubricated by a Newtonian fluid. This feature is generally affected by the value of parameter $R_c$, in the sense that large values of $R_c$ favor the enhancement of the performance of the squeeze-film bearing. Also, the introduction of a rough geometry for the lower plate, characterized by roughness parameter $E$, usually results in enlarging the range of $N$ values for which lubrication with a micropolar fluid gets more efficient than lubrication with a Newtonian fluid.

\section*{Acknowledgements}

The first author has been supported by a public grant as part of the
Investissement d'avenir project, reference ANR-11-LABX-0056-LMH,
LabEx LMH.

The second author has been supported by the {\it Croatian Science Foundation} under the project Multiscale problems in fluid mechanics - MultiFM (IP-2019-04-1140).

The third author has been partially supported by the {\it Ministerio de Econom\'ia y Competitividad} (Spain), under the project Proyecto Excelencia MTM2014-53309-P.

The authors would also like to thank the referees for their helpful remarks and suggestions.


\section*{Appendix}

This section completes Section~\ref{Sec:generalmicropolar} by providing the explicit formulas for the solution $(\tilde u_1,\tilde w_2)$ to system~\eqref{limit_system_1_reynolds}, in the critical case (Lemmas~\ref{lemma_alpha_neq_1} and~\ref{lemma_alpha_igual_1}) and in the super-critical case (Lemma~\ref{lm_super_critical}). In the critical case, we also give the first terms of the developments of $\tilde u_1,\tilde w_2, p$ and the coefficient $\Theta_{\lambda}$ appearing in Reynolds equation~\eqref{reynolds} in powers of $\lambda^2$ (see Corollary~\ref{cor_alpha_neq_1}). We conclude the Appendix by the proofs of these statements.


\begin{lm}\label{lemma_alpha_neq_1}
	For $\alpha\neq 1$, the solutions of system (\ref{limit_system_1_reynolds}) with boundary conditions (\ref{bc_top_system_1_reynolds})--(\ref{bc_botom_case2_system_1_reynolds})  are:
	\begin{equation}\label{tilde_u_1}
	\begin{array}{rl}
	\tilde u_1(y_1,y_3)=&\left(\left[{2N^2\over k}\Big(sh(ky_3)-\eta_\lambda sh(kh)\Big)+\gamma_\alpha(y_3-\eta_\lambda h)-(1-\eta_\lambda)\Big(\gamma_\alpha+{2N^2\over kh}sh(kh)\Big)\right]A
	\right.\\
	\noame
	&+\left[{2N^2\over k}\Big(ch(ky_3)-\eta_\lambda ch(kh)\Big)+(1-\eta_\lambda){2N^2\over k}\Big(-1+{y_3\over h}(1-ch(kh))\Big)\right]B \\
	\noame &
	\left.
	+{1\over 2(1-N^2)}\Big[y_3^2-h^2+(1-\eta_\lambda)(y_3h+h^2)\Big]\right)\partial_{y_1} p(y_1)\,,
	\end{array}
	\end{equation}
	$$\begin{array}{rl}\widetilde w_2(y_1,y_3)=&\left(\left[ch(ky_3)+{\gamma_\alpha\over 2}-(1-\eta_\lambda)\Big({\gamma_\alpha\over 2}+{N^2\over kh}sh(kh)\Big)\right]A\right.\\
	\noame &  \left.
	+\left[sh(ky_3)+(1-\eta_\lambda)(1-ch(kh)){N^2\over kh}\right]B
	+{1\over 2(1-N^2)}\Big[y_3+(1-\eta_\lambda){h\over 2}\Big]\right)\partial_{y_1}  p(y_1)\,,
	\end{array}
	$$
	where 
	\begin{equation}\label{gamma_eta}k=2N\sqrt{{1-N^2\over R_c}},\quad {\gamma_\alpha\over 2}={1-\alpha N^2\over \alpha-1},\quad \eta_\lambda=\left(1+{\alpha h\over \alpha-1}E_\lambda\right)^{-1}\,,
	\end{equation}
	$$\ba{rl}
	A= &\!\!\!\!{L\over 2(1-N^2) }\left({h\over 2}(1+\eta_\lambda)\left[4N^4\eta_\lambda(1-ch(kh))+{R_c\over\beta}k^2\right]\right.
	\\
	\noame &\!\!\!\!\! \left.-
	\left[k sh(kh)+(1-\eta_\lambda)(1-ch(kh)){N^2\over h}\right]\left[{R_c\over\beta}-2N^2h^2\eta_\lambda\right]\right)\,,
	\\
	\\
	B =& \!\!\!\! {L\over 2(1-N^2) }\Big(N^2h(1+\eta_\lambda)\eta_\lambda\left[\gamma_\alpha k h+2N^2sh(kh)\right]
	\\
	\noame&\!\!\!\!  +
	k\left[{R_c\over\beta}-2N^2h^2\eta_\lambda\right]\left[ch(kh)+{\gamma_\alpha\over 2}-(1-\eta_\lambda)\left({\gamma_\alpha\over 2}+{N^2\over k h}sh(kh)\right)\right]\Big)\,, \\
	\\
	\dis L= 
	&\!\!\!\!-\left(\left[{\gamma_\alpha\over 2}+ch(kh)-(1-\eta_\lambda)\left({\gamma_\alpha\over 2}+{N^2\over kh }sh(kh)\right)\right]\left[4N^4\eta_\lambda(1-ch(kh))+{R_c\over \beta}k^2\right]\right.\\
	\noame & \!\!\!\! \left.
	+2N^2\eta_\lambda\left[sh(kh)+(1-\eta_\lambda)(1-ch(kh)){N^2\over k h}\right]
	\left[\gamma_\alpha k h +2N^2 sh(kh)\right]\right)^{-1}\,.
	\ea$$
\end{lm}

\begin{lm}\label{lemma_alpha_igual_1}
	For $\alpha= 1$, the solutions of system (\ref{limit_system_1_reynolds}) with boundary conditions (\ref{bc_top_system_1_reynolds})--(\ref{bc_botom_case2_system_1_reynolds})  are
	\begin{equation}\label{tilde_u_1_alpha_igual_1}
	\begin{array}{rl}
	\tilde u_1(y_1,y_3)=&\!\!\!\!
	\left(\Big[{2N^2\over k}\big(ch(ky_3)-\mu_\lambda ch(kh)\big) -(1-\mu_\lambda){2N^2\over k}\Big(1-(1-ch(kh)){sh(ky_3)\over sh(kh)}\Big)\Big]B'\right.
	\\
	\noame & 
	+{1\over 2(1-N^2)}\big(y_3^2-\mu_\lambda h^2\big)-(1-\mu_\lambda){h^2\over 1-N^2}{sh(ky_3)\over sh(kh)}
	\\
	\noame & 
	\left. +\Big[{1\over 1-N^2}(y_3-\mu_\lambda h)-(1-\mu_\lambda){h\over 1-N^2}{sh(ky_3)\over sh(kh)}\Big]A'\right)\partial_{y_1}  p(y_1)\,,\\
	\noame
	\tilde w_2(y_1,y_3)=&\!\!\!\!\left(\Big[sh(ky_3)+(1-\mu_\lambda)(1-ch(kh)){ch(ky_3)\over sh(kh)}\Big]B'\right.\\
	\noame &
	+{y_3\over 2(1-N^2)}-(1-\mu_\lambda){kh^2\over 2N^2(1-N^2)}{ch(ky_3)\over sh(kh)}\\
	\noame & 
	\left.
	+\Big[{1\over 2(1-N^2)}-(1-\mu_\lambda){kh\over 2N^2(1-N^2)}{ch(ky_3)\over sh(kh)}\Big]A'
	\right)\partial_{y_1} p(y_1)\,,
	\end{array}
	\end{equation}
	where
	\begin{equation}\label{eta_N}
	\mu_\lambda=\left(1-{N^2\over 1-N^2}{sh(kh)\over k}E_\lambda\right)^{-1}\,,
	\end{equation}
	$$\ba{rl}
	A'=&\!\!\!\!L'\left(\Big[4N^4\mu_\lambda(1-ch(kh))+{R_c\over \beta}k^2\Big]\Big[h-(1-\mu_\lambda)coth(kh){h^2k\over 2N^2}\Big]\right.\\
	\noame &\!\!\!\!\left. -k\Big[sh(kh)+(1-\mu_\lambda)coth(kh)(1-ch(kh))\Big]\Big[{R_c\over \beta}-2N^2h^2\mu_\lambda\Big]\right)\,
	\\
	\\
	B' =&\!\!\!\!\!k{L' \over 2(1-N^2)}\left(2N^2h^2\mu_\lambda+{R_c\over \beta}-(1-\mu_\lambda)coth(kh){hk\over N^2}{R_c\over\beta}\right)\,,
	\\
	\\
	L' =&\!\!\!\!-\left(\Big[1-(1-\mu_\lambda)coth(kh){hk\over N^2}\Big]\Big[4N^4\mu_\lambda(1-ch(kh))+{R_c\over \beta}k^2\Big]
	\right.\\
	\noame & \!\!\!\!\left. +4N^2hk\mu_\lambda\Big[sh(kh)+(1-\mu_\lambda)coth(kh)(1-ch(kh))\Big]\right)^{-1}.
	\ea$$
\end{lm}


Let us define
$$C_\alpha={\alpha h\over \alpha-1},\quad C_N={N^2\over 1-N^2}{sh(kh)\over k} ,\quad E=\int_{\widehat Q}|D_z \widehat\phi^{1,1}|^2\, dz.$$
as consequence of the previous results, we get the following developments in powers of $\lambda$, which will be useful for the numerical part.
\begin{cor}\label{cor_alpha_neq_1}
	We obtain the following developments: 
	\begin{equation}\label{develop_u1_w2_alpha_neq_1}\begin{array}{rl}
	\tilde u_1(y_1,y_3)=&v_0(y_3)\partial_{y_1}p_0(y_1)+C_j E\lambda^2\Big(v_1(y_3)\partial_{y_1}p_0(y_1)+v_0(y_3)\partial_{y_1}p_1(y_1)\Big)   +O(\lambda^4),\\
	\tilde w_2(y_1,y_3)= & \varpi_0(y_3)\partial_{y_1}p_0(y_1)+ C_j E\lambda^2\Big(\varpi_1(y_3)\partial_{y_1}p_0(y_1)+\varpi_0(y_3)\partial_{y_1}p_1(y_1)\Big) +O(\lambda^4),\\
	p(y_1)=&p_0(y_1)+C_j E\lambda^2 p_1(y_1)+O(\lambda^4),\\
	\Theta_\lambda=&\Theta_0-C_jE\lambda^2 \Theta_1+O(\lambda^4),
	\end{array}
	\end{equation}
	where $j=\alpha$ if $\alpha\neq 1$ and $j=N$ if $\alpha=1$. Here, $p_0$ satisfies the following equation
	\begin{equation}\label{reynoldsp0_neq1}\int_0^1\Theta_0\partial_{y_1}p_0(y_1)\partial_{y_1}\theta(y_1)\,dy_1=\int_0^1 S\theta(y_1)\,dy_1,\quad\forall \theta\in H^1(0,1),
	\end{equation}
	and $p_1$ the following one
	\begin{equation}\label{reynoldsp1_neq1}\int_0^1\Theta_0\partial_{y_1}p_1(y_1)\partial_{y_1}\theta(y_1)\,dy_1=\int_0^1 \Theta_1\partial_{y_1}p_0(y_1)\partial_{y_1}\theta(y_1)\,dy_1,\quad\forall \theta\in H^1(0,1).
	\end{equation}
	In particular, $p_1$ is given by $p_1(y_1)=\frac{\Theta_1}{\Theta_0}p_0(y_1)$. \\

	\noindent For $\alpha\neq 1$, $v_i$, $\varpi_i$, $\Theta_i$, $i=0,1$, are defined as follows  
	$$\begin{array}{rl}
	v_0(y_3)=&{1\over 2(1-N^2)}\Big(-\Big[{2N^2\over k}\left(sh(ky_3)-sh(kh)\right)+\gamma_\alpha(y_3-h)\Big]L_0A_0\\
	& \qquad \qquad - {2N^2\over k}\left(ch(ky_3)-ch(kh)\right)L_0B_0+y_3^2-h^2\Big),
	\\
	v_1(y_3)=& {1\over 2(1-N^2)}\Big(-\Big[{2N^2\over kh}sh(kh)+\gamma_\alpha\Big](h-1)L_0A_0-{2N^2\over kh}(1-ch(kh))\Big({y_3}-h\Big)L_0B_0
	\\
	&\qquad \qquad -\Big[{2N^2\over k}\left(sh(ky_3)-sh(kh)\right)+\gamma_\alpha(y_3-h)\Big]L_0(A_1+L_0L_1A_0)\\
	& \qquad \qquad -{2N^2\over k}\Big[ ch(ky_3)-ch(kh)\Big]L_0(B_1+L_0L_1B_0)+ y_3h+h^2 \Big),
	\\
	\varpi_0(y_3)=& {1\over 2(1-N^2)}\Big[-\Big(ch(ky_3)+\gamma_\alpha\Big)L_0A_0-sh(ky_3)L_0B_0+y_3\Big],\\
	\\
	\varpi_1(y_3)=& {1\over 2(1-N^2)}\Big[-\Big(ch(ky_3)+\gamma_\alpha\Big)L_0(A_1+L_0L_1A_0)+\Big({\gamma_\alpha\over 2}+{N^2\over kh}sh(kh)\Big)L_0A_0\\
	&\qquad\qquad
	-sh(ky_3)L_0\left(B_1+L_0L_1B_0\right)- (1-ch(kh)){N^2\over kh}L_0B_0+{h\over 2}\Big]
	\\
	\Theta_0=&{h^3\over 3(1-N^2)}\\
	& -{1\over 2(1-N^2)}\Big[\Big({2N^2\over k}\left({ch(kh)-1\over k}-hsh(kh)\right)-{\gamma_\alpha \over 2}h^2\Big)L_0A_0+{2N^2\over k}\left({sh(kh)\over k}-hch(kh)\right)L_0B_0\Big],\\
	\\
	\Theta_1=& {3h^3\over 4(1-N^2)} \\
	& -{1\over 2(1-N^2)}\Big[\Big({2N^2\over k}\left({ch(kh)-1\over k}-hsh(kh)\right)-{\gamma_\alpha \over 2}h^2\Big)L_0(A_1+L_0L_1 A_0)\\
	&\qquad\qquad +\Big(-{2N^2\over k}sh(kh)(1-h)-\gamma_\alpha h(1-h)\Big)L_0A_0\\
	&\qquad\qquad +{2N^2\over k}\left({sh(kh)\over k}-hch(kh)\right)L_0\left(B_1+L_0L_1B_0\right)\\
	&\qquad\qquad +\Big({2N^2\over k}hch(kh)-(1+ch(kh))h{N^2\over k}\Big)L_0B_0\Big],
	\end{array}$$
	where $A_i, B_i, L_i$, $i=0,1$ are given by
	\begin{equation}\label{A0L0B0}
	\begin{array}{l}
	A_0=h\left[4N^4(1-ch(kh))+{R_c\over \beta}k^2\right]-ksh(kh)\left[{R_c\over \beta}-2N^2h^2\right],
	\\ 
	A_1=-h\left[4N^4(1-ch(kh))+{R_c\over \beta}{k^2\over 2}\right]-2N^2h^2k sh(kh)-(1-ch(kh)){N^2\over h}{R_c\over \beta},
	\\
	B_0= 2N^2h\left[\gamma_\alpha kh+2N^2sh(kh)\right]+k\left[{R_c\over \beta}-2N^2h^2\right]\left[ch(kh)+{\gamma_\alpha\over 2}\right],\\
	B_1= -{\gamma_\alpha\over 2}k\left[2N^2h^2+{R_c\over \beta}\right]-sh(kh)\left[4N^4h+{R_c\over \beta}{N^2\over h}\right]+2N^2h^2kch(kh),
	\\
	L_0 = \left(\left[{\gamma_\alpha\over 2}+ch(kh)\right]\left[4N^4(1-ch(kh))+{R_c\over \beta}k^2\right]+2N^2\,sh(kh)\left[\gamma_\alpha kh+2N^2\,sh(kh)\right]\right)^{-1}\!\!\!,\\
	L_1 =   4N^4(1-ch(kh))\left[\gamma_\alpha+ch(kh)+{N^2\over kh}sh(kh)\right]+ {R_c\over \beta}k^2\left[{\gamma_\alpha\over 2}+{N^2\over kh}sh(kh)\right] \\
	\quad\quad +2N^2\left[sh(kh)-{N^2\over kh}(1-ch(kh))\right]\left[\gamma_\alpha kh+2N^2 sh(kh)\right].
	\end{array}
	\end{equation}
	For $\alpha= 1$, $v_i$, $\varpi_i$, $\Theta_i$, $i=0,1$, are defined as follows  
	$$\begin{array}{rl}
	v_0(y_3)=& -{N^2\over (1-N^2)}\big(ch(ky_3)-ch(kh)\big)L'_0B'_0+{1\over 2(1-N^2)}\big(y_3^2-  h^2\big)-{1\over 1-N^2}(y_3-h)L'_0A'_0,\\
	\\
	v_1(y_3)=& -{N^2\over 1-N^2}\big(ch(ky_3)-ch(kh)\big)L'_0(B'_1-L'_0L'_1B'_0)-{N^2\over 1-N^2}(1-ch(kh)){sh(ky_3)\over sh(kh)}L'_0B'_0
	\\
	& 
	+{h^2\over 1-N^2}\Big(-{1\over 2}+{sh(ky_3)\over sh(kh)}\Big)   -{1\over 1-N^2}(y_3-h)L'_0(A'_1-L'_0L'_1A'_0)-L'_0A'_0\Big(-1+{sh(ky_3)\over sh(kh)}\Big),
	\\
	\varpi_0(y_3)=& -{k\over 2(1-N^2)}sh(ky_3)L'_0B'_0+{y_3\over 2(1-N^2)}-{1\over 2(1-N^2)}L'_0A'_0,\\
	\\
	\varpi_1(y_3)=&-{k\over 2(1-N^2)}sh(ky_3)L'_0(B'_1-L'_0L'_1B_0)+ {k\over 2(1-N^2)}(1-ch(kh)){ch(ky_3)\over sh(kh)}L'_0B'_0\\
	& {kh^2\over 2N^2(1-N^2)}{ch(ky_3)\over sh(kh)}-{1\over 2(1-N^2)}L'_0(A'_1-L'_0L'_1A'_0)-{kh\over 2N^2(1-N^2)}{ch(ky_3)\over sh(kh)}L'_0A'_0,\\
	\\
	\Theta_0=& {h^3\over 3(1-N^2)}-{h^2\over 2(1-N^2)}L'_0A'_0+{N^2\over 1-N^2}\Big({sh(kh)\over k}- h\,  ch(kh)\Big)L'_0B'_0,
	\\
	\Theta_1=&-{h^2\over 1-N^2}\Big({h\over 2}+{1-ch(kh)\over k\,sh(kh)}\Big)
	+{h^2\over 2(1-N^2)}L'_0\Big(A'_1-L'_0L'_1A'_0\Big)+{h\over 1-N^2}\Big(h+{1-ch(kh) \over k\,sh(kh)}\Big)L'_0A'_0\\
	&-{N^2\over 1-N^2}\Big({sh(kh)\over k}- h\,  ch(kh)\Big)L'_0\Big(B'_1-L'_0L'_1B'_0\Big)\\
	& -{N^2\over 1-N^2}\Big(- h\, ch(kh)+ \Big[h+{(1-ch(kh))^2\over k\, sh(kh)}\Big]\Big)L'_0B'_0.
	\end{array}$$
	where $A'_i, B'_i, L'_i$, $i=0,1$ are given by
	$$\begin{array}{rl}
	A'_0=& h\Big[4N^4(1-ch(kh))+{R_c\over \beta}k^2\Big]-k\, sh(kh)\Big[{R_c\over \beta}-2N^2h^2\Big],\\
	A'_1=& h\Big[4N^4(1-ch(kh))+2N^2hksh(kh)\Big]+{R_c\over \beta}k\, coth(kh)\Big[1-ch(kh)+{k^2h^2\over 2N^2}\Big],\\
	B'_0=&2N^2h^2+{R_c\over \beta},\\
	B'_1=& 2N^2h^2+coth(kh){hk\over N^2}{R_c\over\beta},\\
	L'_0=&\Big(4N^4(1-ch(kh))+{R_c\over \beta}k^2 + 4N^2hk sh(kh)\Big)^{-1},\\
	L'_1=&4N^4(1-ch(kh)) + coth(kh) {k^3h\over N^2}{R_c\over \beta}+4khN^2sh(kh).
	\end{array}$$
	
\end{cor}


\begin{lm}\label{lm_super_critical} In the super-critical case $1<\delta<{3\over 2}\ell-{1\over 2}$, the solutions of system (\ref{limit_system_1_reynolds}) with boundary conditions
	$$
	\begin{array}{l}
	\widetilde u_1=\widetilde w_2=0\quad\hbox{on }\Gamma^1,\quad \widetilde u_1= \partial_{y_3}\widetilde w_2=0\quad\hbox{on }\Gamma\,,
	\end{array}
	$$
	are
	\begin{equation}\label{lm_super_critical_u1}
	\begin{array}{rl}
	\tilde u_1(y_1,y_3)=& \Big({2N^2\over k}\left[ sh(ky_3)-{y_3\over h}sh(kh)\right]A''\\
	\noame &
	+ {2N^2\over k}\left[ ch(ky_3)-{y_3\over h}(ch(kh)-1)-1 \right]B'' +{1\over 2(1-N^2)}\left[y_3^2-y_3h\right]\Big)\partial_{y_1}  p(y_1)\,,\\
	\\
	\tilde w_2(y_1,y_3)=&  \Big(\left[ch(ky_3)-{N^2\over kh}sh(kh)\right]A''\\
	\noame & 
	+\left[ sh(ky_3)-{N^2\over kh}(ch(kh)-1)\right]B''+{1\over 2(1-N^2)}\left[ y_3-{h\over 2}\right]\Big)\partial_{y_1}  p(y_1)\,,
	\end{array}
	\end{equation}
	where
	$$
	\begin{array}{rl}
	A''=&\!\!\! {1\over 2(1-N^2)}L''\left[-{h\over 2}+{sh(kh)\over k}-{N^2\over k^2h}(ch(kh)-1)\right]\,,
	\\
	\\
	B''=&\!\!\! -{1\over 2(1-N^2)}\Big[{ch(kh)\over k}-{N^2\over k^2h}sh(kh)\Big]\,,
	\\
	\\
	L''=&\left[ch(kh)-{N^2\over kh}sh(kh)\right]^{-1}\,.
	\end{array}
	$$
\end{lm}

In the rest of the appendix, we give the proofs of Lemmas \ref{lemma_alpha_neq_1}, \ref{lemma_alpha_igual_1}, \ref{lm_super_critical} and Corollary \ref{cor_alpha_neq_1}.\\

\noindent {\bf Proof of Lemma \ref{lemma_alpha_neq_1}. }  Let us start with system (\ref{limit_system_1_reynolds})--(\ref{bc_botom_case2_system_1_reynolds}). Integrating (\ref{limit_system_1_reynolds})$_1$ in $y_3$, we obtain 
\begin{equation}\label{partialu1}
\partial_{y_3}\widetilde u_1(y_1,y_3)=\partial_{y_1}  p(y_1)\, y_3 + 2N^2\widetilde w_2(y_1,y_3) + K_1(y_1),
\end{equation}
where $K_1$ is an unknown function. Putting (\ref{partialu1}) into (\ref{limit_system_1_reynolds})$_2$, we obtain
\begin{equation}\label{partial2w2}
\partial_{y_3}^2 \widetilde w_2(y_1,y_3)-{4N^2\over R_c}(1-N^2)\widetilde w_2(y_1,y_3)=-{2N^2\over R_c}\partial_{y_1}  p(y_1)\,y_3-{2N^2\over R_c} K_1(y_1),
\end{equation}
whose solution can be written as
\begin{equation}\label{w2}
\widetilde w_2(y_1,y_3)=c_1(y_1)e^{ky_3}+c_2(y_1)e^{-k y_3}+{y_3\over 2(1-N^2)}\partial_{y_1}  p(y_1) + {K_1(y_1)\over 2(1-N^2)},
\end{equation}
where $c_i(y_1)$, $i=1,2$ are unknown functions and $\displaystyle k=2N\sqrt{{1-N^2\over R_c}}$.

Putting (\ref{w2}) in (\ref{partialu1}), we obtain
\begin{equation}\label{partialu1-2}
\partial_{y_3}\widetilde u_1(y_1,y_3)={1\over 1-N^2}\partial_{y_1}  p(y_1)\,y_3 + 2N^2\left(c_1(y_1)e^{ky_3}+c_2(y_1)e^{-ky_3}\right)+{K_1(y_1)\over 1-N^2}.
\end{equation}

Integrating (\ref{partialu1-2}), we get
\begin{equation}\label{u_1}
\begin{array}{rl}
\dis \widetilde u_1(y_1,y_3)=&\dis {2N^2\over k}\left(\tilde A(y_1) sh(ky_3)+\tilde B(y_1)ch(ky_3)\right)+{y_3^2\over 2(1-N^2)}\partial_{y_1}  p(y_1)\\
\noame &\dis +{K_1(y_1)\over 1-N^2}y_3+K_2(y_1),
\end{array}
\end{equation}
where 
$$\tilde A(y_1)=c_1(y_1)+c_2(y_1),\quad \tilde  B(y_1)=c_1(y_1)-c_2(y_1),$$
and $K_2(y_1)$ is an unknown function. $\widetilde w_2$ in (\ref{w2}) can also be written as follows:
\begin{equation}\label{w2-2}
\widetilde w_2(y_1,y_3)=\tilde A(y_1)ch(ky_3)+\tilde B(y_1)sh(ky_3)+{y_3\over 2(1-N^2)}\partial_{y_1}  p(y_1) + {K_1(y_1)\over 2(1-N^2)}.
\end{equation}

Using (\ref{partialu1-2}) and (\ref{w2-2}) in the boundary condition (\ref{bc_botom_case2_system_1_reynolds})$_1$, we obtain
\begin{equation}\label{systemK1K2}
K_1(y_1)-E_\lambda{\alpha(1-N^2)\over \alpha-1} K_2(y_1)=\gamma_\alpha (1-N^2) \tilde A(y_1) + E_\lambda{\alpha(1-N^2)\over \alpha-1}{2N^2\over k}\tilde B(y_1)\,
\end{equation}
with $\gamma_\alpha$ defined in (\ref{gamma_eta}). Using the condition $\widetilde u_1=0$ on $\Gamma_1$ in (\ref{u_1}),  (\ref{systemK1K2}) and taking into account  $\eta_\lambda$ defined in (\ref{gamma_eta}), we obtain
$$
\begin{array}{l}
\dis K_2(y_1)= -\eta_\lambda\left({2N^2\over k}sh(kh)+h\gamma_\alpha\right) \tilde A(y_1)\\
\noame \dis -\eta_\lambda \left({2N^2\over k}ch(kh)+E_\lambda{\alpha h\over \alpha-1}{2N^2\over k} \right)\tilde  B(y_1) -\eta_\lambda{h^2\over 2(1-N^2)}\partial_{y_1} p(y_1).
\end{array}
$$
Plugging this expression in (\ref{systemK1K2}) yields
$$
\begin{array}{l}
\dis K_1(y_1)=\dis (1-N^2)\left[\gamma_\alpha-E_\lambda{\alpha\over \alpha-1}\eta_\lambda\left({2N^2\over k}sh(kh)+h\gamma_\alpha\right)\right]\tilde A(y_1)\\
\noame \dis +E_\lambda{\alpha(1-N^2)\over \alpha-1}{2N^2\over k}\left[1-\eta_\lambda\left(ch(kh)+E_\lambda{\alpha h\over \alpha-1}\right)\right]\tilde B(y_1)-E_\lambda{\alpha h^2\over 2(\alpha-1)}\eta_\lambda\partial_{y_1}  p(y_1).
\end{array}
$$
Taking into account that 
$$\eta_\lambda=\left(1+{\alpha h\over \alpha-1}E_\lambda\right)^{-1}\quad\Leftrightarrow\quad \eta_\lambda^{-1}-1={\alpha h\over \alpha-1} E_\lambda
\quad\Leftrightarrow\quad 1-\eta_\lambda={\alpha h\over \alpha-1}E_\lambda\eta_\lambda,$$
we rewrite $K_1$ and $K_2$ as follows
\beq\label{K1}
\begin{array}{rl}
\dis K_1(y_1)=&\dis (1-N^2)\left[\gamma_\alpha-{1-\eta_\lambda\over h}\left({2N^2\over k}sh(kh)+h\gamma_\alpha\right)\right]\tilde A(y_1)\\
\noame &\dis +{2N^2\over k h}(1-N^2)(1-ch(kh))(1-\eta_\lambda)\tilde B(y_1)-(1-\eta_\lambda){h\over 2}\partial_{y_1}  p(y_1).
\end{array}
\end{equation}
\beq\label{K2}
\begin{array}{rl}
\dis K_2(y_1)=&\displaystyle  -\eta_\lambda\left({2N^2\over k}sh(kh)+h\gamma_\alpha\right) \tilde A(y_1) - {2N^2\over k}\left(\eta_\lambda(ch(kh)-1)+1\right)\tilde B(y_1)\\
\noame &\dis -\eta_\lambda{h^2\over 2(1-N^2)}\partial_{y_1} p(y_1).
\end{array}
\end{equation}

From condition $\widetilde w_2=0$ on $\Gamma^1$ and (\ref{bc_botom_case2_system_1_reynolds})$_2$, we obtain, using (\ref{K1}) and (\ref{K2}), the following system  
$$Q\left(
\ba{c}\dis \tilde A\\
\noame\dis \tilde B\ea
\right)
=\left(
\ba{c}
\dis  -{h\over 4(1-N^2)}(1+\eta_\lambda)\\
\noame\dis
{{R_c\over \beta}-2N^2h^2\eta_\lambda\over  2(1-N^2)}
\ea
\right)\partial_{y_1}  p(y_1)$$
where $Q$  is the matrix defined by
\begin{equation}Q=\left(\ba{ll}
{\gamma_\alpha\over 2}+ch(kh) -{1-\eta_\lambda\over 2 h}\left(h\gamma_\alpha+{2N^2\over k}sh(kh)\right)& sh(kh)+ {1-\eta_\lambda\over h}(1-ch(kh)){N^2\over k}\\
\noame\dis
2N^2\eta_\lambda\left(\gamma_\alpha h+{2N^2\over k}sh(kh)\right)& -{4N^4\over k}\eta_\lambda(1-ch(kh))-{R_c\over \beta}k
\end{array}\right)\,.
\end{equation}
The, the solution of this system is given by 
$$\tilde A(y_1)=A\partial_{y_1}  p(y_1),\quad \tilde B(y_1)=B\partial_{y_1}  p(y_1),$$
where $A$ and  $B$ are solution of 
$$\ba{l}\dis
Q\left(
\ba{c}\dis A\\
\noame\dis B\ea
\right)
={1\over 2(1-N^2)}\left(
\ba{c}
\dis  -{h\over 2}(1+\eta_\lambda)\\
\noame\dis
{{R_c\over \beta}-2N^2h^2\eta_\lambda}
\ea
\right)\,,
\ea$$
Computing $A$, $B$, then $\widetilde u_1$ and $\widetilde w_2$ are obtained by (\ref{u_1}) and (\ref{w2-2}) as functions of $p$ and of known data. 
\par\hfill$\square$
\\

\noindent {\bf Proof of Lemma \ref{lemma_alpha_igual_1}. }  The beginning of the proof is as in Lemma \ref{lemma_alpha_neq_1}. Using (\ref{partialu1-2}), (\ref{u_1}) and (\ref{w2-2}) in the boundary condition (\ref{bc_botom_case2_system_1_reynolds})$_1$, for $\alpha=1$, we obtain 
\begin{equation}\label{A_alpha_igual_1}
\tilde A(y_1)=-{E_\lambda\over 2(1-N^2)}K_2(y_1)-{N^2\over 1-N^2}{E_\lambda\over k}\tilde B(y_1)\,.
\end{equation}
 with $k$ given in (\ref{gamma_eta}). Using the condition $\widetilde u_1=0$ on $\Gamma_1$ in (\ref{u_1}), taking (\ref{A_alpha_igual_1}) into account, we obtain
\begin{equation}\label{K2_alpha_igual_1}
K_2(y_1)=-\mu_\lambda{h \over 1-N^2}K_1(y_1)+
{2N^2\over k}\Big(\mu_\lambda(1-ch(kh))-1\Big)\tilde B(y_1)
-\mu_\lambda{h^2 \over 2(1-N^2)}\partial_{y_1}  p(y_1)\,,
\end{equation}
with $\mu_\lambda$ defined in (\ref{eta_N}). Using the definition of $\mu_\lambda$, we rewrite $A$ as follows
\begin{equation}\label{A_alpha_igual_1_2}
\begin{array}{rl}
 \tilde A(y_1)=& \dis -(1-\mu_\lambda){1\over 2N^2(1-N^2)}{kh\over sh(kh)}K_1(y_1)+(1-\mu_\lambda){1-ch(kh)\over sh(kh)}\tilde B(y_1)\\
 \noame
 &\dis -(1-\mu_\lambda){1\over 2N^2(1-N^2)}{kh^2\over sh(kh)}\partial_{y_1}p(y_1).
 \end{array}
\end{equation}
From the conditions $\widetilde w_2=0$ on $\Gamma^1$ and (\ref{bc_botom_case2_system_1_reynolds})$_2$, using  (\ref{K2_alpha_igual_1}) and (\ref{A_alpha_igual_1_2}) the following system is obtained
$$Q'\left(
\ba{c}\dis K_1\\
\noame\dis \tilde B\ea
\right)
=\left(
\ba{c}
\dis {1\over 2(1-N^2)}\left(-h+(1-\mu_\lambda)coth(kh){h^2k\over 2N^2}\right)\\
\noame\dis
{{R_c\over\beta}-2N^2 h^2\mu_\lambda \over 2(1-N^2)}
\ea
\right)\partial_{y_1}  p(y_1)\,,$$
where $Q'$ is the matrix defined by 
$$Q'=\left(\begin{array}{cc}
{1\over 2(1-N^2)}\Big(1-(1-\mu_\lambda)coth(kh){k h\over N^2}\Big) & 
sh(kh)+(1-\mu_\lambda)coth(kh)(1-ch(kh))
\\
\noame
{2N^2 h\over 1-N^2}\mu_\lambda & {-4N^4\over k}\mu_\lambda(1-ch(kh))-{R_c\over \beta}k\end{array}\right)\,.$$
The solution of this system is given by
$$K_1(y_1)=A'\partial_{y_1} p(y_1),\quad B(y_1)=B' \partial_{y_1} p(y_1)\,,$$
where $A'$ and $B'$  are solution of 
$$\ba{l}\dis
Q'\left(
\ba{c}\dis A'\\
\noame\dis B'\ea
\right)
=\left(
\ba{c}
\dis {1\over 2(1-N^2)}\left(-h+(1-\mu_\lambda)coth(kh){h^2k\over 2N^2}\right)\\
\noame\dis
{{R_c\over\beta}-2N^2 h^2\mu_\lambda \over 2(1-N^2)}
\ea
\right)\,.
\ea$$
Computing $A'$, $B'$,  then $\widetilde u_1$ and $\widetilde w_2$ are obtained as functions of $p$ and of known data.
\par\hfill$\square$\\

\noindent {\bf Proof of Corollary  \ref{cor_alpha_neq_1}. }  We  first remark that the roughness parameter $E_\lambda$ given in Theorem \ref{thm_effective} satisfies 
\[
E_\lambda=\int_{\widehat Q}|D_z \widehat\phi^{1,\lambda}|\, dz = \lambda^2 \int_{\widehat Q}|D_z \widehat\phi^{1,1}|\, dz = \lambda^2 E. 
\]
We explain the case $\alpha\neq 1$ (for the case $\alpha=1$ proceed similarly). Using power series of $\lambda^2$ and omitting terms of order $O(\lambda^4)$, there holds
\begin{align*}
\eta_\lambda & = (1+C_\alpha E\lambda^2)^{-1} \sim 1- C_\alpha E\lambda^2,\quad \hbox{with}\quad C_\alpha={\alpha h\over \alpha-1},\\
\eta_\lambda^2 & =  (1+C_\alpha E\lambda^2)^{-2} \sim 1- 2C_\alpha E\lambda^2.
\end{align*}
Using the development of $\eta_\lambda$ in terms of $\lambda^2$ in $A, B$ and $L$ given in Lemma \ref{lemma_alpha_neq_1}, we get
$$\begin{array}{l}
A=-{1\over 2(1-N^2)}L_0\Big(A_0+C_\alpha E\lambda^2[A_1+L_0L_1 A_0]\Big),\\
B=-{1\over 2(1-N^2)}L_0\left(B_0+C_\alpha E\lambda^2\left[B_1+L_0L_1B_0\right] \right),\\
L=-L_0\left(1+L_0L_1 C_\alpha E\lambda^2\right),
\end{array}$$
with $A_i, B_i$ and $L_i$, $i=0,1$ given by (\ref{A0L0B0}). Next, developing the pressure as $p(y_1)=p_0(y_1)+C_\alpha E\lambda^2 p_1(y_1)+O(\lambda^4)$ and using previous development of $A, B$ and $L$ in the expressions of $\tilde u_1$ and $\tilde w_2$ given in (\ref{tilde_u_1}), we get the expressions (\ref{develop_u1_w2_alpha_neq_1}). 

Finally, using again the development of $A, B$ and $L$ given above,  the development of $p$ and the development of $\tilde u_1$ in the Reynolds equation (\ref{reynolds}), we deduce that $p_0$ and $p_1$ satisfy (\ref{reynoldsp0_neq1}) and (\ref{reynoldsp1_neq1}), respectively.  Combining the equations on $p_0$ and $p_1$, one easily sees that $\frac{\Theta_0}{\Theta_1}p_1$ and $p_0$ satisfy the same equation, hence they are equal by uniqueness of the solution to~(\ref{reynoldsp0_neq1}).
\hfill$\square$ 
\\

\noindent {\bf Proof of Lemma \ref{lm_super_critical}}. The beginning of the proof is as in Lemma \ref{lemma_alpha_neq_1}. In this case, we consider the boundary conditions given in (\ref{super_critical_bc}). Thus, using  (\ref{u_1}) and boundary conditions $\tilde u_1(y_1,0)=0$ and $\tilde u_1(y_1,h)=0$, respectively, we have
$$\begin{array}{l}
K_2(y_1)=-{2N^2\over k}\tilde B(y_1)\,,\\
\noame\dis
K_1(y_1)=-{2N^2(1-N^2)\over kh}sh(kh) \tilde A(y_1)-{2N^2(1-N^2)\over kh}(ch(kh)-1)\tilde B(y_1)-{h\over 2}\partial_{y_1}  p(y_1)\,,
\end{array}$$
with $k$ given in (\ref{gamma_eta}). From the boundary conditions $\tilde w_2(y_1,h)=0$ and $\partial_{y_3}\tilde w_2(y_1,0)=0$, 
the following system is obtained
$$Q''\left(
\ba{c}\dis \tilde A\\
\noame\dis \tilde B\ea
\right)
=\left(
\ba{c}
\dis {-h\over 4(1-N^2)}\\
\noame\dis
-{1 \over 2(1-N^2)k}
\ea
\right)\partial_{y_1}  p(y_1),$$
where $Q''$ is the matrix defined by 
$$Q''=\left(\begin{array}{cc}
ch(kh)-{N^2\over kh}sh(kh) & sh(kh)-{N^2\over kh}(ch(kh)-1)\\
\noame
0 & 1
\end{array}\right)\,.$$
The solution of this system is given by
$$\tilde A(y_1)=A''\partial_{y_1} p(y_1),\quad \tilde B(y_1)=B'' \partial_{y_1} p(y_1)\,,$$
where $A''$ and $B''$  are solution of 
$$Q''\left(
\ba{c}\dis  A''\\
\noame\dis B''\ea
\right)
=\left(
\ba{c}
\dis {-h\over 4(1-N^2)}\\
\noame\dis
-{1 \over 2(1-N^2)k}
\ea
\right).$$
Computing $A''$, $B''$,  then $\widetilde u_1$ and $\widetilde w_2$ are obtained as functions of $p$ and of known data.
\par\hfill$\square$

\end{document}